\documentclass[11pt]{amsart}
\usepackage{amsmath,amscd,amssymb,amsfonts,amsthm}

 \newtheorem{theorem}{Theorem}[section]
\newtheorem{corollary}[theorem]{Corollary}
\newtheorem{proposition}[theorem]{Proposition}

\newtheorem{lemma}[theorem]{Lemma}
\theoremstyle{remark}
\newtheorem{remark}[theorem]{Remark}
\newtheorem{definition}[theorem]{Definition}
\newtheorem{remarks}[theorem]{Remarks}
\newtheorem{example}[theorem]{Example}

\newcommand\A{\mathcal{A}}

\renewcommand\L{\mathcal{L}}

\newcommand{\V}{\mathbb{V}}
\newcommand{\TM}{\mathbb{T}M}
\newcommand{\TG}{\mathbb{T}G}
\renewcommand{\O}{\mathcal{O}}
\newcommand{\Co}{\mathcal{C}}
\newcommand{\T}{\mathbb{T}}

\newcommand{\n}{\mathfrak{n}}

\newcommand{\R}{\mathbb{R}}
\newcommand{\C}{\mathbb{C}}

\newcommand{\Z}{\mathbb{Z}}

\newcommand{\pr}{\on{pr}}
\newcommand{\Cl}{{\on{Cl}}}
\newcommand{\sfe}{{\sf{e}}}
\newcommand{\sff}{{\sf{f}}}
\newcommand{\sfs}{{\sf{s}}}
\newcommand{\sfp}{{\sf{p}}}
\newcommand{\sfN}{{\sf{N}}}
\newcommand\lie[1]{\mathfrak{#1}}
\renewcommand{\k}{\lie{k}}

\newcommand{\s}{\lie{s}}
\newcommand{\g}{\lie{g}}

\renewcommand{\t}{\lie{t}}

\newcommand{\on}{\operatorname}

\newcommand{\Ad}{ \on{Ad} }
\newcommand{\ad}{ \on{ad} }

\newcommand{\Eul}{ \on{Eul} }
\newcommand{\End}{ \on{End} }

\newcommand{\Hom}{ \on{Hom}} 
\renewcommand{\ker}{ \on{ker}}

\newcommand{\SO}{ \on{SO}}
\newcommand{\GL}{ \on{GL}}

\newcommand{\Mult}{ \on{Mult}}

\newcommand\dirac{/\kern-1.2ex\partial} 
\newcommand\qu{/\kern-.7ex/} 

\renewcommand\a{\mathfrak{a}}
\newcommand{\lra}{\longrightarrow}
\newcommand{\hra}{\hookrightarrow}

\renewcommand{\d}{{\mbox{d}}}
\newcommand{\dd}{\mf{d}}

\newcommand{\ol}{\overline}

\newcommand\Sig{\Sigma}
\newcommand\sig{\sigma}
\newcommand\eps{\epsilon}
\newcommand\Om{\Omega}
\newcommand\om{\omega}

\newcommand{\f}{\frac}

\newcommand{\p}{\partial}
\renewcommand{\l}{\langle}
\renewcommand{\r}{\rangle}
\newcommand{\hh}{{\textstyle \f{1}{2}}}
\newcommand{\ti}{\tilde}
\newcommand{\wti}{\widetilde}

\newcommand\pt{\on{pt}}

\newcommand\Sm{\mathsf{S}} 

\newcommand\beqn{\begin{equation}}
\newcommand\eeqn{\end{equation}}
\newcommand{\ca}{\mathcal}
\newcommand{\wh}{\widehat}
\newcommand{\wt}{\widetilde}
\newcommand{\mf}{\mathfrak}
\newcommand{\beq}{\begin{eqnarray*}}
\newcommand{\eeq}{\end{eqnarray*}}

\newcommand{\tpi}{{2\pi\sqrt{-1}}}

\newcommand{\Cour}[1]      {[\![#1]\!]}
\newcommand{\ts}{\textstyle}
\newcommand{\tri}{\Xi}

\begin{document}

\title[]{Pure Spinors on Lie groups}

\author{A. Alekseev}
\address{University of Geneva, Section of Mathematics,
2-4 rue du Li\`evre, 1211 Gen\`eve 24, Switzerland}
\email{alekseev@math.unige.ch}

\author{H. Bursztyn}
\address{Instituto de Matem\'atica Pura e Aplicada,
Estrada Dona Castorina 110, Rio de Janeiro, 22460-320, Brasil }
\email{henrique@impa.br}

\author{E. Meinrenken}
\address{University of Toronto, Department of Mathematics,
40 St George Street, Toronto, Ontario M4S2E4, Canada }
\email{mein@math.toronto.edu}

\date{\today}
\begin{abstract}
  For any manifold $M$, the direct sum $\TM=TM\oplus T^*M$ carries a
  natural inner product given by the pairing of vectors and
  covectors. Differential forms on $M$ may be viewed as spinors for the
  corresponding Clifford bundle, and in particular there is a notion
  of \emph{pure spinor}.
  In this paper, we study pure spinors and Dirac structures
  in the case when $M=G$ is a Lie group with a bi-invariant
  pseudo-Riemannian metric, e.g. $G$ semi-simple. The
  applications of our theory include the construction of distinguished
  volume forms on conjugacy classes in $G$, and a new approach to the
  theory of quasi-Hamiltonian $G$-spaces.
\end{abstract}
\maketitle
\setcounter{tocdepth}{3}

\begin{quote}
{\it \small Dedicated to Jean-Michel Bismut on the occasion of his
60th birthday.}
\end{quote}

\setcounter{section}{-1}
\vskip.3in\section{Introduction}\label{sec:intro}
For any manifold $M$, the direct sum $\TM=TM\oplus T^*M$ carries a
non-degenerate symmetric bilinear form, extending the pairing between
vectors and covectors. There is a natural Clifford action $\varrho$ of
the sections $\Gamma(\TM)$ on the space $\Om(M)=\Gamma(\wedge T^*M)$
of differential forms, where vector fields act by contraction and
1-forms by exterior multiplication.  That is, $\wedge T^*M$ is viewed
as a spinor module over the Clifford bundle $\Cl(\TM)$.  A form
$\phi\in \Om(M)$ is called a \emph{pure spinor} if the solutions $w\in
\Gamma(\TM)$ of $\varrho(w)\phi=0$ span a Lagrangian subbundle
$E\subset \TM$. Given a closed 3-form $\eta\in\Om^3(M)$, a pure spinor
$\phi$ is called \emph{integrable} (relative to $\eta$)
\cite{al:der,gua:ge} if there exists a section $w\in \Gamma(\TM)$ with
\[ (\d+\eta)\phi=\varrho(w)\phi.\]
In this case, there is a generalized foliation of $M$ with tangent
distribution the projection of $E$ to $TM$. The subbundle $E$ defines a
\emph{Dirac structure} \cite{cou:di,sev:poi} on $M$, and the triple
$(M,E,\eta)$ is called a \textit{Dirac manifold}.

The present paper is devoted to the study of Dirac structures and
pure spinors on Lie groups $G$. We assume that the Lie algebra
$\g$ carries a non-degenerate invariant symmetric bilinear form $B$,
and take $\eta\in\Om^3(G)$ as the corresponding Cartan 3-form.
Let $\ol{\g}$ denote the Lie algebra $\g$ with the opposite
bilinear form $-B$. We will describe a trivialization
\[ \TG\cong G\times (\g\oplus \ol{\g}),\]
under which any Lagrangian Lie subalgebra $\mf{s}\subset
\g\oplus{\ol{\g}}$ defines a Dirac structure on $G$. There is also a
similar identification of spinor bundles
\[ \ca{R}\colon G\times \Cl(\g) \stackrel{\cong}{\longrightarrow} \wedge T^*G,\]
taking the standard Clifford action of $\g\oplus \ol{\g}$ on
$\Cl(\g)$, where the first summand acts by left (Clifford)
multiplication and the second summand by right multiplication, to
the Clifford action $\varrho$. This isomorphism takes the
\emph{Clifford differential} $\d_\Cl$ on $\Cl(\g)$, given as Clifford commutator by a cubic
element \cite{al:no,ko:sy2}, to the the differential $\d+\eta$ on
$\Om(G)$. As a result, pure spinors $x\in \Cl(\g)$ for the
Clifford action of $\Cl(\g\oplus \ol{\g})$ on $\Cl(\g)$ define
pure spinors $\phi=\ca{R}(x)\in \Om(G)$, and the integrability
condition for $\phi$ is equivalent to a similar condition for $x$.
The simplest example $x=1$ defines the \emph{Cartan-Dirac
structure} $E_G$ \cite{bur:di,sev:poi}, introduced by Alekseev,
\v{S}evera and Strobl in the 1990's. In this case, the resulting
foliation of $G$ is just the foliation by conjugacy classes. We
will study this Dirac structure in detail, and examine in
particular its behavior under group multiplication and under the
exponential map. When $G$ is a complex semi-simple Lie group, it
carries another interesting Dirac structure, which we call the
\emph{Gauss-Dirac structure}. The corresponding foliation of $G$
has a dense open leaf which is the `big cell' from the Gauss
decomposition of $G$.

The main application of our study of pure spinors is to the theory
of q-Hamiltonian actions \cite{al:qu,al:mom}. The original
definition of a q-Hamiltonian $G$-space in \cite{al:mom} involves
a $G$-manifold $M$ together with an invariant 2-form $\omega$ and
a $G$-equivariant map $\Phi\colon M\to G$ satisfying appropriate
axioms. As observed in \cite{bur:di,bur:int}, this definition is
equivalent to saying that the `$G$-valued moment map' $\Phi$ is a
suitable morphism of Dirac manifolds (in analogy with classical
moment maps, which are morphisms $M\to \mathfrak{g}^*$ of Poisson
manifolds).
In this paper, we will carry this observation further, and develop all
the basic results of q-Hamiltonian geometry from this perspective. A
conceptual advantage of this alternate viewpoint is that, while the
arguments in \cite{al:mom} required $G$ to be compact, the Dirac
geometry approach needs no such assumption, and in fact works in the
complex (holomorphic) category as well. This is relevant for
applications: For instance, the symplectic form on a representation
variety $\on{Hom}(\pi_1(\Sig),G)/G$ (for $\Sig$ a closed surface) can
be obtained by q-Hamiltonian reduction, and there are many interesting
examples for noncompact $G$. (For instance, the case
$G=\on{PSL}(2,\R)$ gives the symplectic form on Teichm\"uller space.)
Complex q-Hamiltonian spaces appear e.g. in the work of Boalch
\cite{boa:qu} and Van den Bergh \cite{vdb:dou}.

The organization of the paper is as follows. Sections
\ref{sec:prelim} and \ref{sec:man} contain a review of Dirac
geometry, first on vector spaces and then on manifolds. The main
new results in these sections concern the geometry of Lagrangian
splittings $\TM=E\oplus F$ of the bundle $\TM$. If $\phi,\psi \in
\Omega(M)$ are pure spinors defining $E,F$, then, as shown in
\cite{car:spi1,ch:al1}, the top degree part of $\phi^\top\wedge
\psi$ (where $\top$ denotes the standard anti-involution of the
exterior algebra) is nonvanishing, and hence defines a volume form
$\mu$ on $M$. Furthermore, there is a bivector field
$\pi\in\mf{X}^2(M)$ naturally associated with the splitting,
which satisfies
\[  \phi^\top\wedge \psi=e^{-\iota(\pi)}\mu.\]
We will discuss the properties of $\mu$ and $\pi$ in detail,
including their behavior under Dirac morphisms.

In Section \ref{sec:dirgroup} we specialize to the case $M=G$,
where $G$ carries a bi-invariant pseudo-Riemannian metric,
and our main results concern the isomorphism $\TG\cong
G\times(\g\oplus \ol{\g})$ and its properties. Under this
identification, the Cartan-Dirac structure $E_G\subset \TG$
corresponds to the diagonal $\g_\Delta\subset \g\oplus \ol{\g}$, and
hence it has a natural Lagrangian complement $F_G\subset \TG$ defined
by the anti-diagonal. We will show that the exponential map gives rise
to a Dirac morphism $(\g,E_\g,0)\to (G,E_G,\eta)$ (where $E_\g$ is the
graph of the linear Poisson structure on $\g\cong \g^*$), but this
morphism does not relate the obvious complements $F_\g=T\g$ and $F_G$.
The discrepancy is given by a `twist', which is a solution of the
\emph{classical dynamical Yang-Baxter equation}.  For $G$ complex
semi-simple, we will construct another Lagrangian complement of $E_G$,
denoted by $\wh{F}_G$, which (unlike $F_G$) is itself a Dirac
structure. The bivector field corresponding to the splitting
$E_G\oplus \wh{F}_G$ is then a Poisson structure on $G$, which
appeared earlier in the work of Semenov-Tian-Shansky \cite{se:dr}.

In Section \ref{sec:spingroup}, we construct an isomorphism $\wedge
T^*G\cong G\times \Cl(\g)$ of spinor modules, valid under a mild
topological assumption on $G$ (which is automatic if $G$ is simply
connected). This allows us to represent the Lagrangian subbundles
$E_G$, $F_G$ and $\wh{F}_G$ by explicit pure spinors $\phi_G$,
$\psi_G$, and $\wh{\psi}_G$, and to derive the differential equations
controlling their integrability. We show in particular that the
Cartan-Dirac spinor satisfies
\[ (\d+\eta)\phi_G=0.\]

Section \ref{sec:qham} investigates the foundational
properties of q-Hamiltonian $G$-spaces from the Dirac geometry
perspective. Our results on the Cartan-Dirac structure give a
direct construction of the fusion product of q-Hamiltonian spaces. On
the other hand, we use the bilinear pairing of spinors to show
that, for a q-Hamiltonian space $(M,\om,\Phi)$, the top degree
part of $e^\om\Phi^*\psi_G \in \Omega(M)$ defines a volume form
$\mu_M$. This volume form was discussed in \cite{al:du} when $G$
is compact, but the discussion here applies equally well to
non-compact or complex Lie groups. Since conjugacy classes in $G$
are examples of q-Hamiltonian $G$-spaces, we conclude that for any
simply connected Lie group $G$ with bi-invariant pseudo-Riemannian
metric (e.g. $G$ semi-simple), \emph{any conjugacy class in $G$
carries a distinguished invariant volume form}. If $G$ is complex
semi-simple, one obtains the same volume form $\mu_M$ if one
replaces $\psi_G$ with the Gauss-Dirac spinor $\wh{\psi}_G$.
However, the form $e^\om \Phi^*\wh{\psi}_G$ satisfies a nicer
differential equation, which allows us to compute the volume of
$M$, and more generally the measure $\Phi_*|\mu_M|$, by
Berline-Vergne localization \cite{be:ze}. We also explain in this Section how
to view the more general q-Hamiltonian q-Poisson spaces
\cite{al:qu} in our framework.

Lastly, in Section \ref{sec:kstar}, we revisit the theory of
$K^*$-valued moment maps in the sense of Lu \cite{lu:mo} and its
connections with $P$-valued moment maps \cite[Sec.~10]{al:mom}
from the Dirac geometric standpoint.

\vskip.2in \noindent {\bf Acknowledgments.} We would like to thank
Marco Gualtieri and Shlomo Sternberg for helpful comments on an
earlier version of this paper.  The research of A.A. was supported in
part by the Swiss National Science Foundation, while E.M. was
supported by an NSERC Discovery Grant and a Steacie Fellowship. H.B.
acknowledges financial support from CNPq, and thanks University of
Toronto, the Fields Institute and Utrecht University for hospitality
at various stages of this project.  We also thank the Erwin
Schr\"odinger Institut and the Forschungsinstitut Oberwolfach where
part of this work was carried out.

\vskip.2in \noindent {\bf Notation.} Our conventions for Lie group
actions are as follows: Let $G$ be a Lie group (not necessarily
connected), and $\g$ its Lie algebra.  A $G$-action on a manifold
$M$ is a group homomorphism $\A\colon G\to \on{Diff}(M)$ for which
the action map $G\times M\to M,\ (g,m)\mapsto \A(g)(m)$ is smooth.
Similarly, a $\g$-action on $M$ is a Lie algebra homomorphism
$\A\colon \g\to \mf{X}(M)$ for which the map $\g\times M\to TM,\
(\xi,m)\mapsto \A(\xi)_m$ is smooth. Given a $G$-action $\A$, one
obtains a $\g$-action by the formula $\A(\xi)(f)=\f{\p}{\p
t}\big|_{t=0}\A(\exp(-t\xi))^*f$, for $f\in C^\infty(M)$ (here
vector fields are viewed as derivations of the algebra of smooth
functions).
\newpage

{\small \tableofcontents \pagestyle{headings}}

\vskip.3in\section{Linear Dirac geometry} \label{sec:prelim}
The theory of Dirac manifolds was initiated by Courant and
Weinstein in \cite{cou:di,couwein:beyond}. We briefly review this
theory, developing and expanding the approach via pure spinors
advocated by Gualtieri \cite{gua:ge} (see also Hitchin
\cite{hi:gen} and Alekseev-Xu \cite{al:der}). All vector spaces in
this section are over the ground field $\mathbb{K}=\R$ or $\C$. We
begin with some background material on Clifford algebras and
spinors (see e.g. \cite{ch:al1} or \cite{me:lec}.)

\subsection{Clifford algebras}\label{subsec:clif}
Suppose $V$ is a vector space with a non-degenerate symmetric
bilinear form $B$. We will sometimes refer to such a bilinear form
$B$ as an \textit{inner product} on $V$. The \emph{Clifford
algebra} over $V$ is the associative unital algebra
generated by the elements of $V$, with relations
\[ vv'+v'v=B(v,v')\,1. \]
It carries a compatible $\Z_2$-grading and $\Z$-filtration, such
that the generators $v\in V$ are odd and have filtration degree
$1$. We will denote by $x\mapsto x^\top$ the canonical
anti-automorphism of exterior and Clifford algebras, equal to the identity
on $V$. For any $x\in \Cl(V)$, we denote by $l^\Cl(x),r^\Cl(x)$
the operators of graded left and right multiplication on $\Cl(V)$:
\[ l^\Cl(x)x'=xx',\ \ \ r^\Cl(x)x'=(-1)^{|x||x'|}x'x.\]
Thus $l^\Cl(x)-r^\Cl(x)$ is the operator of graded commutator
$[x,\cdot]_\Cl$.

The \emph{quantization map} $q\colon \wedge V\to
\Cl(V)$ is the isomorphism of vector spaces defined by
$q(v_1\wedge\cdots\,\wedge v_r)=v_1\cdots v_r$ for pairwise orthogonal
elements $v_i\in V$.  Let
\[ \on{str}\colon \Cl(V)\to \det(V):=\wedge^{\mathrm{top}}(V)\]
be the \emph{super-trace}, given by $q^{-1}$, followed by taking the top degree part.
It has the property $\on{str}([x,x']_\Cl)=0$.

A \emph{Clifford module} is a vector space $\Sm$ together with an algebra
homomorphism $\varrho\colon \Cl(V)\to \End(\Sm)$. If $\Sm$ is a
Clifford module, one has a dual Clifford module given by the dual
space $\Sm^*$ with Clifford action
$\varrho^*(x)=\varrho(x^\top)^*$.

Recall that $\on{Pin}(V)$ is the subgroup of $\Cl(V)^\times$ generated
by all $v\in V$ whose square in the Clifford algebra is $vv=\pm 1$. It
is a double cover of the orthogonal group $\on{O}(V)$, where $g\in
\on{Pin}(V)$ takes $v\in V$ to $(-1)^{|g|}gvg^{-1}$, using Clifford
multiplication. The \emph{norm homomorphism} for the Pin group is the
group homomorphism
\begin{equation}\label{eq:normhom}
\sfN\colon \on{Pin}(V)\to \{-1,+1\},\ \ \
\sfN(g)=g^\top g=\pm 1. \end{equation}

Let $\{\cdot,\cdot\}$ be the graded Poisson bracket on $\wedge V$,
given on generators by $\{v_1,v_2\}=B( v_1,v_2)$. Then $\wedge^2
V$ is a Lie algebra under the Poisson bracket, isomorphic to
$\mf{o}(V)$ in such a way that $\varepsilon \in\wedge^2 V$
corresponds to the linear map $v\mapsto \{\varepsilon,v\}$.
The Lie algebra $\mf{pin}(V)\cong \mf{o}(V)$ is realized as the Lie subalgebra
$q(\wedge^2(V))\subset \Cl(V)$.

A subspace $E\subset V$ is called \emph{isotropic} if $E\subset
E^\perp$ and \emph{Lagrangian} if $E=E^\perp$. The set of Lagrangian
subspaces is non-empty if and only if the bilinear form is
\emph{split}.  If $\mathbb{K}=\C$, this just means that $\dim V$ is
even, while for $\mathbb{K}=\R$ this requires that the bilinear form
has signature $(n,n)$.  From now on, we will reserve the letter $W$
for a vector space with split bilinear form $\l\cdot,\cdot\r$. We
denote by $\on{Lag}(W)$ the Grassmann manifold of Lagrangian subspaces
of $W$. It carries a transitive action of the orthogonal group
$\on{O}(W)$.

\begin{remark}
Suppose $\mathbb{K}=\R$, and identify $W\cong \R^{2n}$ with the
standard bilinear form of signature $(n,n)$. The group $\on{O}(W)\cong
\on{O}(n,n)$ has maximal compact subgroup $\on{O}(n)\times
\on{O}(n)$. Already the subgroup $\on{O}(n)\times \{1\}$ acts
transitively on $\on{Lag}(W)$, and in fact the action is free. It
follows that $\on{Lag}(W)$ is diffeomorphic to $\on{O}(n)$. Further
details may be found in \cite{me:cli}.
\end{remark}

\subsection{Pure spinors}\label{subsec:purespinor}
An irreducible module $\Sm$ over the Clifford algebra $\Cl(W)$ is
called a \emph{spinor module}.
Any $E\in\on{Lag}(W)$ defines a spinor module $\Sm=\Cl(W)/\Cl(W)E$.
The choice of a Lagrangian complement $F$ to $E$ identifies $\Sm=\wedge E^*$,
where the generators in $E\subset W$ act by contraction and the
generators in $F\subset W$ act by exterior multiplication. (Here
$F$ is identified with $E^*$, using the pairing defined by
$\l\cdot,\cdot\r$.) The dual spinor module is $\Sm^*=\wedge E$, with generators
in $E$ acting by exterior multiplication and those in $F$ by
contraction.

For any non-zero element $\phi\in \Sm$ of a spinor module, its
null space
\[N_\phi=\{w\in W|\,\varrho(w)\phi=0\}\]
is easily seen
to be isotropic.  The element $\phi\in \Sm$ is a \emph{pure
spinor} \cite{car:spi1} provided $N_\phi$ is Lagrangian.  One can
show that any Lagrangian subspace $E\in\on{Lag}(W)$ arises in this
way: in fact, $\Sm^E=\{\phi\in \Sm|\,\varrho(E)\phi=0\}$ is a
one-dimensional subspace, with non-zero elements given by the pure
spinors defining $E$.  Any spinor module $\Sm$ admits a
$\Z_2$-grading (unique up to parity inversion) compatible with the
Clifford action.  Pure spinors always have a definite parity,
either even or odd.

\begin{example} \label{ex:example1}
  Let $V$ be a vector space with inner product $B$.  We denote by
$\ol{V}$ the same vector space with the opposite bilinear form
$-B$. Then $W=V\oplus \ol{V}$ is a vector space with split
bilinear form. The space $\Sm=\Cl(V)$ is a spinor module over
$\Cl(W)=\Cl(V)\otimes\Cl(\ol{V})$, with Clifford action given on
generators by $\varrho(v\oplus v')= l^\Cl(v)-r^\Cl(v')$. The
element $1\in \Cl(V)$ is a pure spinor, with corresponding
Lagrangian subspace the diagonal $V_\Delta\subset V\oplus \ol{V}$.
\end{example}

\subsection{The bilinear pairing of spinors}
 For any two spinor
modules $\Sm_1,\Sm_2$ over $\Cl(W)$, the space
$\Hom_{\Cl(W)}(\Sm_1,\Sm_2)$ of intertwining operators is
one-dimensional. Given a spinor module $\Sm$, let
\[ K_{\Sm}=\Hom_{\Cl(W)}(\Sm^*,\Sm)\]
be the \emph{canonical line}. There is a bilinear pairing
\cite{car:spi1}
\[ \Sm\otimes \Sm \to K_{\Sm},\ \phi\otimes \psi\mapsto
(\phi,\psi)_{\scriptscriptstyle{\Sm}},\]
defined by the isomorphism $\Sm\otimes \Sm\cong \Sm\otimes \Sm^*\otimes
\Hom_{\Cl(W)}(\Sm^*,\Sm)$ followed by the duality pairing
$\Sm\otimes\Sm^*\to \mathbb{K}$. The pairing satisfies
\begin{equation}\label{eq:chevprop}
 (\varrho(x^\top)\phi,\psi)_{\scriptscriptstyle{\Sm}}=
 (\phi,\varrho(x)\psi)_{\scriptscriptstyle{\Sm}},\ \
 x\in\Cl(W),\end{equation}
and is characterized by this property up to a scalar.
\eqref{eq:chevprop} implies the following  invariance property under the action of the group
$\on{Pin}(V)$, involving the norm homomorphism \eqref{eq:normhom},
\[ (g\phi,g\psi)_{\scriptscriptstyle{\Sm}}
={\sfN(g)} (\phi,\psi)_{\scriptscriptstyle{\Sm}},\ \
g\in\on{Pin}(V).
\]
\begin{theorem}[E. Cartan \cite{car:spi1}] \label{th:chevalley}
  Let $\Sm$ be a spinor modules over $\Cl(W)$, and let $\phi,\psi\in
  \Sm$ be pure spinors. Then the corresponding Lagrangian subspaces
  $N_\phi,N_\psi$ are transverse if and only if
  $(\phi,\psi)_{\scriptscriptstyle{\Sm}}\not=0$.
\end{theorem}
A simple proof of this result is given in Chevalley's book
\cite[III.2.4]{ch:al1}, see also \cite[Section 3.5]{me:lec}.
\begin{example}\label{ex:example2}
  Suppose $V$ is a space with inner product $B$, and take $\Sm=\Cl(V)$
  as a spinor module over $\Cl(V\oplus \ol{V})$ (cf.  Example
  \ref{ex:example1}). Then $K_{\Sm}=\det(V)$, with bilinear pairing on
  spinors given as
\begin{equation}\label{eq:clifpair}
 (x,x')_{\Cl(V)}=\on{str}(x^\top x')\in\det(V).\end{equation}
Using the isomorphism $q: \wedge(V)\to \Cl(V)$ to identify
$\Sm\cong\wedge(V)$, the bilinear pairing becomes
\begin{equation}\label{eq:clifpair2}
 (y,y')_{\wedge(V)}=(y^\top \wedge y')^{\on{[top]}}\in\det(V).
\end{equation}
\end{example}
%

\subsection{Contravariant spinors}\label{subsec:contra}
For any vector space $V$, the direct sum $\V:=V\oplus V^*$ carries
a split bilinear form given by the pairing between $V$ and $V^*$:
\begin{equation}\label{eq:bilin}
\l w_1,\, w_2 \r=\l\alpha_1,v_2\r+\l\alpha_2,v_1\r,\ \  \ \ w_i=v_i\oplus\alpha_i\in\V.
\end{equation}
Every vector space $W$ with split bilinear form is of this form,
by choosing a pair of transverse Lagrangian subspaces $V,V'$, and
using the bilinear form to identify $V'=V^*$. Then $\Sm=\wedge
V^*$, with Clifford action given on generators $w=v\oplus\alpha\in\V$ by
\[ \varrho(w)=\eps(\alpha)+\iota(v)\]
(where $\eps(\alpha)=\alpha\wedge \cdot$),
is a natural choice of spinor module for $\Cl(\V)$. The
restriction of $\varrho$ to $\wedge V^*\subset \Cl(\V)$ is given
by exterior multiplication, while the restriction to $\wedge
V\subset \Cl(\V)$ is given by contraction \footnote{We are using
the convention that
  $\iota\colon \wedge(V)\to \on{End}(\wedge V^*)$ is the extension of
  the map $v\mapsto \iota(v)$ as an \emph{algebra homomorphism}. Note
  that some authors use the extension as an algebra
  anti-homomorphism.}.  The line $K_{\Sm}=\on{Hom}_{\Cl(\V)}(\Sm^*,\Sm)$ is
canonically isomorphic to $\det(V^*)=\wedge^{\mathrm{top}}V^*$,
and the bilinear pairing on spinors is simply
\[ (\phi,\psi)_{\wedge(V^*)}=(\phi^\top\wedge \psi)^{[\mathrm{top}]}\in \det(V^*),\]
similar to Example \ref{ex:example2}. Theorem \ref{th:chevalley}
shows that if $\phi,\psi$ are pure spinors for transverse
Lagrangian subspaces, the pairing $(\phi,\psi)_{\wedge(V^*)}$ defines a
\emph{volume form} on $V$.
\begin{remarks}\label{rem:reference}
We mention the following two facts for later reference.
\begin{enumerate}
\item\label{it:aaa}
 We have the identity
\[
(-1)^{|\phi|}(-(\varrho(w)\phi)^\top \wedge\psi +\phi^\top \wedge
(\varrho(w)\psi))=\iota(v)(\phi^\top \wedge \psi)
,\ \ w=v\oplus \alpha\in \V,
\]
which refines property \eqref{eq:chevprop} of the bilinear
pairing. \item \label{it:rem-b} One can also consider the
\emph{covariant spinor module} $\wedge(V)$, obtained by reversing
the roles of $V$ and $V^*$.  Suppose $\mu\in\det(V)$ is non-zero,
and let $\star\colon\wedge(V^*)\to\wedge(V)$ be the corresponding
star operator, defined by $\star \phi=\iota(\phi)\mu$.  Let
$\mu^*$ be the dual generator defined by $\star((\mu^*)^\top)=1$ .
Then $\star$ is an isomorphism of $\Cl(\V)$-modules. Furthermore,
using $\mu,\mu^*$ to trivialize $\det(V),\det(V^*)$, the
isomorphism intertwines the bilinear pairings:
\[
(\phi,\psi)_{\wedge(V^*)} = (\star\phi,\star\psi)_{\wedge(V)},\ \ \phi,\psi\in\wedge(V^*).
\]
\end{enumerate}
\end{remarks}

Any 2-form $\omega\in \wedge^2 V^*$ defines a pure spinor
$\phi=e^{-\omega}$, with $N_\phi$ the graph of $\omega$:
\[ \on{Gr}_\omega=\{v\oplus \alpha|\ v\in V,\ \alpha=\iota(v)\omega\}.\]
Note that, in accordance with Theorem \ref{th:chevalley},
$\on{Gr}_\omega\cap V=\{0\}$ if and only if $\omega$ is
non-degenerate, if and only if $(e^\omega)^{[\mathrm{top}]}$ is
non-zero. The most general pure spinor $\phi\in\wedge V^*$ can be
written in the form
\begin{equation}\label{eq:general}
 \phi=e^{-\omega_Q}\wedge\theta,
\end{equation}
where $\omega_Q\in \wedge^2 Q^*$ is a 2-form on a subspace
$Q\subset V$ and $\theta\in\det(\on{Ann}(Q))\backslash \{0\}$ is a
volume form on $V/Q$. To write \eqref{eq:general}, we have chosen
an extension of $\omega_Q$ to a 2-form on $V$. (Clearly, $\phi$
does not depend on this choice.)  The corresponding Lagrangian
subspace is
\[ N_\phi=\{v\oplus\alpha|\ v\in Q,\ \alpha|_Q=\iota(v)\omega_Q \}.\]
The triple $(Q,\omega_Q,\theta)$ is uniquely determined by $\phi$, see
e.g.  \cite[III.1.9]{ch:al1}. A simple consequence is that any pure
spinor has definite parity, that is, $\phi$ is either even or odd
depending on the parity of $\dim(V/Q)$.  For any $E\in\on{Lag}(\V)$ we
define subspaces $\ker(E)\subset \on{ran}(E)\subset V$ by
\[ \ker(E)=E\cap V,\ \ \ \ \on{ran}(E)=\on{pr_V}(E),\]
where $\pr_V\colon \V\to V$ is the projection along $V^*$. For any
pure spinor $\phi$, written in the form \eqref{eq:general}, we have
$\on{ran}(E_\phi)=Q$ and $\ker(E_\phi)=\ker(\omega_Q)$. In particular,
$\phi^{[\mathrm{top}]}$ is non-zero if and only if $\ker(E_\phi)=0$.
Similarly, $\on{ran}(E_\phi)=V$ if and only if $\phi^{[0]}$ is
non-zero, if and only if $\phi=e^{-\omega}$ for a global 2-form $\omega$.

\subsection{Action of the orthogonal group}\label{subsec:orthaction}
%

Recall the identification $\wedge^2(W)\cong \mf{o}(W)$ (see
Section \ref{subsec:clif}).  For any Lagrangian subspace $E\subset
W$, the space $\wedge^2(E)$ is embedded as an Abelian subalgebra
of $\mf{o}(W)$. The inclusion map exponentiates to an injective
group homomorphism,
\begin{equation}
 \wedge^2(E)\to \on{O}(W),\ \ \varepsilon \mapsto A^\varepsilon,\ \ \ \
A^\varepsilon(v\oplus \alpha)=v\oplus (\alpha-\iota(v)\varepsilon),
\end{equation}
with image the orthogonal transformations fixing $E$ pointwise.
The subgroup $\wedge^2(E)$ acts freely and transitively on the
subset of $\on{Lag}(W)$ of Lagrangian subspaces transverse to $E$,
which therefore becomes an affine space. Observe that
$A^\varepsilon$ has a distinguished lift
$\wt{A}^\varepsilon=\exp(\varepsilon) \in\on{Pin}(W)$ (exponential
in the subalgebra $\wedge(E)\subset \Cl(W)$).

For any spinor module
$\Sm$ over $\Cl(W)$, the induced representation of the group
$\on{Pin}(W)\subset \Cl(W)^\times$ preserves the set of pure
spinors, and the map $\phi\mapsto N_\phi$ is equivariant. That is,
if $\wt{A}\in\on{Pin}(W)$ lifts $A\in \on{O}(W)$, then
\[ N_{\varrho(\wt{A})\phi}=A(N_\phi).\]
Consider again the case $W=\V$. Then 2-forms $\omega\in \wedge^2
V^*$ and bivectors $\pi\in \wedge^2(V)$ define orthogonal
transformations
\[ A^{-\omega}(v\oplus\alpha)= v\oplus(\alpha+\iota_v\om),\ \
A^{-\pi}(v\oplus\alpha)=(v+\iota_\alpha\pi)\oplus \alpha.
\]
Their lifts act in the spin representation as follows:
\begin{equation}
 \varrho(\wt{A}^{-\omega})\phi=e^{-\omega}\,\phi,\ \
 \varrho(\wt{A}^{-\pi})\phi=e^{-\iota(\pi)}\phi.
\end{equation}

\subsection{Morphisms}
It is easy to see that the
group of orthogonal transformations of $\V$ preserving the `polarization'
\begin{equation}\label{eq:exact}
 0\lra V^*\lra \V\lra V\lra 0\end{equation}
(i.e., taking the subspace $V^*$ to itself) is the semi-direct product
$\wedge^2 V^*\rtimes \GL(V)\subset \on{O}(\V)$, where
$\omega\in\wedge^2V^*$ acts as $A^{-\omega}$ and $\GL(V)$ acts in the
natural way on $V$ and by the conjugate transpose on $V^*$.

More generally, for vector spaces $V$ and $V'$, we define the set
of \emph{morphisms from $\V$ to $\V'$} \cite{hu:ext} to be
\[ \on{Hom}(V,V')\times \wedge^2 V^*, \]
with the following composition law:
\begin{equation}\label{eq:compo}
 (\Phi_1,\omega_1)\circ   (\Phi_2,\omega_2)=(\Phi_1\circ \Phi_2,\ \om_2+\Phi_2^*\om_1).
\end{equation}
Given $w= v\oplus\alpha\in \V$ and $w'=v'\oplus\alpha'\in \V'$, we
write
\[ w\sim_{(\Phi,\om)} w'\ \ \Leftrightarrow \ \ v'=\Phi(v),\ \Phi^*\alpha'=\alpha+\iota_v\om.
\]
In particular, taking $V'=V$ and $\Phi=\on{id}$ we have
$w\sim_{(\on{id},\omega)}w'$ if and only $w'=A^{-\om}(w)$.
The \emph{graph of a morphism} $(\Phi,\om)$ is the subspace
\begin{equation}\label{eq:graph}
\Gamma_{(\Phi,\om)}=\{(w',w)\in \V'\times {\V}\,|\ w\sim_{(\Phi,\om)} w'\}.
 \end{equation}
We have $\Gamma_{(\Phi_1,\om_1)\circ
  (\Phi_2,\om_2)}=\Gamma_{(\Phi_1,\om_1)}\circ
\Gamma_{(\Phi_2,\om_2)}$ under composition of relations.
The morphisms $(\Phi,\om)$ are `isometric', in the sense
that
\begin{equation}\label{eq:isometric}
w_1\sim_{(\Phi,\om)} w_1',\ \ w_2\sim_{(\Phi,\om)} w_2'
\Rightarrow \l w_1,w_2\r=\l w_1',w_2'\r.\end{equation}
Equivalently, $\Gamma_{(\Phi,\om)}$ is Lagrangian in $\V' \oplus \ol{\V}$.
We write
\[\begin{split}
\ker(\Phi,\om)&=\{w\in \V|\, w\sim_{(\Phi,\om)} 0\},\ \ \
\\\on{ran}(\Phi,\om)&=\{w'\in \V'|\, \exists w\in \V\colon w\sim_{(\Phi,\om)} w'\}.\end{split}\]
Thus $\ker(\Phi,\om)=\{(v,-\iota_v\om)|\, v\in \ker(\Phi)\}$
while $\on{ran}(\Phi,\om)=\on{ran}(\Phi)\oplus (V')^*$.
\begin{definition}
Let $(\Phi,\om)\colon \V\to \V'$ be a morphism, and $E\in
\on{Lag}(\V)$. We define the \emph{forward image} $E'\in
\on{Lag}(\V')$ to be the Lagrangian subspace
\[  E':=\Gamma_{(\Phi,\om)}\circ E = \{ w'\in \V'\,|\,
\exists \, w\in E\colon w\sim_{(\Phi,\om)} w'  \}.
\]
Similarly, for $F'\in \on{Lag}(\V')$  the \textit{backward image}
is defined as the Lagrangian subspace
$$
F:=F'\circ \Gamma_{(\Phi,\om)} = \{w\in \V\,|\ \exists\, w'\in
F'\colon \, w\sim_{(\Phi,\om)} w'\}.
$$
\end{definition}
The proof that forward and backward images of Lagrangian subspaces
are Lagrangian is parallel to the similar statement in the
symplectic category of Guillemin-Sternberg \cite{gui:lin} (see also
Weinstein \cite{we:le}).
It is simple to check that the
composition $E'=\Gamma_{(\Phi,\om)}\circ E$ is transverse if and
only if $\ker(\Phi,\om)\cap E=\{0\}$. Similarly, the
composition $F= F'\circ \Gamma_{(\Phi,\om)}$ is transverse if and
only if $\on{ran}(\Phi,\om)+F'=\V'$ (equivalently, if and only if
$\on{ran}(\Phi)+\on{ran}(F')=V'$).
\begin{remark}
  As in the symplectic category \cite{gui:lin,we:le}, one could
  consider morphisms given by arbitrary \emph{Lagrangian relations},
  i.e. Lagrangian subspaces $\Gamma\subset \V' \oplus \ol{\V}$ (see
  e.g.  \cite{bur:gauge}). The graphs \eqref{eq:graph} of morphisms
  $(\Phi,\om)$ are exactly those Lagrangian relations preserving the
  `polarization' \eqref{eq:exact}, in the sense that $\Gamma\circ
  V^*=(V')^*$ (where the composition is transverse), see
  \cite{hu:ext}.
\end{remark}

The $(\Phi,\om)$-relation may also be interpreted in terms of the
spinor representations of $\Cl(\V)$ and $\Cl(\V')$:
\begin{lemma} Suppose $(\Phi,\omega)\colon \V\to \V'$ is a morphism,
  and $w\in\V,\ w'\in \V'$. Then
\begin{equation}\label{eq:interspin}
w\sim_{(\Phi,\om)} w'\ \  \Leftrightarrow\ \
\varrho(w)(e^{\om} \Phi^*\psi')=e^{\om} \Phi^*(\varrho(w')\psi'),\ \ \
\psi'\in\wedge(V')^*.
\end{equation}
\end{lemma}
\begin{proof} This follows from
$(\eps(\alpha)+\iota_v)(e^{\om}
\Phi^*\psi')
=e^\om (\eps(\alpha+\iota_v\om)+\iota_v)\Phi^*\psi'$, for $v\oplus\alpha\in\V$.
\end{proof}

\begin{lemma}\label{lem:psipull}
Suppose $(\Phi,\om)\colon \V\to \V'$ is a morphism, and
$\psi'$ is a pure spinor defining a Lagrangian subspace $F'$. Then
$\psi=e^{\om} \Phi^*\psi'$ is non-zero if and only if the composition
$F=F'\circ \Gamma_{(\Phi,\om)}$ is
transverse, and in that case it is a pure spinor defining $F$.
\end{lemma}
\begin{proof}
Suppose $w\in F$, i.e. $w\sim_{(\Phi,\om)} w'$ with $w'\in
F'=N_{\psi'}$. Then $w\in N_\psi$ by Equation \eqref{eq:interspin}.
Thus $F\subset N_\psi$. For $\psi\not=0$, this is an equality since $F$ is Lagrangian.
\end{proof}

\begin{example}\label{ex:wedge}
Suppose $E,F\subset \V$ are Lagrangian, with defining pure spinors
$\phi,\psi$. Let $E^\top$ be the image of $E$ under the map
$v\oplus\alpha\mapsto v\oplus(-\alpha)$. Then $\phi^\top$ is a
pure spinor defining $E^\top$. Consider the diagonal inclusion
$\on{diag}\colon V\to V\times V$, so that $\on{diag}^*(\phi^\top
\otimes \psi)=\phi^\top\wedge\psi$ is just the wedge product. The
wedge product is non-zero if and only if the composition $E^\top
\wedge F:=(E^\top \times F)\circ \Gamma_{\on{diag}}$ is
transverse. This is the case, for instance, if $E$ and $F$ are
transverse (since the top degree part of $\phi^\top\wedge\psi$ is
non-zero in this case). Explicitly,
\[ E^\top \wedge F
=\{v\oplus\alpha\,|\, \exists v\oplus\alpha_1\in E,\
v\oplus\alpha_2 \in F\colon \alpha=\alpha_2-\alpha_1\}.
\]
Note that $\on{ran}(E^\top\wedge F)=\on{ran}(E)\cap \on{ran}(F)$,
with 2-form the difference of the restrictions of the 2-forms on
$\on{ran}(E)$ and $\on{ran}(F)$. Note also that $(A^{-\om}(E))^\top\wedge
(A^{-\om}(F))=E^\top\wedge F$ for all $\om\in\wedge^2V^*$.
%
\end{example}
This ``wedge product'' operation of Lagrangian subspaces was
noticed independently by Gualtieri, see \cite{gua:ge2}.

\subsection{Dirac spaces}

A \emph{Dirac space} is a pair $(V,E)$, where $V$ is a vector
space and $E\subset \V$ is a Lagrangian subspace. As remarked in
Section \ref{subsec:contra}, $E$ determines a subspace
$Q=\on{ran}(E)=\pr_V(E)\subset V$ together with a 2-form $\om_Q\in
\wedge^2Q^*$,
\begin{equation}\label{eq;2form}
 \om_Q(v,v')=\l \alpha,v'\r=-\l\alpha',v\r
 \end{equation}
for arbitrary lifts $v\oplus\alpha, v'\oplus\alpha'\in E$ of $v,v'
\in Q$. The kernel of $\om_Q$ is the subspace $\ker(E)=E\cap V$.
Conversely, any subspace $Q$ equipped with a 2-form $\om_Q$
determines a Lagrangian subspace $E=\{v\oplus\alpha \in\V|\ v\in
Q,\ \alpha|_Q=\omega_Q(v,\cdot)\}$. The gauge transformation
$A^{-\om}$ by a 2-form $\omega\in\wedge^2 V^*$ preserves $Q$,
while $\om_Q$ changes by the pull-back of $\om$ to $Q$.

\begin{definition}
Let $(\V,E)$ and $(\V',E')$ be Dirac spaces. A \emph{Dirac
morphism} $(\Phi,\om)\colon (V,E)\to (V',E')$ is a morphism
$(\Phi,\om)$ with $E'=\Gamma_{(\Phi,\om)}\circ E$. It is called a
\emph{strong Dirac morphism}\footnote{In the particular
case when $\omega=0$, Dirac morphisms are also called
\emph{forward Dirac maps} \cite{bur:int,bur:gauge}, and strong
Dirac morphisms are called \emph{Dirac realizations} \cite{bur:di}.} if this composition is transverse,
i.e.,
\[\ker(\Phi,\om)\cap E=\{0\}.\]
\end{definition}
Clearly, the composition of strong Dirac morphisms is again a
strong Dirac morphism. Note that the definition of a Dirac
morphism $(\Phi,\omega)\colon (V,E)\to (V',E')$ amounts to the
existence of a linear map $\wh{\a}\colon E'\to E$, assigning to
each $w'\in E'$ an element of $E$ to which it is
$(\Phi,\omega)$-related:
\begin{equation}\label{eq:AA}
 \wh{\a}(w')\sim_{(\Phi,\omega)} w'\ \ \ \ \ \  \forall w'\in E'.\end{equation}
The map $\wh{\mf{a}}$ is completely determined by its $V$-component
\[ \mf{a}=\pr_V\circ \,\wh{\mf{a}}\colon E'\to V,\]
since
$\wh{\mf{a}}(v'\oplus \alpha')=v\oplus (\Phi^*\alpha'+\iota_v
\omega)$
where $v=\mf{a}(v'\oplus\alpha')$. Hence $(\Phi,\omega)$ is a
Dirac morphism if and only if there exists a map $\mf{a}\colon
E'\to V$, such that the corresponding map $\wh{\mf{a}}$ takes
values in $E$.

\begin{lemma}
For a strong Dirac morphism $(\Phi,\om)\colon (V,E)\to (V',E')$,
the map $\wh{\a}$ satisfying \eqref{eq:AA} is unique. Its range is
given by
\begin{equation}\label{eq:therange}
\on{ran}(\wh{\mf{a}})=E\cap \ker(\Phi,\om)^\perp. \end{equation}
%
\end{lemma}

\begin{proof}
  The map $\wh{\a}$ associated to a Dirac morphism is unique up to
  addition of elements in $E\cap \ker(\Phi,\om)$. Hence, it is unique
  precisely if the Dirac morphism is strong. Its range consists of all
  $w\in E$ which are $(\Phi,\om)$-related to some element of $w'\in
  E'$.  By \eqref{eq:isometric}, the subspace $\{w\in \V\,|\, \exists
  w'\in\V'\colon w\sim_{(\Phi,\om)} w'\}$ is orthogonal to
  $\ker(\Phi,\omega)$. Hence, by a dimension count it coincides with
  $\ker(\Phi,\om)^\perp$. On the other hand, if $w\in E$ lies in this
  subspace,
  it is automatic that $w'\in E'$ since $E'=\Gamma_{(\Phi,\om)}\circ
  E$.
\end{proof}
\begin{example}
  Let $E\subset\V$ be a Lagrangian subspace, and let $\omega_Q$ be the
  corresponding 2-form on $Q=\on{ran}(E)$.  Let $\iota_Q\colon Q\to V$
  be the inclusion. Then $(\iota_Q,\om_Q)\colon (Q,Q)\to (V,E)$ is a
  strong Dirac morphism. Equivalently $(\iota_Q,0)\colon (Q,\on{Gr}_{\omega_Q})\to
  (V,E)$ is a strong Dirac morphism. Here $\mf{a}(v\oplus\alpha)=\iota_Q(v)$.
\end{example}

\begin{example}
  Suppose $\pi\in\wedge^2V$ and $\pi'\in\wedge^2V'$.
  Then $(\Phi,0)\colon (V,\on{Gr}_\pi)\to (V',\on{Gr}_{\pi'})$ is a
  Dirac morphism if and only if $\Phi(\pi)=\pi'$. It is automatically
  strong (since $\ker(\on{Gr}_\pi)=0$), with
  $\mf{a}(v'\oplus\alpha')=\pi^\sharp(\Phi^*\alpha')$.
\end{example}

\begin{proposition}\label{prop:transI}
Suppose $(\Phi,\omega)\colon (V,E)\to (V',E')$ is a Dirac morphism,
and that $F'$ is a Lagrangian subspace transverse to $E'$. Let
$\phi$ be a pure spinor defining $E$, and $\psi'$ a pure spinor
defining $F'$. Then $\psi:=e^\om
\Phi^*\psi'$ is non-zero, and is a pure spinor defining the backward
image $F=F'\circ \Gamma_{(\Phi,\om)}$. Moreover,
the following are equivalent:
\begin{enumerate}
\item $(\Phi,\omega)$ is a \emph{strong} Dirac morphism,
\item the
backward image $F$ is transverse to $E$,
\item
The pairing $(\phi,\psi)_{\wedge(V^*)}\in \det(V^*)$ is non-zero, that is, it is a volume form
on $V$.
\end{enumerate}
\end{proposition}
\begin{proof}
By \eqref{eq:general}, we may
write $\psi'=e^{-\omega_{Q'}}\theta'$, where $\omega_{Q'}$ is a
2-form on $Q'=\on{ran}(F')$, and $\theta'\in
\wedge^{\mathrm{top}}(V'/\on{ran}(F'))^*$. Identifying
$(V'/\on{ran}(F'))^*$ with the annihilator of $\on{ran}(F')$, this
gives
\[ \begin{split}
\psi\not=0&\Leftrightarrow
\Phi^*\theta'\not=0\\
&\Leftrightarrow \ker(\Phi^*)\cap \on{ann}(\on{ran}(F'))=0\\
&\Leftrightarrow \{w'\in F'|\ 0\sim_{(\Phi,\om)} w'\}=\{0\}.\\
\end{split}\]
(Indeed, $0\sim_{(\Phi,\om)} w'$ if and only if $w'=0\oplus
\alpha'$ with $\Phi^*\alpha'=\{0\}$. Moreover $w'\in
F'=(F')^\perp$ if and only if $\alpha'\in
\on{ann}(\on{ran}(F'))$.) But the condition $0\sim_{(\Phi,\om)}
w'$ implies that $w'\in E'$. Since $E'\cap F'=0$
it follows that $\{w'\in F'|\ 0\sim_{(\Phi,\om)} w'\}=\{0\}$, hence
$\psi\not=0$. Lemma \ref{lem:psipull} shows that
it is a pure spinor defining the backward image $F$.

$(a) \Leftrightarrow (b)$. By definition, $E\cap F$ consists of all $w\in E$ such that
$w\sim_{(\Phi,\om)}w'$ for some $w'\in F'$. Since
$E'=\Gamma_{(\Phi,\om)}\circ  E$, this element $w'$ also lies in $E'$, and
hence $w'=0$. Thus,
\[ E\cap F=E\cap \ker(\Phi,\om),\]
which is zero precisely if the Dirac morphism $(\Phi,\omega)$ is
strong. $(b) \Leftrightarrow (c)$ is immediate from
Theorem~\ref{th:chevalley}.
\end{proof}

%

\subsection{Lagrangian splittings}

Suppose $W$ is a vector space with split bilinear form.  By a
\emph{Lagrangian splitting} of $W$ we mean a direct sum
decomposition $W=E\oplus F$ into transverse Lagrangian subspaces.
\begin{lemma}
Let $W$ be a vector space with split bilinear form $\l\cdot,\cdot\r$.
There is a 1-1 correspondence between projection operators
$\sfp\in \End(W)$ with the property $\sfp+\sfp^t=1$, and
Lagrangian splittings $W=E\oplus F$. (Here $\sfp^t$ is the
transpose with respect to the inner product on $W$.)
\end{lemma}
\begin{proof}
A Lagrangian splitting of $W$ into transverse Lagrangian subspaces
is equivalent to a projection operator whose kernel and range are
isotropic. For any projection operator $\sfp=\sfp^2$, the
range $\on{ran}(\sfp)$ is isotropic if and only if
$\sfp^t\sfp=0$, while $\ker(\sfp)=\on{ran}(1-\sfp)$ is
isotropic if and only if $(1-\sfp)^t(1-\sfp)=0$. If both the
kernel and the range of $\sfp$ are isotropic, then
\[
1-(\sfp+\sfp^t)=(1-\sfp)^t(1-\sfp)-\sfp^t\sfp=0.
\]
Conversely, if $\sfp$ is a projection operator with
$\sfp+\sfp^t=1$, then $\sfp^t\sfp=(1-\sfp)\sfp=0$, and
similarly $ (1-\sfp)^t(1-\sfp)=0$.
\end{proof}

Again, we specialize to the case $W=\V$. Suppose $\V=E\oplus F$ is
a Lagrangian splitting, with associated projection operator
$\sfp$. The property $\sfp+\sfp^t=1$ implies that there is a
bivector $\pi\in \wedge^2 V$ defined by
%
%
\begin{equation}\label{eq:pisharp}
\pi^\sharp(\alpha)=-\pr_V(\sfp(\alpha)),\ \ \alpha\in V^*,
\end{equation}
that is,
$\pi(\alpha,\beta)=-\l
\sfp(\alpha),\,\beta\r=\l\alpha,\sfp(\beta)\r,\ \ \alpha,\beta\in V^*$.
If $\{e_i\}$ is a basis of $E$ and $\{f^i\}$ is the dual basis of
$F$, then
\begin{equation}\label{eq:pibasis}
 \pi=\hh \pr_V(e_i)\wedge \pr_V(f^i).
 \end{equation}
The graph of the bivector $\pi$ was
encountered in Example \ref{ex:wedge} above:
%

\begin{proposition}\label{prop:bi}
The graph of the bivector $\pi$ is given by
\begin{equation}\label{eq:graph1}
 \on{Gr}_\pi=E^\top \wedge F.
\end{equation}
In particular, $\on{ran}(\pi^\sharp)=\on{ran}(E)\cap \on{ran}(F)$, and
the symplectic 2-form on $\on{ran}(\pi^\sharp)$ is the difference of
the restrictions of the 2-forms on $\on{ran}(E),\on{ran}(F)$.
If $\phi,\psi$ are pure spinors defining $E,F$, then
\[
\phi^\top\wedge\psi=e^{-\iota(\pi)}(\phi^\top\wedge\psi)^{[\mathrm{top}]}.
\]
\end{proposition}
\begin{proof}
  Since both sides of \eqref{eq:graph1} are Lagrangian subspaces, it
  suffices to prove the inclusion $\supset$. Let $v\oplus \alpha\in
  E^\top \wedge F$.  Hence, there exist $\alpha_1,\alpha_2$ with
  $\alpha=\alpha_2-\alpha_1$ and $v\oplus \alpha_1\in E$, $v\oplus
  \alpha_2\in F$. Thus $v\oplus \alpha_1=-\sfp(\alpha)$, which implies
  that $\pi^\sharp(\alpha)=-\pr_V \sfp(\alpha)= v$.  The description
  of $\on{ran} \pi^\sharp=\on{ran}(\on{Gr}_\pi)$ is immediate from
  \eqref{eq:graph1}, see the discussion in Example \ref{ex:wedge}. The
  formula for $\phi^\top\wedge\psi$ follows since both sides are pure
  spinors defining the Lagrangian subspace $\on{Gr}_\pi$, with the
  same top degree part.
\end{proof}

\begin{proposition}\label{prop:gaugee}
Suppose $\V=E\oplus F$ is a Lagrangian splitting, defining a
bivector $\pi$. If $\varepsilon \in \wedge^2 E$, so that
$F_\varepsilon=A^{-\varepsilon}F$ is a new
      Lagrangian complement to $E$, the bivector $\pi_\varepsilon$ for the
splitting $E\oplus F_\varepsilon$ is given by
\[ \pi_\varepsilon=\pi+\pr_V(\varepsilon),\]
where $\pr_V\colon \wedge E \to \wedge V$ is the algebra homomorphism
extending the projection to $V$.
\end{proposition}
\begin{proof}
Let $\phi,\psi$ be pure spinors defining $E$, $F$. Then
$F_\varepsilon$ is defined by the pure spinor
$\psi_\varepsilon=\varrho(e^{-\varepsilon})\psi$. Using
Remark \ref{rem:reference}\eqref{it:aaa}, we obtain
\[ \phi^\top\wedge \psi_\varepsilon=\phi^\top \wedge \varrho(e^{-\varepsilon})\psi
=e^{-\iota(\pr_V(\varepsilon))} \phi^\top\wedge\psi.\]
The claim now
follows from \eqref{prop:bi}.
\end{proof}

\begin{proposition}\label{prop:transversality}
Let $(\Phi,\omega)\colon (V,E)\to (V',E')$ be a strong Dirac
morphism. Suppose $F'\in\on{Lag}(\V')$ is transverse to $E'$, and
$F$ is its backward image under $(\Phi,\omega)$. Then the
bivectors for the Lagrangian splittings $\V=E\oplus F$ and
$\V'=E'\oplus F'$ are $\Phi$-related:
\[\Phi(\pi)=\pi'.\]
\end{proposition}
\begin{proof}
To prove $\Phi(\pi)=\pi'$, we have to show that $(\Phi,0)\colon (V,\on{Gr}_\pi)\to
  (V',\on{Gr}_{\pi'})$ is a Dirac morphism:
\[ \Gamma_{(\Phi,0)}\circ (E^\top \wedge F)= (E')^\top \wedge F' .\]
Since both sides are Lagrangian, it suffices to prove the
inclusion $\supset$. If $v'\oplus \alpha'\in (E')^\top \wedge F'$,
then $\alpha'=\alpha_2'-\alpha_1'$, where $v'\oplus \alpha_1'\in E'$
and $v'\oplus \alpha_2'\in F'$. Since $(\Phi,\omega)$ is a strong Dirac
morphism for $E,E'$, there is a unique element $v\oplus \alpha_1\in E$
such that $v'=\Phi(v),\ \Phi^*(\alpha_1')=\alpha_1+\iota_v\om$.
Let $\alpha_2=\Phi^*(\alpha_2')-\iota_v\om$. Then $v\oplus \alpha_2\in F$
since $v\oplus \alpha_2\sim_{(\Phi,\omega)} v'\oplus \alpha_2$. Hence
$v\oplus \Phi^*(\alpha')=v\oplus (\alpha_2-\alpha_1)\in E^\top\wedge F$, proving
that $v'\oplus \alpha' \in \Gamma_{(\Phi,0)}\circ (E^\top \wedge F
)$.
\end{proof}

We next explain how a splitting $\V'=E'\oplus F'$ may be `pulled back'
under a linear map $\Phi\colon V\to V'$, given a bivector
$\pi\in\wedge^2V$ and a linear map $\mf{a}\colon E'\to V$ satisfying
suitable compatibility relations.

\begin{theorem}\label{prop:reconstruct}
Suppose that $\Phi\colon V\to V'$ is a linear map and
$\omega\in\wedge^2 V^*$ a 2-form. Given a Lagrangian splitting
$\V'=E'\oplus F'$, with associated projection
$\sfp'\in\End(\V')$, there is a 1-1 correspondence between
\begin{enumerate}
\item[(i)] Lagrangian subspaces $E\subset \V$ such that
$(\Phi,\omega)\colon (V,E)\to (V',E')$ is a strong Dirac morphism, and
\item[(ii)] Bivectors $\pi\in\wedge^2 V$ together with linear maps
$\mf{a}\colon E'\to V$, satisfying $\Phi\circ \mf{a}=\pr_{V'}\big|_{E'}$ and
\begin{equation}\label{eq:relation}
\pi^\sharp\circ \Phi^*=-\mf{a}\circ \sfp'\big|_{(V')^*}.
\end{equation}
\end{enumerate}
Under this correspondence, $\pi$ is the bivector defined by the
splitting $\V=E\oplus F$, where $F$ is the backward image of $F'$, and
$\mf{a}$ is the linear map defined by the strong Dirac morphism
$(\Phi,\omega)$ (see \eqref{eq:AA}).
\end{theorem}
\begin{proof}

``$(i) \Rightarrow (ii)$''.  By Proposition \ref{prop:transI}, we
know that the backward image $F$ of $F'$ is transverse to $E$. Let
$\sfp$ and $\sfp'$ be the projections defined by the
Lagrangian splittings $\V=E\oplus F$ and $\V'=E'\oplus F'$, and
$\pi,\pi'$ the corresponding bivectors.  As in \eqref{eq:AA}, the
strong Dirac morphism $(\Phi,\omega)$ defines a linear map
$\wh{\mf{a}}\colon E'\to E$, taking $w'\in E'$ to the unique
element $w\in E$ such that $w\sim_{(\Phi,\om)} w'$.  We claim that
for all $w\in \V,\ w'\in V'$,
\begin{equation}\label{eq:rel}
 w\sim_{(\Phi,\omega)} w' \Rightarrow \sfp(w)=\wh{\mf{a}}(\sfp'(w')).
\end{equation}
Indeed, let $w_1=\sfp(w)\in E$, so that $w_2=w-w_1\in F$.
There is a (unique) element $w_2'\in F'$ with
$w_2\sim_{(\Phi,\omega)} w_2'$, so let $w_1'=w'-w_2'$. Since
$w_2\sim_{(\Phi,\omega)} w_2'$, it follows that
$w_1\sim_{(\Phi,\omega)} w_1'$. Hence $w_1'\in E'$ by definition
of $E'$. It follows that
$ \sfp(w)=w_1=\wh{\mf{a}}(w_1')=\wh{\mf{a}}(\sfp'(w'))$,
as claimed. In particular, since
$\Phi^*\alpha'\sim_{(\Phi,\om)}\alpha'$ for $\alpha'\in V'$,
\eqref{eq:rel} implies that
\[\pi^\sharp(\Phi^*\alpha')=-\pr_V(\sfp(\Phi^*\alpha'))
=-\pr_V(\wh{\mf{a}}(\sfp'(\alpha')))
=-\mf{a}(\sfp'(\alpha')),\ \ \alpha'\in (V')^*\]
where $\mf{a}=\pr_V\circ \wh{\mf{a}}$.

``$(i) \Leftarrow (ii)$''. Our aim is to
construct the projection $\sfp$ with kernel
$F:=F'\circ\Gamma_{(\Phi,\om)}$ and range $E$. We define $\sfp$ by the
following equations, for $v,v_1,v_2\in V$ and
$\alpha,\alpha_1,\alpha_2\in V^*$:
\[ \begin{split}
\l\sfp(v_1),v_2\r &=\l \sfp'(\Phi(v_1)),\Phi(v_2)\r,
\\
\l\sfp(\alpha_1),\alpha_2\r&=-\pi(\alpha_1,\alpha_2),
\\ \l\sfp(v),\alpha\r&=\l\mf{a}^*\alpha,\,\Phi(v)\r+\pi(\iota_v\om,\alpha),
\\ \l\sfp(\alpha),v\r&=\l\alpha,v\r -\l\mf{a}^*\alpha,\,\Phi(v)\r-\pi(\iota_v\om,\alpha),
\end{split}\]
where $\mf{a}^*\colon V^*\to (E')^*=F'$ is the dual map to
$\mf{a}$. The linear map $\sfp$ defined in this way has the
property $\sfp+\sfp^t=1$. We claim that this linear map
satisfies \eqref{eq:rel}, where $\wh{\mf{a}}\colon E'\to \V$ is
defined as follows,
\[ \wh{\mf{a}}(w')=\mf{a}(w')\oplus
(\Phi^* \pr_{(V')^*}(w')-\iota_{\mf{a}(w')}\om).\]
For $w=v\oplus \iota_v\om,\ w'=\Phi(v)\oplus 0$, \eqref{eq:rel} is
easily checked using the definition of $\sfp$. Hence it suffices to
consider the case $w=\Phi^*\alpha',\ w'=\alpha'$ with $\alpha'\in
(V')^*$. For all $v\in V$, using the definition of $\sfp$ and $\Phi\circ
\mf{a}=\pr_{V'}|_{E'}$, i.e., $\mf{a}^*\circ
\Phi^*={(\sfp')^t}|_{(V')^*}$, we have:
\[ \begin{split}
\l\sfp(\Phi^*\alpha'),v\r&=
\l \alpha',\Phi(v)\r-\l (\sfp')^t \alpha',\Phi(v)\r
-\pi(\iota_v\om,\Phi^*\alpha')
\\&=\l
\sfp'\alpha',\Phi(v)\r
+\pi(\Phi^*\alpha',\iota_{v}\om)\\
\l\wh{\mf{a}}(\sfp'(\alpha')),v\r&=
\l \Phi^*\pr_{(V')^*}\sfp'(\alpha'),v\r-\om(\mf{a}(\sfp'(\alpha')),v)
\\&=\l\sfp'\alpha',\Phi(v)\r+\om(\pi^\sharp(\Phi^*\alpha'),v)
\end{split}\]
which shows
$\l\sfp(\Phi^*\alpha'),v\r=\l\wh{\mf{a}}(\sfp'(\alpha')),v\r$.
Similarly, for $\beta\in V^*$ we have, by \eqref{eq:relation},
\[ \l \sfp(\Phi^*\alpha'),\beta\r=
-\l \pi^\sharp (\Phi^*\alpha'),\beta\r= \l
\wh{\mf{a}}(\sfp'(\alpha')),\beta\r.\]
This proves \eqref{eq:rel}. Equation \eqref{eq:rel}
applies in particular to all elements $w\in F$, since these are by
definition $(\Phi,\om)$-related to elements $w'\in F'$. We hence see
that $\sfp(w)=0$ for all $w\in F$. This proves that $F\subset
\ker(\sfp)$.  Taking orthogonals, $\on{ran}(\sfp^t)\subset
F$. In particular, the range of $\sfp^t$ is isotropic, i.e.
$\sfp\sfp^t=0$, and hence $\sfp-\sfp^2
=\sfp(1-\sfp) =\sfp\sfp^t=0$. Thus $\sfp$ is a
projection. As before, we see that
$\ker\sfp=\on{ran}(1-\sfp)$ is isotropic as well, hence
$F=\ker(\sfp)$ since $F$ is maximal isotropic. It remains to
show that the Lagrangian subspace $E:=\on{ran}(\sfp)$ satisfies
$\Gamma_{(\Phi,\om)}\circ E\subset E'$. Suppose $w\sim_{(\Phi,\om)} w'$
for some $w\in E$. By \eqref{eq:rel}, we also have
$w=\sfp(w)\sim_{(\Phi,\om)} \sfp'(w')$. Thus
$0\sim_{(\Phi,\om)} (w'-\sfp'(w'))=(\sfp')^t(w')$.  Observe
that $\on{ran}(\Phi)\supset \Phi(\mf{a}(E'))=\on{ran}(E')$.  Hence
$\ker(\Phi^*)\subset \on{ann}(\on{ran}(E'))$.  Since $E'\cap
F'=0$, it follows that
\begin{equation}\label{eq:Ftrans}
\ker(\Phi^*)\cap \on{ann}(\on{ran}(F'))=0.
\end{equation}
Using Equation \eqref{eq:Ftrans}, the relation $0\sim_{(\Phi,\om)}
(\sfp')^t(w')\in F'$ implies that $(\sfp')^t(w')=0$, i.e.
$w'\in E'$.
\end{proof}

The proof shows that $\sfp|_V=\wh{\mf{a}}\circ \sfp'\circ \Phi$,
whereas $h:=\sfp|_{V^*}\colon V^*\to E$ is given by
\begin{equation}\label{eq:hdef}
h(\alpha)= (-\pi^\sharp(\alpha))\oplus
(\alpha-\Phi^*\pr_{(V')^*}\mf{a}^*(\alpha)-\iota(\pi^\sharp(\alpha))\om).
\end{equation}
It follows that $E=\on{ran}(\wh{\mf{a}})+\on{ran}(h)$. Projecting to
$V$, it follows in particular that
\begin{equation}\label{eq:Erange}
 \on{ran}(E)=\on{ran}(\mf{a})+\on{ran}(\pi^\sharp).
\end{equation}

\section{Pure spinors on manifolds}\label{sec:man}
A pure spinor on a manifold is simply a differential form whose
restriction to any point is a pure spinor on the tangent space.
The following discussion is carried out in the category of real
manifolds and $C^\infty$ vector bundles, but works equally well
for complex manifolds with holomorphic vector bundles.

\subsection{Dirac structures}\label{subsec:diracstru}
For any manifold $M$, we denote by $\TM=TM\oplus T^*M$ the direct
sum of the tangent and cotangent bundles, with fiberwise inner
product $\l\cdot,\cdot\r$. The fiberwise Clifford action
defines a bundle map
\begin{equation}\label{eq:spinoraction}
 \varrho\colon \Cl(\TM)\to \End(\wedge T^*M).\end{equation}
The same symbol will denote the action of sections of $\Cl(\TM)$ on
sections of $\wedge T^*M$, i.e. differential forms. The bilinear
pairing will be denoted by
\begin{equation}\label{eq:spinorpairing}
(\cdot,\cdot)_{\wedge T^*M} \colon \wedge T^*M\otimes\wedge T^*M\to \det(T^*M),
\end{equation}
and the same notation will be used for sections. Thus
$(\phi,\phi')_{\wedge T^*M}=(\phi^\top\wedge \phi')^{[\on{top}]}$ for
differential forms $\phi,\phi'\in\Gamma(\wedge T^*M)=\Om(M)$.  An
\emph{almost Dirac structure} on $M$ is a smooth Lagrangian subbundle
$E\subset \TM$. The pair $(M,E)$ is called an almost Dirac manifold.
A \emph{pure spinor defining $E$} is a nonvanishing differential form
$\phi\in\Omega(M)$ such that $\phi|_m$ is a pure spinor defining
$E_m$, for all $m$. Equivalently, $\phi$ is a nonvanishing section of
the line bundle $(\wedge T^*M)^E$.  Thus $E$ is globally represented
by a pure spinor if and only if the line bundle $(\wedge T^*M)^E$ is
orientable.  (Otherwise, one may still use pure spinors to describe
$E$ locally.)

Let $\eta\in \Om^3(M)$ be a closed 3-form. A direct computation shows
that the spinor representation defines a bilinear bracket
$\Cour{\cdot,\cdot}_\eta\colon \Gamma(\TM)\times \Gamma(\TM)\to \Gamma(\TM)$
by the condition:
\begin{equation}\label{eq:Courantbr}
\varrho(\Cour{x_1,x_2}_\eta)\psi=[[\d+\eta,\varrho(x_1)],\varrho(x_2)]\psi,\
\ \ \psi\in \Om(M),\ x_i\in\Gamma(\TM),
\end{equation}
where the brackets on the right-hand side are graded commutators of
operators on $\Omega(M)$. The bracket $\Cour{\cdot,\cdot}_\eta$ is the
\textit{$\eta$-twisted Courant bracket}
\cite{kli:wzw,sev:poi}. \footnote{This definition agrees with the
non-skew symmetric version of the Courant bracket
\cite{liu:ma,sev:poi}, called the \emph{Dorfman bracket} in
\cite{gua:ge}; the $\eta$-term in the bracket, however, differs from
the one in \cite{sev:poi} by a sign.}  (For more on the definition of
$\Cour{\cdot,\cdot}_\eta$ as a `derived bracket', see e.g.
\cite{al:der,kos:der,roy:co}.)  The operator on $\Omega(M)$ defined by
$$
[\varrho(x_1),[\varrho(x_2),[\varrho(x_3),\d+\eta]]]
$$
is multiplication by a function
\begin{equation}\label{eq:Upsilon}
\Upsilon(x_1,x_2,x_3)=-\l \Cour{x_3,x_2}_\eta,\,x_1\r\in C^\infty(M).
\end{equation}
Given an almost Dirac structure $E \subset \TM$, let $\Upsilon^E$
denote the restriction of the trilinear form $(x_1,x_2,x_3) \mapsto
\Upsilon(x_1,x_2,x_3)$ to the sections of $E$. In
contrast to $\Upsilon$, the trilinear form $\Upsilon^E$ is
\emph{tensorial} and
\emph{skew-symmetric}. The resulting element
\[\Upsilon^E\in \Gamma(\wedge^3 E^*)\]
is called the \textit{$\eta$-twisted Courant tensor} of
$E$.
\begin{definition}
  A \emph{Dirac structure} on a manifold $M$ is an almost Dirac
  structure $E$ together with a closed 3-form $\eta$ such that its
  $\eta$-twisted Courant tensor vanishes: $\Upsilon^E=0$.
   The triple $(M,E,\eta)$ is called a \emph{Dirac manifold}.
\end{definition}

For $E$ an almost Dirac structure one can always choose a
complementary almost Dirac structure $F$ such that $E \oplus
F=\TM$. (This is parallel to a well-known fact from symplectic
geometry \cite[Proposition 8.2]{ca:le}, with a similar proof.) As a
vector bundle, $F\cong E^*$ with pairing induced by the inner product
on $\TM$. We have:
\begin{proposition}
  Let $E$ be an almost Dirac structure on $M$, and  $F$ be
  a complementary almost Dirac structure.
  Suppose $E$ is represented by a pure spinor $\phi\in \Omega(M)$.
  Then there is a unique section $\sig^E\in \Gamma(E^*)$ (depending on
  $\phi$) such that
\[ (\d+\eta)\phi= \varrho(-\Upsilon^E+\sig^E)\phi .\]
Here we view $\Upsilon^E$ and $\sig^E$ as sections of
$\wedge F\subset \Cl(\TM)$.
\end{proposition}
\begin{proof}
  Choose a Lagrangian subbundle $F$ complementary to $E$. Since
\[ \Gamma(\wedge F)\to \Omega(M),\ x\mapsto \varrho(x)\phi\]
is an isomorphism, there is a unique odd element
$x\in\Gamma(\wedge F)\subset \Gamma(\T M)$ such that
$(\d+\eta)\phi=\varrho(x)\phi$. To see that $x$ has filtration
degree $3$, let $x_1,x_2,x_3$ be three sections of $E$. Since
$\varrho(x_i)\phi=0$, it follows that
\[\begin{split}
&\varrho([x_1,[x_2,[x_3,x]]])\phi=
[[[\varrho(x_1),[\varrho(x_2),[\varrho(x_3),\varrho(x)]]]\phi =\varrho(x_1x_2x_3)\varrho(x)\phi\\ &=  \varrho(x_1x_2x_3)(\d+\eta)\phi=[[[\varrho(x_1),[\varrho(x_2),[\varrho(x_3),\d+\eta]]]\phi =
\Upsilon^E(x_1,x_2,x_3)\phi, \end{split}
\]
proving that the Clifford commutator
$[x_1,[x_2,[x_3,x]]]=\iota(x_1)\iota(x_2)\iota(x_3)x$ (contraction
of $x\in\Gamma(\wedge(E^*))$ with sections of $E$) is a scalar.
This implies that $x$ has filtration degree 3, and that the degree
3 part of $x$ is $- \Upsilon^E$.
\end{proof}
We hence see that an almost Dirac structure $E\subset \TM$
is integrable if and only if
\[ (\d+\eta)\phi\in \varrho(\TM)\phi,\]
for any pure spinor $\phi\in\Om(M)$ (locally) representing $E$.
The characterization of the integrability condition $\Upsilon^E=0$
in terms of pure spinors was observed by Gualtieri \cite{gua:ge}, see
also \cite{al:der}.

Examples of Dirac structures (for a given $\eta$) include graphs of
2-forms $\omega\in \Omega(M)$ with $\d\omega=\eta$, as well as graphs
of bivector fields $\pi\in\mf{X}^2(M)$ defining $\eta$-\emph{twisted
  Poisson structures} \cite{kli:wzw,sev:poi} in the sense that
$\frac{1}{2}[\pi,\pi]_{\on{Sch}}+\pi^\sharp(\eta)=0$.  One may also
consider \emph{complex} Dirac structures on $M$, given by complex
Lagrangian subbundles $E\subset \TM^\C$ satisfying $\Upsilon^E=0$.
The defining pure spinors are complex-valued differential forms $\phi$
on $M$, given as nonvanishing sections of $(\wedge T^*M^\C)^E$. If $E$
is a Dirac structure, then its image $E^c$ under the complex
conjugation mapping is a Dirac structure defined by the complex
conjugate spinor $\phi^c$. $E$ is called a \emph{generalized complex
  structure} \cite{hi:gen,gua:ge} if $E\cap E^c=0$.

Suppose $E\subset \TM$ is a Dirac structure.  The vanishing of the
Courant tensor implies that $E$ is a \emph{Lie algebroid}, with anchor
given by the natural projection on $TM$, and Lie bracket
$[\cdot,\cdot]_E$ on $\Gamma(E)$ given by the restriction of the
Courant bracket $\Cour{\cdot,\cdot}_\eta$.
%
From the theory of Lie algebroids, it follows that the generalized
distribution $\on{ran}(E)$ is integrable (in the sense of Sussmann)
\cite{zun:po}. The generalized foliation having $\on{ran}(E)$ as its
tangent distribution is called the \emph{Dirac foliation}.  For any
leaf $Q\subset M$ of the Dirac foliation, the collection of 2-forms on
$T_mQ$ (defined as in \eqref{eq;2form}) defines a smooth 2-form
$\omega_Q\in\Om^2(Q)$ with
\[ \d\omega_Q=i_Q^*\eta,\]
where $i_Q\colon Q\to M$ is the inclusion (for a proof, see e.g.
\cite[Proposition 6.10]{me:lec}). If $E$ is the graph of a Poisson
bivector $\pi$ (with $\eta=0$), this is the usual symplectic
foliation.

\subsection{Dirac morphisms}\label{subsec:morph}
Suppose $\Phi\colon M\to M'$ is a smooth map, and $\omega\in \Om^2(M)$
is a 2-form.  As in the linear case, we view the pair $(\Phi,\om)$ as
a `morphism', with composition rule \eqref{eq:compo}.  Given sections
$x\in \Gamma(\TM)$ and $x'\in \Gamma(\TM')$, we will write
\[ x\sim_{(\Phi,\omega)} x'\  \Leftrightarrow \ \ \forall m\in M\colon
x_m\sim_{((d\Phi)_m,\omega_m)} x'_{\Phi(m)}.\]
In terms of the spinor representation, this is equivalent to the condition
\[ e^{\omega} \Phi^*(\varrho(x')\psi')=\varrho(x)(e^{\omega}\Phi^*(\psi')),\ \
\psi' \in \Omega(M').\] Using the definition \eqref{eq:Courantbr}
of the Courant bracket as a derived bracket, one obtains:
\begin{lemma}[Stienon-Xu]\cite[Lemma~2.2]{sti:red} \label{lem:sti}
  Let $M,M'$ be manifolds with closed 3-forms $\eta,\eta'$,
  $\Phi\colon M\to M'$ a smooth map, and $\om\in \Om^2(M)$ a 2-form
  such that $\Phi^*\eta'=\eta+\d\om$. Then
\[ x_i\sim_{(\Phi,\omega)} x_i',\ i=1,2\ \Rightarrow
\Cour{x_1,x_2}_\eta\sim_{(\Phi,\omega)} \Cour{x_1',x_2'}_{\eta'}.\]
\end{lemma}
That is, the morphism $(\Phi,\om)\colon M\to M'$ intertwines both
the inner product and the ($\eta$- resp. $\eta'$-twisted) Courant brackets
on $\TM$ and $\TM'$.
\begin{definition}
\begin{enumerate}
\item
Suppose $(M,E)$ and $(M',E')$ are almost Dirac manifolds. A
morphism $(\Phi,\om)\colon M\to M'$ is called a
\emph{(strong) almost Dirac morphism}  $(\Phi,\om)\colon (M,E)\to
(M,E')$ if  %
$((\d\Phi)_m,\omega_m)\colon (T_mM,E_m)\to
(T_{\Phi(m)}M',E'_{\Phi(m)})$
is a linear (strong) Dirac morphism for all $m\in M$.
\item
Suppose $(M,E,\eta)$ and $(M',E',\eta')$ are Dirac manifolds. A (strong)
almost Dirac  morphism $(\Phi,\om)\colon M\to M'$ is called a
\emph{(strong) Dirac morphism} $(\Phi,\omega)\colon (M,E,\eta) \to
(M',E',\eta')$ if $\eta+\d\om=\Phi^*\eta'$.
\end{enumerate}
\end{definition}
For $\omega=0$, strong Dirac morphisms coincide with the \emph{Dirac
  realizations} of \cite{bur:di}.

\begin{example}
  If $(M,E,\eta)$ is a Dirac manifold, then so is
  $(M,A^{-\omega}(E),\eta+\d\omega)$, for any 2-form $\om$, and
  $(\on{id}_M,\om)$ is a Dirac morphism between the two.  The Dirac
  structures $E$ and $A^{-\omega}(E)$ are isomorphic as Lie
  algebroids; in particular, they define the same Dirac foliation.
  However, the 2-forms on the leaves of this foliation change by the
  pull-back of $\omega$.
\end{example}
\begin{example}
Any manifold $M$ can be trivially viewed as a Dirac manifold
$M=(M,TM,0)$. A strong Dirac morphism from $M$ to $\pt$ is then
the same thing as a symplectic 2-form on $M$. More generally,
strong Dirac morphisms $M\to N$ are (special types of) symplectic
fibrations.
\end{example}

\begin{example}\label{ex:leaves}
If $(M,E,\eta)$ is a Dirac manifold, and $Q\subset M$ is a leaf of the
associated foliation of $M$, then the inclusion map defines a strong
Dirac morphism $(\iota_Q,\om_Q)\colon (Q,TQ,0)\to (M,E,\eta)$.
\end{example}
From the linear case, it follows that a strong almost Dirac morphism
gives rise to a bundle map
\[
\wh{\mf{a}}\colon \Phi^* E'\to E.
\]
This is indeed a smooth bundle map: the projection $\T M\oplus
\Phi^*\T M'\to \Phi^*\T M'$ restricts to a bundle isomorphism
$\Gamma_{\Phi}\cap (E\oplus\Phi^*\T M')\to \Phi^* E'$, and
$\wh{\mf{a}}$ is the inverse of this bundle isomorphism followed by
the projection to $\T M$. We let
\begin{equation}\label{eq:A}
\mf{a}=\pr_{TM}\circ \wh{\mf{a}}
\colon\  \Phi^* E'\to \on{ran}(E)\subset TM
\end{equation}
\begin{proposition}\label{prop:comor}
  Suppose $(\Phi,\omega)\colon (M,E,\eta)\to (M',E',\eta')$ is a
  strong Dirac morphism. Then the induced
  bundle map $\wh{\mf{a}}\colon \Phi^*E'\to E$ is a \emph{comorphism of Lie
    algebroids} \cite{mac:gen}. That is,
  it is compatible with the anchor maps in the sense that
  \[\d\Phi\circ \mf{a}=\pr_{\Phi^*TM'}|_{\Phi^*E'},\] and the
  induced map on sections
\[\wh{\mf{a}}\colon \Gamma(E')\to \Gamma(E),\ (\wh{\mf{a}}(x'))_{m}=\wh{\mf{a}}(x'_{\Phi(m)})
\]
preserves brackets.
\end{proposition}
\begin{proof}
  Compatibility with the anchor is obvious.  If $x_1', x_2'$ are
  section of $E'$, then
(using Lemma \ref{lem:sti}) both
$\wh{\mf{a}}(\Phi^*[x_1',x_2']_{E'})$ and
$[\wh{\mf{a}}(\Phi^*x_1'),\wh{\mf{a}}(\Phi^*x_2')]_E$ are sections of
$E$ which are $(\Phi,\omega)$-related to $[x_1',x_2']_{E'}$. Hence their
difference is $(\Phi,\omega)$-related to $0$. Since $(\Phi,\omega)$ is
a strong Dirac morphism, it follows that the difference is in fact $0$.
\end{proof}

The second part of
Proposition \ref{prop:comor} shows that
\eqref{eq:A} defines a Lie algebra homomorphism $\mf{a}\colon
\Gamma(E')\to \mf{X}(M)$. That is, the strong Dirac morphism defines
an `\emph{action}' of the Lie algebroid $E'$ on the manifold $M$.

\subsection{Bivector fields}
From the linear theory, we see that any Lagrangian splitting
$\TM=E\oplus F$ defines a bivector field $\pi$ on $M$.
Furthermore,
\[
e^{-\iota(\pi)}(\phi^\top\wedge\psi)^{[\mathrm{top}]}=\phi^\top\wedge\psi
\]
for any pure spinors $\phi,\psi$ defining $E,F$. Recall that
$(\phi^\top\wedge\psi)^{[\mathrm{top}]}$ is a volume form on $M$.

For an arbitrary volume form $\mu$ on $M$, and any bivector field
$\pi\in\mf{X}^2(M)$, one has the formula \cite{ev:po}
\begin{equation}\label{eq:evlu}
 \d (e^{-\iota(\pi)} \mu)=\iota\big(-\hh [\pi,\pi]_{\on{Sch}}+
X_\pi\big) (e^{-\iota(\pi)} \mu).\end{equation}
Here $[\cdot,\cdot]_{\on{Sch}}$ is the Schouten bracket on
multivector fields, and $X_\pi$ is the vector field on $M$
defined by $\d\iota(\pi)\mu=-\iota(X_\pi)\mu$.  If $\pi$ is a
Poisson bivector field, then $X_\pi\in\mf{X}(M)$ is called the
\emph{modular vector
  field} of $\pi$ with respect to the volume form $\mu$
\cite{we:mod}.
(See \cite{kos:mod} for modular vector fields for \emph{twisted}
Poisson structures.)

\begin{theorem}\label{th:bivector}
Let $\pi$ be the bivector field
defined by the Lagrangian splitting $\TM=E\oplus F$. Let
$\Upsilon^E\in\Gamma(\wedge^3 F)$ and $\Upsilon^F \in
\Gamma(\wedge^3 E)$ be the Courant tensor fields of $E,F$.
\begin{itemize}
\item[a)] The Schouten bracket of $\pi$ with itself is given by
the formula
\[\hh [\pi,\pi]_{\on{Sch}}= \pr_{TM}(\Upsilon^E)+\pr_{TM}(\Upsilon^F),\]
where $\pr_{TM}\colon \wedge E\to \wedge TM$ is the algebra homomorphism extending the projection $E\to TM$, and similarly for
$\pr_{TM}\colon \wedge F\to \wedge TM$.
\item[b)] Given pure spinors $\phi,\psi\in\Om(M)$ defining $E,F$, let
$\sig^E\in \Gamma(F)$ and $\sig^F\in \Gamma(E)$ be the unique sections
such that
\[ (\d+\eta)\phi=\varrho(-\Upsilon^E+\sig^E)\phi,\ \ \ \
(\d+\eta)\psi=\varrho(-\Upsilon^F+\sig^F)\psi.\]
Then the vector field $X_\pi$ defined using the volume form
$\mu=(\phi^\top\wedge\psi)^{[\mathrm{top}]}$ is given by
\[X_\pi=\pr_{TM}(\sig^F)-\pr_{TM}(\sig^E).\]
\end{itemize}
\end{theorem}
\begin{proof}
  We may assume that $E,F$ are globally defined by pure spinors
  $\phi,\psi$.  Using Remark \ref{rem:reference}\eqref{it:aaa}, we have
\[\begin{split}
\d(\phi^\top\wedge\psi)&=(-1)^{|\phi|} \big(\phi^\top \wedge\d\psi+(\d\phi)^\top\wedge \psi \big)\\
&=(-1)^{|\phi|} \big(\phi^\top \wedge(\d+\eta)\psi+((\d+\eta)\phi)^\top\wedge \psi \big)\\
&=(-1)^{|\phi|} \big(\phi^\top\wedge
(\varrho(-\Upsilon^F+\sig^F)\psi)+(\varrho(-\Upsilon^E+\sig^E)\phi)^\top\wedge\psi \big)\\
&=\iota(\pr_{TM}(-\Upsilon^F+\sig^F)+\pr_{TM}(-\Upsilon^E-\sig^E))(\phi^\top\wedge\psi).
\end{split}\]
On the other hand, $\phi^\top\wedge\psi=e^{-\iota(\pi)}\mu$ gives, by \eqref{eq:evlu},
\[\d(\phi^\top\wedge \psi)= \iota(-\hh
[\pi,\pi]_{\on{Sch}}+X_\pi)(\phi^\top\wedge \psi).\] Applying the star
operator $\star$ for $\mu$, and using that $\star(\phi^\top\wedge \psi)$ is
invertible, it follows that
\[\pr_{TM}(-\Upsilon^F+\sig^F)+\pr_{TM}(-\Upsilon^E-\sig^E)=-\hh
[\pi,\pi]_{\on{Sch}}+X_\pi.\]
\end{proof}
As a special case, if both $E,F$ are Dirac structures (i.e.
integrable), then the corresponding bivector field $\pi$ satisfies
$[\pi,\pi]_{\on{Sch}}=0$, i.e., it is a Poisson structure.  The
symplectic leaves of $\pi$ are the intersections of the leaves of the
Dirac structures $E$ with those of $F$.  The fact that transverse
Dirac structures (or equivalently \emph{Lie bialgebroids}) define
Poisson structures goes back to Mackenzie-Xu \cite{mac:lie}.
\begin{proposition}\label{prop:diracmaps}
  Suppose $(\Phi,\om)\colon (M,E)\to (M',E')$ is an almost Dirac
  morphism, and let $F'\subset \TM'$ be a Lagrangian subbundle
  complementary to $E'$. Then there is a smooth Lagrangian
  subbundle $F\subset \TM$ complementary to $E$, with the property
  that for all $m\in M$, $F_m$ is the backward image of $F'_{\Phi(m)}$
  under $(\d_m\Phi,\om_m)$. Furthermore:
\begin{enumerate}
\item The bivector fields $\pi,\pi'$ defined by the splittings
      $\TM=E\oplus F$ and $\TM'=E'\oplus F'$ satisfy
      \[ \pi\sim_\Phi \pi',\]
i.e. $(\d \Phi)_m\pi_m=\pi'_{\Phi(m)}$ for all $m\in M$.
\item The Courant tensors $\Upsilon^{F}\in\Gamma(\wedge^3 E)$ and
      $\Upsilon^{F'}\in\Gamma(\wedge^3 E')$ are related by
\[ \Upsilon^F=\wh{\mf{a}}(\Phi^*\Upsilon^{F'}),\]
using the extension of $\wh{\mf{a}}\colon \Gamma(\Phi^* E')\to
\Gamma(E)$ to the exterior algebras.
\item The bivector field $\pi$ satisfies
\[ \hh [\pi,\pi]_{\on{Sch}}= \mf{a}(\Phi^*\Upsilon^{F'})+\pr_{TM}(\Upsilon^E),\]
using the extension of $\mf{a}\colon \Gamma(\Phi^* E')\to
\Gamma(TM)$ to the exterior algebras.
\item
\[ \pi^\sharp\circ \Phi^*=-\mf{a}\circ \sfp'\colon T^*M'\to TM,\]
where $\sfp'\colon \TM'\to E'$ is the projection along $F'$.
\item If $\psi'$ is a pure spinor defining $F'$, and $\psi=e^\om
      \Phi^*\psi'$ the corresponding pure spinor defining $F$, the
      sections $\sig^F,\ \sig^{F'}$ are related by
      $\sig_F=\wh{\mf{a}}(\Phi^*\sig_F)$, that is,
\[ \sig^F\sim_{(\Phi,\om)}\sig^{F'}.\]
\end{enumerate}
\end{proposition}
\begin{proof}
  Let $\psi'\in \Om(M')$ be a pure spinor (locally) representing $F'$.
  From the linear case (Proposition \ref{prop:transI}), it follows
  that $\psi=e^\om \Phi^*\psi'$ is non-zero everywhere, and is a
  pure spinor representing a Lagrangian subbundle $F\subset \TM$
  transverse to $E$.  Now (a) follows from the linear case, see
  Proposition \ref{prop:transversality}. We next verify (b), at any
  given point $m\in M$. Let $m'=\Phi(m)$. Given $(x_i)_m\in F_m$ for
  $i=1,2,3$, let  $(x_i')_{m'}\in F'_{m'}$ with
  \[(x_i)_m\sim_{((d\Phi)_m,\ \omega_m)} (x_i')_{m'}.\] Choose sections
  $x_i\in\Gamma(F),\ x_i'\in\Gamma(F')$ extending the given values at
  $m,\ m'$. We have to show
$ \Upsilon^F(x_1,x_2,x_3)|_m=\Upsilon^{F'}(x_1',x_2',x_3')|_{m'}$.
We  calculate:
\[
\Upsilon^F(x_1,x_2,x_3)\ \psi
=\varrho(x_1 x_2 x_3)\,(\d+\eta)(e^\om \Phi^* \psi')
=\varrho(x_1 x_2 x_3)\,e^\om \Phi^* (\d+\eta')\psi'\]
On the other hand,
\[
(\Phi^* \Upsilon^{F'}(x_1',x_2',x_3'))\ \psi=
e^\om \Phi^*  \Upsilon^{F'}(x_1',x_2',x_3')\ \psi'=e^\om \Phi^* \varrho(x_1' x_2' x_3')\,(\d+\eta')\psi'
\]
These two expressions coincide at $m$, proving (b). Theorem
\ref{th:bivector} together with (b) implies the statement (c). Part
(d) follows from Proposition \ref{prop:reconstruct}. Part (e) follows
from (b) together with the definition of $\sigma^F$, $\sigma^{F'}$.
\end{proof}

Part (b) shows in particular that if $F'$ is a Dirac structure,
transverse to $E'$, then its backward image is again a Dirac
structure.

\subsection{Dirac cohomology}\label{subsec:diraccoh}
In this Section, we will discuss certain cohomology groups associated
with any pair of transverse Dirac structures $E,F\subset \TM$ and a
given volume form $\mu$ on $M$. We assume that $E,F$ are given by pure
spinors $\phi,\psi$, normalized by the condition $(\phi,\psi)_{\wedge
  T^*M}=\mu$. Let $\sigma^E \in \Gamma(F), \sigma^F \in \Gamma(E)$ be
sections defined as in Theorem \ref{th:bivector}, and denote
\[\sigma=\sigma^F-\sigma^E\in\Gamma(\TM).\]
Replacing $\phi,\psi$ with $\ti{\phi}=f\phi,\ \ti{\psi}=f^{-1}\psi$,
for $f$ a
nonvanishing function on $M$, this section changes by a closed 1-form:
\begin{equation}\label{eq:sigmachange}
\ti{\sigma}=\sig-f^{-1}\d f.
\end{equation}
Indeed, letting let $\mathsf{p}$ be the projection from $\TM$ to $E$
along $F$ we have $\tilde{\sigma}^F=\sigma^F-\mathsf{p}(f^{-1}\d f),\
\ \tilde{\sigma}^E=\sigma^E+(I-\mathsf{p})(f^{-1}\d f)$.

We define the {\em Dirac cohomology
  groups} associated to a triple $(E,F,\mu)$ as the cohomology of the
operators
\[
\dirac_+=d+\eta+\varrho(\sigma),\ \  \ \dirac_-=d+\eta-\varrho(\sigma)
\]
on $\Omega(M)$, restricted to the subspace on which they square to
zero:
\begin{equation}\label{eq:diraccoh}
H_\pm(E,F,\mu):=\ker(\dirac_\pm)/\ker(\dirac_\pm)\cap
\on{im}(\dirac_\pm)\equiv H(\ker\dirac_\pm^2,\dirac_\pm).
\end{equation}
The pure spinors $\phi,\psi$ define classes in $ H_+(E,F,\mu)$ and
$H_-(E,F,\mu)$, respectively, since
$\dirac_+\phi=0$ and $\dirac_- \psi=0$.
The Dirac cohomology groups
are independent of the choice of defining spinors $\phi,\psi$:
Changing the pure spinors by a function $f$ as above,
\eqref{eq:sigmachange} shows that the operators $\dirac_\pm$ change
by conjugation,
$\tilde{\dirac}_+=f\dirac_+ f^{-1}$ and
$\tilde{\dirac}_-=f^{-1} \dirac_- f$.

\begin{example}
  Let $M$ be a manifold with volume form $\mu$. Consider transverse
  Dirac structures $E=\on{Gr}_\omega$ for some closed 2-form $\omega$,
  and $F=T^*M$. In this case, one can choose $\phi=e^{-\omega},
  \psi=\mu$.  We obtain $\eta=0$, $\sigma=0$, $\dirac_\pm
  =\d$, and the Dirac cohomology groups $H_\pm (TM, T^*M,
  \mu)$ coincide with the de Rham cohomology of
  $M$.
\end{example}

\begin{example}
  Let $M$ be a manifold with volume form $\mu$ and with a Poisson
  bivector $\pi$.  Let $E=TM, F=\on{Gr}_\pi$. The choice $\phi=1,
  \psi=e^{-\iota(\pi)}\mu$ gives $\dirac_-=d-\iota(X_\pi)$, where
  $X_\pi$ is the modular vector field.  The operator
  $\dirac_-^2=-\L(X_\pi)$ vanishes on differential forms invariant
  under the flow generated by $X_\pi$.  The Dirac cohomology $H_-(TM,
  \on{Gr}_\pi, \mu)= H(\Omega(M)^{X_\pi}, d-\iota(X_\pi))$ resembles
  the Cartan model of equivariant cohomology for circle actions.
\end{example}

Let $\pi$ be the Poisson structure defined by the splitting
$\TM=E\oplus F$, and $X_\pi=\pr_{TM}\sig$ the modular vector field.
Let
\begin{equation}
H_\pi(M)=H(\Omega(M)^{X_\pi}, d-\iota(X_\pi)).
\end{equation}
By
Remark \ref{rem:reference}\eqref{it:aaa} there is a pairing
\[
H_+(E,F,\mu) \otimes H_-(E,F,\mu) \to H_\pi(M)
\]
given on representatives by the formula $u \otimes v \mapsto u^\top \wedge v$.
The pure spinors $\phi, \psi$ define cohomology classes $[\phi]\in H_+(E,F,\mu),
[\psi]\in H_-(E,F,\mu)$, and $[\phi^\top \wedge \psi] \in H_\pi(M)$.
If $M$ is compact, the integration map $\int_M\colon
\Omega(M)^{X_\pi}\to \R$ descends to $H_\pi(M)$. Hence
\[ \int_M \phi^\top \wedge\psi=\int_M\mu>0\]
shows that the cohomology classes $[\phi]\in H_+(E,F,\mu), [\psi]\in
H_-(E,F,\mu)$ are both nonzero.

There is the following version of
functoriality with respect to strong Dirac morphisms for Dirac
cohomology.

\begin{proposition}  \label{prop:psipsi'}
  Let $(\Phi, \omega)\colon (M, E, \eta) \to (M', E', \eta')$ be a
  strong Dirac morphism, and let $F'\subset \TM'$ be a Dirac
  structure transverse to $E'$, with backward image $F$. Assume that
  $E,E'$ are defined by pure spinors $\phi,\phi'$ such that the
  corresponding sections $\sig^{E}$ and $\sig^{E'}$ vanish.  Let
  $\psi'$ and $\psi=e^\om \Phi^*\psi'$ be pure spinors defining $F'$
  and $F$, and let $\mu'$ and $\mu$ be the resulting volume forms.
  Then $e^\omega \circ \Phi^*$ intertwines $\dirac_-$ and $\dirac'_-$,
  and hence induces a map in Dirac cohomology $e^\omega \Phi^*\colon
  H_-(E',F',\mu') \to H_-(E,F,\mu)$ taking $[\psi']$ to $[\psi]$.
\end{proposition}

\begin{proof}
Since $\sig^E,\sig^{E'}$ vanish we have
$\sig=\sig^F$ and $\sig'=\sig^{F'}$. By Proposition
\ref{prop:diracmaps} (e), the map $e^\om \Phi^*$ intertwines the
Clifford actions of $\sig^F$ and $\sig^{F'}$, while on the other hand
this map also intertwines $\d+\eta$ with $\d+\eta'$. Hence  it
intertwines $\dirac_-$ with $\dirac'_-$.
\end{proof}

\subsection{Classical dynamical Yang-Baxter equation}\label{subsec:cdybe}
The following result describes the Courant tensor of Lagrangian
subbundles defined by elements in $\Gamma(\wedge^2E)$.

\begin{proposition}[Liu-Weinstein-Xu \cite{liu:ma}]\label{prop:LWX}
  Let $\T M=E\oplus F$ be a splitting into Lagrangian subbundles,
  where both $E,F$ are integrable relative to the closed 3-form
  $\eta$, and let us identify $F^*=E$.  Given a section $\varepsilon \in\Gamma(\wedge^2 E)$,
  defining a section $A^{-\varepsilon}\in\Gamma(\on{O}(\TM))$, let
  $F_\varepsilon=A^{-\varepsilon}(F)$ be the Lagrangian subbundle spanned by the sections
  $x+\iota_x \varepsilon$ for $x\in \Gamma(F)=\Gamma(E^*)$.
Then the Courant tensor $\Upsilon_\varepsilon\in\Gamma(\wedge^3
E)$ of $F_\varepsilon$ is given by the formula:
\[   \Upsilon_\varepsilon=\d_F \varepsilon+\hh [\varepsilon,\varepsilon]_E.\]
Here $[\cdot,\cdot]_E$ is the Lie algebroid bracket of $E$, and
$\d_F\colon \Gamma(\wedge^\bullet F^*)\to \Gamma(\wedge^{\bullet + 1}
F^*)$ is the Lie algebroid differential of $F$.
\end{proposition}

\begin{remark}
  The result in \cite{liu:ma} is stated only for $\eta=0$. However, since
  the statement is local, one may use a gauge transformation by a
  local primitive of $\eta$ to reduce to this case.
\end{remark}

We are interested in the following special case: Let $M=\g^*$,
with its standard linear Poisson structure $\pi_{\g^*}\in
\Gamma(\wedge^2T\g^*)= C^\infty(\g^*)\otimes \wedge^2\g^*$, and
put $F=T\g^*$ and $E=\on{Gr}_{\pi_\g^*}$. The bundle $E$ is
spanned by sections $\A_0(\xi)\oplus \l\theta_0,\xi\r$ for
$\xi\in\g$, where $\A_0(\xi)$ is the generating vector fields for
the co-adjoint action, and $\l\theta_0,\xi\r\in\Om^1(\g^*)$ is the
`constant' 1-form defined by $\xi$. The trivialization
$E=\g^*\times\g$ defined by these sections identifies $E$ with the
action algebroid for the co-adjoint action: The bracket on
$\Gamma(E)=C^\infty(\g^*,\g)$ is defined by the Lie bracket on
$\g$ via the Leibniz rule, and the anchor map is given by the
action map $\A_0\colon\g\to T\g^*$. For
$\varepsilon\in\Gamma(\wedge^2 E)$, the bracket
$[\varepsilon,\varepsilon]_E$ is given by the Schouten bracket on
$\wedge\g$. On the other hand we may view $\varepsilon\in
C^\infty(\g^*,\wedge^2 \g)$ as a 2-form on $\g^*$, and then $\d
\varepsilon=\d_F \varepsilon$ is just its exterior differential.
The resulting equation reads
\[ \d \varepsilon+\hh [\varepsilon,\varepsilon]_{\on{Sch}}=\Upsilon_\varepsilon.\]
If $\Upsilon_\varepsilon$ is a multiple of the structure constants
tensor, this is a special case of the \emph{classical dynamical
Yang-Baxter equation} (CDYBE) \cite{al:cli,et:ge}. We will see
below how a solution arises from the Cartan-Dirac structure on
$G$.

For more information on the relation between Dirac structures and
the CDYBE, see the work of Liu-Xu \cite{liu:dir} and
Bangoura-Kosmann-Schwarzbach \cite{ban:eq}.


\vskip.3in\section{Dirac structures on Lie groups}\label{sec:dirgroup}
In this Section, we will study Dirac structures over Lie groups
$G$ with bi-invariant pseudo-Riemannian metrics. This will be
based on the existence of a canonical isomorphism
\[ \TG\cong G\times (\g\oplus \ol{\g})\]
preserving scalar products and Courant brackets.  In the
subsequent section, we will describe a corresponding isomorphism
of spinor modules.

\subsection{The isomorphism $\TG\cong
G\times(\g\oplus\ol{\g})$}\label{subsec:isom} Let $G$ be a Lie
group (not necessarily connected), and let $\g$ be its Lie
algebra. We denote by $\xi^L,\xi^R\in \mf{X}(G)$ the left-,
right-invariant vector fields on $G$ which are equal to $\xi \in
\g=T_eG$ at the group unit.  Let $\theta^L,\theta^R\in
\Om^1(G)\otimes\g$ be the left-, right-Maurer-Cartan forms, i.e.
$\iota({\xi^L})\theta^L=\iota(\xi^R)\theta^R=\xi$. They are related by
$\theta^R_g=\Ad_g(\theta^L_g)$, for all $g\in G$. The adjoint action
of $G$ on itself will be denoted $\A_{{\sf{ad}}}$ (or simply $\A$, if
there is no risk of confusion). The
corresponding infinitesimal action is given by the vector fields
\[ \A_{{\sf{ad}}}(\xi)=\xi^L-\xi^R.\]
Suppose that the Lie algebra $\g$ of $G$ carries an \emph{invariant
  inner product}. By this we mean an $\Ad$-invariant, non-degenerate
symmetric bilinear form $B$, not necessarily positive definite.
Equivalently, $B$ defines a bi-invariant pseudo-Riemannian metric
on $G$. Given $B$, we can define the bi-invariant 3-form
$\eta\in\Om^3(G)$,
\[ \eta:=\f{1}{12}B(\theta^L,[\theta^L,\theta^L])
\]
Since $\eta$ is bi-invariant, it is closed, and so it defines an
$\eta$-twisted Courant bracket $\Cour{\cdot,\cdot}_{\eta}$ on $G$. The
conjugation action $\A_{{\sf{ad}}}$ extends to an action of $D=G\times
G$ on $G$, by
\begin{equation}\label{eq:actA}
\A\colon D\to \on{Diff}(G),\ \  \A(a,a')=l_{a'} \circ r_{a^{-1}},
\end{equation}
where $l_a(g)=ag$ and $r_a(g)=ga$. The corresponding infinitesimal action
\[
\A\colon \mf{d}\to \mf{X}(G),\ \
\A(\xi,\xi')=\xi^L-(\xi')^R
\]
lifts to a map
\begin{equation}\label{eq:s}
\sfs \colon \mf{d}\to \Gamma(\TG),\ \ \ \ \sfs(\xi,\xi')=\sfs^L(\xi)+\sfs^R(\xi'),
\end{equation}
where
\[ \sfs^L(\xi)=\xi^L\oplus \hh B(\theta^L,\xi),\ \
\sfs^R(\xi')=-(\xi')^R\oplus \hh B(\theta^R,\xi').
\]

Let us equip $\dd$ with the bilinear form $B_\dd$ given by $+B$ on
the first $\g$-summand and $-B$ on the second $\g$-summand. Thus
$\dd=\g\oplus\ol{\g}$ is an example of a Lie algebra with
invariant split bilinear form.
\begin{proposition}
  The map $\sfs\colon \dd\to \Gamma(\TG)$ is $D$-equivariant, and satisfies
\begin{equation}\label{eq:scom}
 \l \sfs(\zeta_1),\sfs(\zeta_2)\r=B_{\dd}(\zeta_1,\zeta_2),\ \ \
\Cour{\sfs(\zeta_1),\sfs(\zeta_2)}_{\eta}=\sfs([\zeta_1,\zeta_2])\end{equation}
for all $\zeta_1,\zeta_2\in\dd$. Furthermore,
\begin{equation}\label{eq:Ups1}
\Upsilon(\sfs(\zeta_1),\sfs(\zeta_2),\sfs(\zeta_3))= B_{\dd}\big(\zeta_1,[\zeta_2,\zeta_3]\big)
\end{equation}
for all $\zeta_i\in\dd$, where $\Upsilon\colon
\Gamma(\TG)^{\otimes 3}\to C^\infty(G)$ was defined in
\eqref{eq:Upsilon}.
\end{proposition}
\begin{proof}
The $D$-equivariance of the map $\sfs$ is clear.
Let $\varrho$ be the Clifford action of $\TG$ on $\wedge T^*G$. We
have $ [\varrho(\sfs^L(\xi)),\d+\eta]=\L(\xi^L)$ and
$[\varrho(\sfs^R(\xi)),\d+\eta]=-\L(\xi^R)$, thus
\[ [\d+\eta,\,\varrho(\sfs(\zeta))]=\L(\A(\zeta))\]
for all $\zeta\in\dd$. This proves the second Equation in \eqref{eq:scom}, while the first Equation is obvious. Finally, \eqref{eq:Ups1} follows from
\eqref{eq:scom} and the definition of $\Upsilon$.
Hence,
\[ \varrho(\Cour{\sfs(\zeta_1),\sfs(\zeta_2)}_{\eta})=
[[\d+\eta,\varrho(\sfs(\zeta_1)],\,\varrho(\sfs(\zeta_2))]=
\varrho(\sfs([\zeta_1,\zeta_2])).
\]
%
%
\end{proof}
Put differently, the map $s$ defines a $D$-equivariant isometric isomorphism
\begin{equation}\label{eq:isom}
\TG\cong G\times\dd,\end{equation}
identifying the $\eta$-twisted Courant bracket on $\TG$ with the
unique Courant bracket on $ G\times \dd$ which agrees with the Lie
bracket on $\dd$ on constant sections.

\subsection{$\eta$-twisted Dirac structures on $G$}
Using \eqref{eq:isom}, we see that any Lagrangian
subspace $\s\subset \dd$ defines a Lagrangian subbundle
\[ E^\s\cong G\times\s,\]
spanned by the sections $\sfs(\zeta)$ with $\zeta\in\s$.  The
Lagrangian subbundle $E^\s$ is invariant under the action of the
subgroup of $D$ preserving $\s$. Let
$\Upsilon^\s\in\wedge^3\s^*$ be defined as
\begin{equation}\label{eq:Upsh}
 \Upsilon^\s(\zeta_1,\zeta_2,\zeta_3)=B_{\dd}\big(\zeta_1,[\zeta_2,\zeta_3]\big),\ \zeta_i\in\s.\end{equation}
By \eqref{eq:Ups1}, the Courant tensor $\Upsilon^{E^\s}$ is just
$\Upsilon^\s$, using the sections $\sfs$ to identify
$(E^\s)^*\cong G\times \s^*$.  In particular, we see that $\s$
defines a Dirac structure if and only if $\s$ is a Lie subalgebra.
To summarize:
\begin{quote}\emph{Any Lagrangian
subalgebra  $\s\subset \dd$ defines an $\eta$-twisted Dirac
structure $E^\s$.} \end{quote}

The Dirac structure $E^\s$ is invariant under the action of any Lie
subgroup normalizing $\s$, and in particular under the action of the
subgroup $S\subset D$ integrating $\s$. As a Lie algebroid, $E^\s$ is
just the action algebroid for this $S$-action. In particular, its
leaves are just the components of the $S$-orbits on $G$.  The 2-form
on the orbit $\O=\A(S)g$ of an element $g\in G$ is the $S$-invariant
form $\om_\O$ given as follows: for $\zeta_i=(\xi_i,\xi_i')\in\s$,
\begin{equation}\label{eq:genform}
 \begin{split} \om_\O(\A(\zeta_1),\A(\zeta_2))|_g&= \hh \l B(\theta^L,
\xi_1)+B(\theta^R,\xi_1'),\ \xi_2^L-(\xi_2')^R\r\\ &=\hh
B(\xi_2-\Ad_{g^{-1}}\xi_2',\xi_1+\Ad_{g^{-1}}\xi_1')\\ &=\hh (B(\Ad_g
\xi_2,\xi_1')-B(\xi_2',\Ad_g\xi_1)),
\end{split}
\end{equation}
using $B(\xi_1,\xi_2)=B(\xi_1',\xi_2')$ since $\mf{s}$ is Lagrangian.
By the general theory from Section \ref{subsec:diracstru}, these 2-forms satisfy
$\d\om_\O=\iota_\O^*\eta$, where $\iota_\O\colon\O\to G$ is the
inclusion. The kernel of $\om_\O$ equals $\ker(E^\s)$, i.e. it is
spanned by all $\A(\zeta)$ such that the $T^*G$-component of
$\sfs(\zeta)$ is zero:
\begin{equation}\label{eq:kernelformula}
\ker(\om_\O|_g)=\{\A(\zeta)|_g\ |\ \zeta=(\xi,\xi')\in\s,\ \ \Ad_g\xi+\xi'=0\}. \end{equation}
\begin{remark}
  For $\g$ a complex semi-simple Lie algebra, a complete
  classification of Lagrangian subalgebras of $\dd$ was obtained by
  Karolinsky \cite{kar:cla}.   The Poisson geometry of the variety of Lagrangian
  subalgebras of $\dd$ was studied in detail by Evens--Lu
  \cite{ev:on2}.
\end{remark}

\begin{remark}
  If $\dd=\s_1\oplus\s_2$ is a splitting into two Lagrangian
  subalgebras (i.e., $(\dd,\s_1,\s_2)$ is a Manin triple), one obtains
  two transverse Dirac structures $E^{\s_1},E^{\s_2}$. As discussed
  after Theorem \ref{th:bivector}, such a pair of transverse Dirac
  structures gives rise to a Poisson structure on $G$, with symplectic
  leaves the intersections of the orbits of $S_1,S_2$.  For $\g$ a
  complex semi-simple Lie algebra, the Manin triples were classified
  by Delorme \cite{del:cla}. See Evens--Lu \cite{ev:on2} for a wealth
  of information regarding Poisson structures obtained from Lagrangian
  subalgebras. An example will be worked out in Section
  \ref{subsec:gauss} below.
\end{remark}

\begin{remark}
  We may also use this construction to obtain generalized complex (and
  K\"ahler) structures \cite{gua:ge} on even-dimensional real Lie
  groups $K$, with complexification $G=K^\C$.  Indeed, let $\s\subset
  \dd=\g\oplus \overline{\g}$ be a Lagrangian subalgebra such that
\begin{equation}
\s\cap \s^c=\{0\},
\end{equation}
where $\s^c$ denotes the
  complex conjugate of $\s$. Then the associated Dirac structure $E^\s
  \subset \T G$ satisfies $E^\s\cap (E^\s)^c=\{0\}$ along $K$.
  Hence it defines a
  generalized complex structure on $K$. For a concrete example,
  suppose $K$ is compact, and let $\g=\n_-\oplus \t\oplus \n_+$ be a triangular
  decomposition. (That is, $\t=\t_K^\C$ is the complexification of a
  maximal Abelian subalgebra, and $\n_+,\n_-$ are the sums of the
  positive, negative root spaces). Then
\[
\s= (\n^+\oplus 0)\oplus \mathfrak{l} \oplus (0\oplus \n^-)\subset
\dd=\g\oplus \bar{\g}
\]
has the desired property, for any Lagrangian subspace $\mathfrak{l}
\subset \t\oplus \bar{\t}$ with
$\mathfrak{l}\cap \mathfrak{l}^c=\{0\}$ (i.e., $\mathfrak{l}$ is a linear
generalized complex structure on the vector space $\t_K$). The
generalized complex structures on Lie groups considered in
Gualtieri \cite[Example~6.39]{gua:ge} are examples of this
construction.
\end{remark}

\subsection{The Cartan-Dirac structure}
The simplest example of a Lagrangian subalgebra is the diagonal
$\s= \g_\Delta\hra \dd$, with corresponding $S$ the diagonal
subgroup $G_\Delta\subset D$. The associated Dirac structure $E_G$
is spanned by the sections ${\sfe}(\xi):=\sfs(\xi,\xi)$:
\begin{equation}\label{eq:cartan-dirac}
 E_G=\on{span}\,\{{\sfe}(\xi)\,|\,\xi\in\g\}\subset \TG,
\end{equation}
\[ \begin{split}
{\sfe}(\xi)=(\xi^L-\xi^R, B(\textstyle{\f{\theta^L+\theta^R}{2}},\xi)\big).
\end{split}
\]
We
call $E_G$ the \emph{Cartan-Dirac structure}, see
\cite{bur:int,sev:poi,ko:dir}. This Dirac structure was introduced
independently by Alekseev, \v{S}evera, and Strobl in the
mid-1990's. The $G_\Delta\cong G$-action is just the action by conjugation
on $G$, hence the Dirac foliation is given by the conjugacy classes $\Co\subset G$.
The formula \eqref{eq:genform} specializes
 to the 2-form on
conjugacy classes introduced in \cite{gu:gr}:
\[ \om_\Co(\A_{{\sf{ad}}}(\xi_1),\A_{{\sf{ad}}}(\xi_2))=-\hh B((\Ad_g-\Ad_{g^{-1}})\xi_1,\xi_2),\]
The kernel at $g\in\Co$ is the span of vector fields
$\A_{{\sf{ad}}}(\xi)|_g$ with $\Ad_g\xi+\xi=0$. The anti-diagonal in
$\g\oplus \ol{\g}$ is a $G$-invariant Lagrangian complement to the
diagonal, and hence defines a $G$-invariant Lagrangian subbundle $F_G$
complementary to $E_G$, spanned by $\sff(\xi)=\sfs(\xi/2,-\xi/2)$:
\begin{equation}\label{eq:FG}
F_G=\on{span}\{\sff(\xi)|\ \xi\in\g\}\subset \TG,
\end{equation}
\[\sff(\xi)=\big(\textstyle{\f{\xi^L+\xi^R}{2}},B(\textstyle{\f{\theta^L-\theta^R}{4}},\xi)\big).\]
The $1/2$ factors in the definition of $\sff(\xi)$ are introduced so that
$\l \sfe(\xi),\sff(\xi')\r=B(\xi,\xi')$.

Let $\tri\in\wedge^3(\g)$
be the structure constants tensor of $\g$, normalized as follows:
\begin{equation}\label{eq:Xi}
 \iota(\xi_3)\iota(\xi_2)\iota(\xi_1)\tri={\ts \frac{1}{4}}\,B(\xi_1,[\xi_2,\xi_3]_\g).\end{equation}
Let $\sfe\colon \wedge\g\to \Gamma(\wedge E_G)$ be the extension of
$\sfe\colon \g\to \Gamma(E_G)$ as an algebra
homomorphism. Thus $\sfe(\tri)$ is a section of $\wedge^3(E_G)$.
\begin{lemma}\label{lem:upsilon}
The Courant tensor of $F_G$ is given by :
\[ \Upsilon^{F_G} =\sfe(\tri).\]
\end{lemma}
\begin{proof}
This follows from \eqref{eq:Ups1} since
$B_\dd(\zeta_1,[\zeta_2,\zeta_3]_\dd)=\f{1}{4}B(\xi_1,[\xi_2,\xi_3]_\g)$
for $\zeta_i=(\xi_i/2,\,-\xi_i/2)$.
\end{proof}
The element $\tri$ also defines a trivector field,
$\A_{{\sf{ad}}}(\tri)\in\mf{X}^3(G)$. Theorem \ref{th:bivector}
implies that the bivector field $\pi_G\in \mf{X}^2(G)$ defined by
the Lagrangian splitting $\TG=E\oplus F$ satisfies
\[
\hh [\pi_G,\pi_G]_{\on{Sch}}= \A_{{\sf{ad}}}(\tri).
\]
To give an explicit formula for $\pi_G$, let $v_i,v^i$ be $B$-dual bases of $\g$,
i.e. $B(v_i,v^j)=\delta_i^j$.

\begin{proposition}\label{prop:pig}
The bivector field $\pi_G$ is given by
\begin{equation}\label{eq:pig}
\pi_G=\hh \sum_i v^{i,L}\wedge v_i^R.\end{equation}
\end{proposition}
\begin{proof}
By \eqref{eq:pibasis}, we have
\[\pi_G=\hh \sum_i \big((v_i)^L-(v_i)^R\big)\wedge \f{(v^i)^L+(v^i)^R}{2}.\]
Since $\sum_i
v^{i,L}\wedge v_i^L=\sum_i  v^{i,R}\wedge v_i^R$, this simplifies to
the expression in \eqref{eq:pig}
\end{proof}
The bivector field $\pi_G$ was first considered in \cite{al:ma,al:qu}.


\subsection{Group multiplication}
In this Section, we will examine the behavior of the Cartan-Dirac
structure under group multiplication,
\[ \Mult\colon G\times G\to G,\ \ (a,b)\mapsto ab.\]
For any differential form $\beta\in\Om(G)$, we will denote by
$\beta^i\in\Omega(G\times G)$ its pull-back to the $i$'th factor,
for $i=1,2$.  We will use similar notation for vector fields on
$G\times G$, and for sections of the bundle $\T(G\times G)$. Let
$\varsigma\in \Om^2(G\times G)$ denote the 2-form
\begin{equation}\label{eq:varsigma}
 \varsigma=-\hh B(\theta^{L,1},\theta^{R,2}).
\end{equation}
A direct computation shows that
\begin{eqnarray}
\Mult^* \eta&=&\eta^1+\eta^2+\d \varsigma, \label{eq:id1}
\end{eqnarray}
hence we have a multiplication morphism
\[(\on{Mult},\varsigma)\colon
(G,\eta)\times(G,\eta) =(G\times G,\eta^1+\eta^2)\to
(G,\eta).\]

\begin{remark}
This is expressed more conceptually in terms of the simplicial
model $B_pG=G^p$ of the classifying space $BG$. Let $\partial_i
\colon G^p\to G^{p-1},\ 0\le
 i\le p$ be the `face maps' given as
$ \partial_i(g_1,\ldots,g_p)=(g_1,\ldots,g_ig_{i+1},\ldots,g_p),$
while $\partial_0$ omits the first entry $g_1$, and $\partial_p$
omits the last entry $g_p$. Let
$\delta=\sum_{i=0}^p\partial_i^*\colon \Om^\bullet(G^{p-1}) \to
\Om^\bullet(G^p)$. Then $\delta$ commutes with the de-Rham
differential, turning $\bigoplus_{p,q} \Om^q(G^p)$ into a double
complex. The total differential on $\Om^q(G^p)$ is
$\d+(-1)^q\delta$. Then $\eta\in\Om^3(G)$ and $\varsigma\in
\Om^2(G^2)$ define a cocycle of degree 4 (see \cite{we:sy}):
\begin{equation}\label{eq:4cycle}
 \d\eta=0,\ \ \partial\eta=-\d\varsigma,\ \ \partial\varsigma=0.\end{equation}
(If $G$ is compact, simple, and simply connected, and $B$ the
basic inner product, this pair is the Bott-Shulman representative
of the generator of $H^4(BG)\cong H^3(G)$.) The second condition
is just the property \eqref{eq:id1} used above. Using the third
property, one may verify that the multiplication morphism is
associative, in the sense that
\[ (\on{Mult},\varsigma)\circ \big((\on{Mult},\varsigma)\times (\on{id}_G,0)\big)
=(\on{Mult},\varsigma)\circ \big((\on{id}_G,0)\times
(\on{Mult},\varsigma)\big).\]
\end{remark}

We will compare the morphism $(\on{Mult},\varsigma)$ with the
groupoid multiplication of $\dd$, viewed as the pair groupoid over
$\g$: writing $\zeta=(\xi,\xi'),\ \zeta_i=(\xi_i,\xi_i'),\ i=1,2$,
the groupoid multiplication is
\[ \zeta=\zeta_2\circ \zeta_1\  \Leftrightarrow\
\xi=\xi_2,\ \xi'=\xi_1',\ \xi_2'=\xi_1.
\]
\begin{proposition}\label{prop:related}
The isomorphism $G\times\dd \to \TG$ defined by $\sfs$ intertwines
the groupoid multiplication of $\dd$ with the morphism
$(\on{Mult},\varsigma)$, in the sense that
\begin{equation}\label{eq:cute}
 \zeta_2\circ \zeta_1=\zeta\ \  \Leftrightarrow\ \
\sfs^1(\zeta_1)+\sfs^2(\zeta_2)\sim_{(\on{Mult},\varsigma)} \sfs(\zeta),
\end{equation}
for $\zeta,\zeta_1,\zeta_2\in \dd$.
\end{proposition}
\begin{proof}
Spelling out the relations \eqref{eq:cute}, we have to show that,
for all $\xi\in \g$,
\begin{equation}\label{eq:rrelated}
\begin{split}
\sfs^{R,1}(\xi)&\sim_{(\on{Mult},\varsigma)} \sfs^R(\xi),\ \ \ \
\sfs^{L,2}(\xi)\sim_{(\on{Mult},\varsigma)} \sfs^L(\xi),\\
\sfs^{L,1}(\xi)&+\sfs^{R,2}(\xi)\sim_{(\on{Mult},\varsigma)}0.
\end{split}
\end{equation}
The equivariance properties
\[ \begin{split} &\on{Mult}(g a, b)=g\on{Mult}(a,b),\ \
\on{Mult}(a,bg^{-1})=\on{Mult}(a,b)g^{-1},\\
&\on{Mult}(ag^{-1},gb)=\on{Mult}(a,b)\ \  \end{split} \] of the
multiplication map imply the following relations of generating
vector fields:
\[ -\xi^{R,1}\sim_{\on{Mult}}-\xi^R,\ \
\xi^{L,2}\sim_{\on{Mult}}\xi^L,\ \
\xi^{L,1}-\xi^{R,2}\sim_{\on{Mult}} 0.\]
This proves the `vector field part' of the relations
\eqref{eq:rrelated}.  The 1-form part is equivalent to the following
three identities, which are verified by a direct computation:
\[ \begin{split}
&\hh B(\theta^{R,1},\xi) + \iota(-\xi^{R,1})\varsigma=\hh \on{Mult}^* B(\theta^R,\xi),
\\
&
\hh B(\theta^{L,2},\xi) +\iota(\xi^{L,2})\varsigma=\hh \on{Mult}^* B(\theta^L,\xi)\, , \\
&\hh B(\theta^{L,1}+\theta^{R,2} ,\xi) +  \iota(\xi^{L,1}-\xi^{R,2})\varsigma=0 .
\end{split} \]
\end{proof}

\begin{theorem}\label{prop:multcardir}
The multiplication map $\on{Mult}\colon G\times G\to G$
extends to a strong Dirac morphism
\[ (\on{Mult},\varsigma)\colon (G,E_G,\eta)\times (G,E_G,\eta)\to
(G,E_G,\eta),\]
with $\varsigma\in\Om^2(G\times G)$ as defined above. In terms of
the trivialization $E_G=G\times\g$, the map $\wh{\a}\colon
\on{Mult}^*E_G\to E_G\times E_G$ associated with the strong Dirac
morphism is given by the diagonal embedding $\g\to \g\times\g$.
Similarly, the inversion map $\on{Inv}\colon G\to G,\ g\mapsto
g^{-1}$ extends to a Dirac morphism
\[(\on{Inv},0)\colon (G,E_G,\eta)\to (G,E_G^\top,-\eta).\]
\end{theorem}
\begin{proof}
By Proposition \ref{prop:related}, the sections $\sfe(\xi)=\sfs(\xi,\xi)$ satisfy
\[\sfe^1(\xi)+\sfe^2(\xi)\sim_{(\on{Mult},\varsigma)} \sfe(\xi).\]
This shows that $(\on{Mult},\varsigma)$ is a Dirac morphism. For any
given point $(a,b)\in G\times G$, no non-trivial linear combination of
$\sfe^1(\xi)|_a,\ \sfe^2(\xi')|_b$ is $(\on{Mult},\varsigma)$-related
to $0$. Hence, the Dirac morphism $(\on{Mult},\varsigma)$ is strong.

We have $\on{Inv}^*B(\theta^L+\theta^R,\xi)=-B(\theta^L+\theta^R,\xi)$
and $\xi^L-\xi^R\sim_{\on{Inv}}(\xi^L-\xi^R)$. Hence
\[ \sfe(\xi)\sim_{(\on{Inv},0)} \sfe(\xi)^\top\]
where $\sfe(\xi)^\top$ is the image of $\sfe(\xi)$ under the map
$(v,\alpha)\to (v,-\alpha)$. Since $\on{Inv}^*\eta=-\eta$,
this shows that $(\on{Inv},0)\colon (G,E_G,\eta)\to
(G,E_G^\top,-\eta)$ is a Dirac morphism.
\end{proof}

\begin{remark}
More generally, suppose that $\mf{s}\subset \dd$ is a Lagrangian
subalgebra, defining a Dirac structure $E^{\mf{s}}$. Since
$\g_\Delta\circ \mf{s}=\mf{s}$, the same argument as in the proof
above shows that $(\on{Mult},\varsigma)$ is a strong Dirac
morphism from $(G,E_G,\eta)\times (G,E^{\mf{s}},\eta)$ to
$(G,E^{\mf{s}},\eta)$.
\end{remark}


Let $\wti{F}_{G\times G}\subset \T(G\times G)$ be the backward
image of the Lagrangian subbundle $F_G$ under
$(\on{Mult},\varsigma)$.  Since $F_G$ is spanned by the sections
$\sff(\xi)=\hh(\sfs^L(\xi)-\sfs^R(\xi))$, \eqref{eq:rrelated} shows that
$\wti{F}_{G\times G}$ is spanned by the sections
\begin{equation}\label{eq:spanned}
 \hh (\sfs^{L,2}(\xi)-\sfs^{R,1}(\xi)),\ \ \
 \hh(\sfs^{L,1}(\xi)+\sfs^{R,2}(\xi)).\end{equation}
Since $F_G$ is a complement to $E_G$, its backward image
$\wti{F}_{G\times G}$ is a complement to $E_G^1\oplus E_G^2$ (see
Proposition \ref{prop:transI}). Let us describe the element of
$\wedge^2(E_G^1\oplus E_G^2)$ relating $\wti{F}_{G\times G}$ to
the standard complement $F_G^1\oplus F_G^2$. Let $v_i\in\g$ and
$v^i\in\g$ be $B$-dual bases, and put
\begin{equation}\label{eq:gamma} \gamma=\hh (v_i)^1\wedge
  (v^i)^2\in \wedge^2(\g\oplus \g). \end{equation}
Let
\begin{equation}\label{eq:egamma}
\sfe(\gamma)=\hh \sum_i \sfe^1(v_i)\wedge \sfe^2(v^i)\in
\Gamma(\wedge^2(E^1_G\oplus E^2_G))
\end{equation}
be the corresponding section.
\begin{proposition}\label{prop:egamma}
The Lagrangian complement $\wti{F}_{G\times G}=F_G\circ
\Gamma_{(\Mult,\varsigma)}$ is obtained from $F_G^1\oplus F_G^2$
by the bivector $\sfe(\gamma)$:
\[ \wti{F}_{G\times G}=A^{-\sfe(\gamma)}(F_G^1\oplus F_G^2).\]
\end{proposition}
\begin{proof}
We compute $\iota(\sff^1(\xi))\sfe(\gamma)=\sfe(\iota^1(\xi)\gamma)=\hh
\sfe^2(\xi)=\hh(\sfs^{L,2}(\xi)+\sfs^{R,2}(\xi))$.  Thus
\[ \sff^1(\xi)+\iota(\sff^1(\xi))\sfe(\gamma)
=\hh(\sfs^{L,1}(\xi)-\sfs^{R,1}(\xi)+ \sfs^{L,2}(\xi)+\sfs^{R,2}(\xi))\]
is the sum of the sections in \eqref{eq:spanned}. Similarly, we
find that $\sff^2(\xi)+\iota(\sff^2(\xi))\sfe(\gamma)$ is the difference of the
sections in \eqref{eq:spanned}.
\end{proof}

The bivector field on $G\times G$ corresponding to the splitting
$(E_G^1\times E_G^2) \oplus A^{-\sfe(\gamma)} (F_G^1\times F_G^2)$ of
$\T(G\times G)$ is given by (see Proposition \ref{prop:gaugee}(i)),
\begin{equation}
\wti{\pi}=\pi^1_G+\pi^2_G+\A^{12}_{{\sf{ad}}}(\gamma),
\end{equation}
where $\pi_G$ is the bivector field for the splitting $\T
G=E_G\oplus F_G$, and $\A^{12}_{{\sf{ad}}}=\A^1_{{\sf{ad}}}\oplus
\A^2_{{\sf{ad}}}\colon \g\oplus \g\to \mf{X}(G\times G)$.
By Proposition \ref{prop:diracmaps}(c) we have
$\ti{\pi}\sim_{\on{Mult}}\pi$.
Furthermore, Proposition \ref{prop:diracmaps}(c) and Lemma
\ref{lem:upsilon}, imply that the Schouten bracket $\hh
[\wti{\pi},\wti{\pi}]_{\on{Sch}}$ equals the trivector field
$\A^{\on{diag}}_{{\sf{ad}}}(\tri)$, where $\A^{\on{diag}}_{{\sf{ad}}}$
is the diagonal action on $G\times G$.

\subsection{Exponential map}\label{subsec:exp1}
We will now discuss the behavior of the Cartan-Dirac structure under the exponential map,
\[ \exp\colon \g\to G.\]
Let $\g_\natural\subset \g$ denote the set of regular points of the
exponential map, that is, all points where $\d\exp$ is an isomorphism.
We begin with some preliminaries concerning $\T\g^*$, not using the
inner product on $\g$ for the time being. Let $\A_0$ be the action of
$D_0:=\g^*\rtimes G$ on $\g^*$ by
\[\A_0(\beta,g)\nu=(\Ad_{g^{-1}})^*\nu-\beta.\]
This action lifts to an action by automorphisms of $\T\g^*$,
preserving the inner product as well as the (untwisted) Courant
bracket. Let $\dd_0=\g^*\rtimes \g$ be the Lie algebra of $D_0$,
equipped with the canonical inner product defined by the pairing, and
let $\A_0\colon \dd_0\to \mf{X}(\g^*)$ be the infinitesimal action. To
simplify notation, we denote the constant vector field defined by
$\beta\in\g^*$ by $\beta_0=\A_0(\beta,0)$, and write
$\A_0(\xi)=\A_0(0,\xi)$. Let $\theta_0\in \Om^1(\g^*)\otimes \g^*$ be
the tautological 1-form, defined by $\iota(\beta_0)\theta_0=\beta$.
Consider the $D_0$-equivariant map
\begin{equation}\label{eq:s0}
\sfs_0\colon \dd_0\to \Gamma(\T\g^*),\ \
 \sfs_0(\beta,\xi)=\A_0(\beta,\xi)\oplus \l\theta_0,\xi\r.\end{equation}
Then $\l \sfs_0(\zeta),\sfs_0(\zeta')\r=B_{\dd_0}(\zeta,\zeta')$,
showing that $\sfs_0$ defines a $D_0$-equivariant isometric
isomorphism
\[ \T\g^*\cong \g^*\times \dd_0.\]
A direct computation shows that this isomorphism is compatible
with the Courant bracket $\Cour{\cdot,\cdot}_0$ on $\T\g^*$ and the
Lie bracket on $\dd_0$.

Since $\g\subset \dd_0$ is a Lagrangian Lie subalgebra, the
sections $\sfe_0(\xi):=\sfs_0(0,\xi)$ span a Dirac structure
$E_{\g^*}\subset \T\g^*$.  Since $E_{\g^*}\cap T\g^*=0$, this Dirac
structure is of the form $E_{\g^*}=\on{Gr}_{\pi_{\g^*}}$ for a
Poisson bivector field $\pi_{\g^*}$ satisfying
\begin{equation}\label{eq:kirillov}
\iota(\l\theta_0,\xi\r)\pi_{\g^*}=\A_0(\xi),\ \ \xi\in\g.
\end{equation}
The Poisson structure $\pi_{\g^*}$ is just the standard linear
Poisson structure on $\g^*$. Similarly, the sections
$\sff_0(\beta):=\sfs_0(\beta,0)$ span the Lagrangian subspace
$F_{\g^*}=T\g^*$, which is complementary to $E_{\g^*}$.

Let us now use the invariant inner product $B$ on $\g$ to identify
$\g^*\cong\g$. Let
\begin{equation}\label{eq:varpi}
\varpi\in\Om^2(\g),\ \ \d\varpi=\exp^*\eta
\end{equation}
be the primitive of $\exp^*\eta\in\Om^3(\g)$ defined by the de
Rham homotopy operator for the radial homotopy.
\begin{proposition}\label{prop:exp}
The sections $\sfe_0(\xi)$ and $\sfe(\xi)$ are
$(\exp,\varpi)$-related:
\begin{equation}\label{eq:erel}
 \sfe_0(\xi)\sim_{(\exp,\varpi)} \sfe(\xi).
\end{equation}
Similarly, over the subset $\g_\natural\subset\g$, one has
\begin{equation}\label{eq:frel}
\sff_0(\xi)+\sfe_0(C\,\xi)\sim_{(\exp,\varpi)} \sff(\xi),
\end{equation}
where $C\colon \g_\natural\to \End(\g)$ is given by the formula:
\begin{equation}\label{eq:sform}
 C|_\nu=\big(1/2\coth(z/2)-1/z\big)\big|_{z=\ad_\nu},\ \
 \nu\in\g_\natural.
\end{equation}
\end{proposition}
\begin{proof}
  Recall that $\beta_0$ denotes the `constant vector field'
  $\A_0(\beta,0)$. We extend the notation $(\cdot)_0$ to
  $\g^*\cong\g$-valued functions on $\g^*\cong\g$: For instance, the
  vector field corresponding to the function $\nu\mapsto
  -\ad_\xi\nu=\ad_\nu\xi$ is $(\ad_\nu \xi)_0=\A_0(\xi)$.

The vector field part of the relation \eqref{eq:erel} says that
$\A_0(\xi)\sim_{\exp} \xi^L-\xi^R=\A_{{\sf{ad}}}(\xi)$, which follows by the
$G$-equivariance of $\exp$. The 1-form part of \eqref{eq:erel} is
equivalent to the following property \cite{al:mom} of $\varpi$:
\[ \iota(\A_0(\xi))\varpi=\hh
\exp^* B(\theta^L+\theta^R,\xi)-B(\theta_0,\xi).\]
Since $\exp$ is a local diffeomorphism over $\g_\natural$, the
section $\sff(\xi)$ of $\TG$ is $(\exp,\varpi)$-related to a
unique section $\wt{\sff}(\xi)$ of $\T\g|_{\g_\natural}$. Since
inner products are preserved under the $(\exp,\varpi)$-relation
(see \eqref{eq:isometric})  we have
\[
\l\sfe_0(\xi'),\wt{\sff}_0(\xi)\r=\l
\sfe(\xi'),\sff_0(\xi)\r=B(\xi',\xi)=\l \sfe_0(\xi'),\sff_0(\xi)\r
\]
for all $\xi'\in \g$, showing that the $F_\g$-component of
$\wt{\sff}_0(\xi)$ is equal to $\sff_0(\xi)$.  It follows that
$\wt{\sff}_0(\xi)=\sff_0(\xi)+\sfe_0(C(\xi))$, where $C$ is
defined by $B(\xi',C (\xi))=\l \sff_0(\xi'),\wt{\sff}_0(\xi)\r$.
To compute $C$, we re-write \eqref{eq:frel} in the equivalent form
(using \eqref{eq:erel}):
\[ \sff_0(\xi)\sim_{(\exp,\varpi)} \sff(\xi)-\sfe(C(\xi))\]
Again, we write out the vector field and 1-form parts of this
relation:
\begin{equation}\label{eq:showthis}
 \begin{split}
\xi_0&=\hh\exp^*(\xi^L+\xi^R)-\A_0(C(\xi)),\\
\iota(\xi_0)\varpi&=\ts{\f{1}{4}}
\exp^*B(\theta^L-\theta^R,\xi)-\hh
\exp^*B(\theta^L+\theta^R,C(\xi)).
\end{split}\end{equation}
We now verify that $C$ given by \eqref{eq:sform} satisfies these
two equations. Let $T,U^L,U^R\colon \g\to \End(\g)$ be the functions defined by
\[ \iota(\xi_0)\varpi=B(\theta_0,T\,\xi),\ \ \
\exp^*\theta^L=U^L\,\theta_0,\ \ \ \exp^*\theta^R=U^R\,\theta_0.\]
It is known that (for the first identity, see e.g. \cite{me:ge})
%
\[T|_\nu=\left.\left(\f{\sinh(z)-z}{z^2}\right )\right|_{z=\ad_\nu},\
\ \
U^L|_\nu= \left.\left(\f{1-e^{-z}}{z}\right )\right|_{z=\ad_\nu},\
\ U^R|_\nu= \left.\left(\f{e^z-1}{z}\right )\right|_{z=\ad_\nu}.
\]

Note that $U^L$ and $U^R$ are transposes relative to the inner
product on $\g$, and that they are invertible over $\g_\natural$.
Their definitions imply that
\[
\exp^*\xi^L=((U^L)^{-1}\,\xi)_0,\ \
\exp^*\xi^R=((U^R)^{-1}\,\xi)_0.
\]
The first equation in \eqref{eq:showthis} becomes
\[
\xi_0=\left(\left(\f{(U^L)^{-1}+(U^R)^{-1}}{2}-\ad_\nu
\,C\right)\xi\right)_0
\]
which follows from the identity
\[ 1=\f{1}{2}\left(\f{z}{1-e^{-z}}+\f{z}{e^z-1}\right)-z\left(\f{1}{2}\coth\left(\f{z}{2}\right)
-\f{1}{z}\right).\]
In a similar fashion, the second equation in \eqref{eq:showthis}
follows from the identity
\[
\f{\sinh(z)-z}{z^2}=
\f{1}{4}\left(\f{e^z-1}{z}-\f{1-e^{-z}}{z}\right)
-\f{1}{2}\left(\f{e^z-1}{z}+\f{1-e^{-z}}{z}\right)\left(\f{1}{2}\coth(\f{z}{2})
-\f{1}{z}\right).
\]
\end{proof}
As an immediate consequence of \eqref{eq:erel}, we obtain
\begin{theorem}\label{th:exp}
The exponential map and the 2-form $\varpi$ define a Dirac
morphism
\[
(\exp,\varpi)\colon (\g,E_\g,0)\to (G,E_G,\eta).
\]
It is a strong Dirac morphism over the open subset
$\g_\natural\subset \g$.
\end{theorem}

Let $\wti{F}_\g$ be the backward image (defined over
$\g_\natural$) of $F_G$ under $(\exp,\varpi)$, and let
\[
\varepsilon\in C^\infty(\g_\natural,\wedge^2\g)
\]
be the unique map such that the associated orthogonal
transformation $A^{-\sfe_0(\varepsilon)}\in
\Gamma(\on{O}(\T\g_\natural))$ takes $F_\g$ to $\wti{F}_\g$. By
\eqref{eq:frel}, this section is given by $\iota_\xi
\varepsilon=C(\xi)$, with $C$ given by \eqref{eq:sform}.

Let $[\varepsilon,\varepsilon]_{\on{Sch}}\in
C^\infty(\g_\natural,\wedge^3\g)$ be defined using the Schouten
bracket on $\wedge\g$, and $\d \varepsilon\in
C^\infty(\g_\natural,\wedge^3\g)$ the exterior differential of
$\varepsilon$, viewed as a 2-form on $\g_\natural$.
\begin{proposition}\label{prop:cdybe}
The map $\varepsilon$ satisfies the classical dynamical
Yang-Baxter equation:
\begin{equation}\label{eq:cdybe}
\d \varepsilon+\hh [\varepsilon,\varepsilon]_{\on{Sch}}=\tri.
\end{equation}
\end{proposition}
\begin{proof}
  Proposition \ref{prop:LWX} and the discussion following it
show that $\d \varepsilon+\hh [\varepsilon,\varepsilon]_{\on{Sch}}$
  equals the Courant tensor of $\wti{F}_\g$ (relative to the
  complementary subbundle $E_\g$). By Lemma \ref{lem:upsilon}, together with
  Proposition \ref{prop:diracmaps}, $\Upsilon^{\wti{F}_\g}=\tri$.
\end{proof}
This solution of the classical dynamical Yang-Baxter equation was
obtained in \cite{al:cli}, using a different argument.  As a special
case of Proposition \ref{prop:gaugee}, the map $\varepsilon$ relates
the linear Poisson bivector $\pi_\g$ on $\g\cong \g^*$ with the
pull-back $\exp^*\pi_G\in\mf{X}^2(\g_\natural)$ of the bivector field
\eqref{eq:pig} on $G$:
\[ \exp^*\pi_G= \pi_\g+\A_0(\varepsilon).\]

\subsection{The Gauss-Dirac structure}\label{subsec:gauss}
In this Section we assume that $G=K^\C$ is a complex Lie group, given
as the complexification of a compact, connected Lie group $K$ of rank
$l$.  Thus the Cartan-Dirac structure $E_G$ will be regarded as a
holomorphic Dirac structure on the complex Lie group $G$. We will show
that $G$ carries another interesting Dirac structure besides the
Cartan-Dirac structure. An important feature of this Dirac structure
is that the corresponding Dirac foliation has an open dense leaf.

Take the bilinear form $B$ on $\g$ to be the complexification
of a positive definite invariant inner product on $\k$. Let $T_K$ be a maximal torus
in $K$, with complexification $T=T_K^\C$. Let
\begin{equation}\label{eq:triang}
\g=\n_-\oplus \t \oplus \n_+
\end{equation}
be the triangular decomposition relative to some choice of
positive Weyl chamber, where $\n_+$ (resp. $\n_-$) is the
nilpotent subalgebra given as the sum of positive (resp. negative)
root spaces. For every root $\alpha$, let $e_\alpha$ be a
corresponding root vector, with the normalization
$B(\ol{e}_\alpha,e_\alpha)=1$ and $e_{-\alpha}=\ol{e}_\alpha$. The
unipotent subgroups corresponding to $\n_\pm$ are denoted $N_\pm$. Recall that the
multiplication map
\begin{equation}\label{eq:j}
 j\colon N_-\times T \times N_+\to G,\ (g_-,g_0,g_+)\mapsto
 g_-g_0g_+
\end{equation}
is a diffeomorphism onto its image $\O\subset G$, called the \emph{big
Gauss cell}. The big Gauss cell is open and dense in $G$, and the
inverse map $j^{-1}\colon \O\to N_-\times T \times N_+$ is known as
the Gauss decomposition.  Consider $\dd=\g\oplus \ol{\g}$ as Section
\ref{subsec:isom}. Then
\begin{equation}\label{eq:Gaussh}
\s=\{(\xi_+ +\xi_0)\oplus (\xi_- -\xi_0)\in \dd\, | \, \ \xi_-\in
\n_-,\ \xi_0\in\t,\ \xi_+\in\n_+\}\end{equation}
is a Lagrangian subalgebra of $\dd$, corresponding to the subgroup
\[ S=\{(g_+t,\,g_-t^{-1})\in G\times G\,|\, g_-\in N_-,\ t\in T,\ g_+\in
N_+\}\]
of $D=G\times G$.  Since $\mf{s}$ is transverse to the diagonal
$\g_\Delta$, the corresponding Lagrangian subbundle
$\wh{F}_G:=E^\s$ is transverse to the Cartan-Dirac structure
$E_G$:
\[ \TG=E_G\oplus \wh{F}_G.\]
We shall refer to it as to {\em Gauss-Cartan} splitting.

Unlike the complement $F_G$ defined by the anti-diagonal,
$\wh{F}_G$ is integrable (since $\mf{s}$ is a subalgebra), and it
defines a Dirac manifold $(G,\wh{F}_G,\eta)$. We refer to
$\wh{F}_G$ as the \emph{Gauss-Dirac structure}. Its leaves are the
orbits of $S$ as a subgroup of $D$,
\begin{equation}\label{eq:saction}
\A(g_+t,g_- t)(g)=g_-t^{-1} g t^{-1} g_+^{-1}.\end{equation}
The $S$-orbit of the group unit is exactly the big Gauss cell. Let
$\om_\O$ be the 2-form on $\O$, and $j^*\om_\O$ its pull-back to
$N_-\times T\times N_+$.
\begin{proposition}\label{prop:gauss}
The pull-back of the 2-form $\om_\O$ on the big
  Gauss cell
$N_-\times T\times N_+$
is given
by:
\begin{equation}\label{eq:omformula}
 j^*\om_\O=-\hh B(\theta^L_-,\Ad_{g_0} \theta^R_+).
\end{equation}
Here $\theta^L_\pm,\ \theta^R_\pm$ are the Maurer-Cartan-forms on
$N^\pm$, and  $g_0$ is the $T$-component (i.e. projection of
$N_-\times T\times N_+$ to the middle factor).
\end{proposition}

\begin{proof}
Let $\om\in\Om^2(N_-\times T\times N_+)$ denote the 2-form on the
right hand side of \eqref{eq:omformula}. Since both $\om$ and $\om_\O$
are $S$-invariant, it suffices to check that $j^*\om_\O=\om$ at the
group unit $g=e$.  At the group unit, the formula \eqref{eq:genform}
for $\om_\O$ simplifies to
\begin{equation}\label{eq:simpler}
 \om_\O(\A(\zeta_1),\A(\zeta_2))|_e=
\hh (B(\xi_1',\xi_2)-B(\xi_2',\xi_1)),\end{equation}
for $\zeta_1=(\xi_1,\xi_1'),\ \zeta_2=(\xi_2,\xi_2')\in \mf{s}\subset
\dd$. Its kernel is
\[ \ker(\om_\O|_e)=\{\A(\zeta)\big|_e\, |\ \zeta=(\xi_0,-\xi_0),\xi_0\in\t\}=T_e(T)\]
which coincides with the kernel of $-\hh B(\theta^L_-,\theta^R_+)|_e$.
Moreover, it is clear that $T_e(N_+)$ and $T_e(N_-)$ are isotropic
subspaces for both 2-forms. Hence it is enough to compare on tangent
vectors $\A(\zeta_1),\A(\zeta_2)$ for $\zeta_i$ of the form
$\zeta_1=(0,\xi_-)$ with
$\xi_-\in\n_-$, and $\zeta_2=(\xi_+,0)$ with $\xi_+\in \n_+$.
\eqref{eq:simpler} gives,
\[ \om_\O(\A(0,\xi_-),\A(\xi_+,0))|_e= \hh B(\xi_+,\xi_-).\]
Since $j^*\A(\xi_+,0)|_e=(0,0,\xi_+)\in \n_+ \subset \g=T_eG$ and
$j^*\A(0,\xi_-)|_e=(-\xi_-,0,0)$, the right hand side of
\eqref{eq:omformula} gives exactly the same answer.
\end{proof}

Since $F_G$ and $\wh{F}_G$ are both complements to the Cartan-Dirac
structure $E_G$, they are related by an element in $\Gamma(\wedge^2
E_G)$. To compute this element, let $\mf{p}$ be the anti-diagonal in
$\dd=\g\oplus\ol{\g}$, and let $\g_\Delta\cong\g$ be the diagonal.
Let
\begin{equation}
\mf{r}=\sum e_{-\alpha}\wedge e_{\alpha}\in \wedge^2\g
\end{equation}
be the classical $r$-matrix.
\begin{lemma}\label{lem:rmatrix}
  The bivector taking $\mf{p}$ to $\s$ is the image $\mf{r}_\Delta\in
  \wedge^2\g_\Delta$ of the classical $\mf{r}$-matrix under the
  diagonal embedding $\g\to \g_\Delta\subset\dd$.
\end{lemma}
\begin{proof}
Let $\g\oplus\g^*$ carry the bilinear form defined by the pairing,
and consider the isometric isomorphism
\[ \g\oplus \g^*\to \dd=\g\oplus \ol{\g},\ \ \ \xi\oplus \mu\mapsto
(\xi+{\ts \f{B^\sharp(\mu)}{2}})\oplus
(\xi-{\ts \f{B^\sharp(\mu)}{2}}).\]
This isomorphism takes $\g=\g\oplus 0$ to the diagonal $\g_\Delta$,
and $\g^*$ to the anti-diagonal, $\mf{p}$. The graph
$\on{Gr}_{\mf{r}}\subset \g\oplus\g^*$ of the bivector $\mf{r}$ is
spanned by vectors of the form
\[ 0\oplus B^\flat(\xi_0),\ \ e_{\alpha}\oplus  B^\flat(e_\alpha),\ \
  e_{-\alpha} \oplus (-B^\flat(e_{-\alpha})),\]
for $\xi_0\in\t$ and positive roots $\alpha$. The isomorphism
$\g\oplus \g^*\cong\dd$ takes these vectors to
\[ \xi_0/2\oplus (-\xi_0/2),\ \ 0\oplus e_{-\alpha},\ e_{\alpha}\oplus
0.\]
Hence, it defines an isomorphism $\on{Gr}_{\mf{r}}\cong \s$.
\end{proof}

\begin{corollary}
The orthogonal transformation $A^{-\sfe(\mf{r})}\in\Gamma(\on{O}(\TG))$ takes
$F_G$ to $\wh{F}_G$.
\end{corollary}
\begin{proof}
This follows from Lemma \ref{lem:rmatrix} and the isomorphism
$\TG\cong G\times\dd$.
\end{proof}

The Gauss-Cartan splitting $\TG=E_G\oplus \wh{F}_G$ also defines a bivector
field $\wh{\pi}_G$, and Proposition\ref{prop:gaugee} implies that it is
related to the bivector field $\pi_G$ \eqref{eq:pig} by
\[ \wh{\pi}_G=\pi_G+\A_{{\sf{ad}}}(\mf{r}).\]
Since $\wh{F}_G$ is integrable, this bivector field is in fact a
Poisson structure on $G$ -- see the remarks before Proposition
\ref{prop:diracmaps}. (On the other hand, unlike $\pi_G$, the Poisson
structure is not invariant under the full adjoint action, but is only
$T$-invariant.)

\begin{proposition}
The Poisson structure $\wh{\pi}_G$ associated with the
Gauss-Cartan splitting
$\TG=E_G\oplus \wh{F}_G$ is given by the formula:
\[ \wh{\pi}_G=
\hh \sum_i e_i^L\wedge (e^i)^R -\sum_{\alpha\succ
0}e_{-\alpha}^L\wedge e_\alpha^R +\hh \mf{r}^L+\hh \mf{r}^R.
\]
Here  $e_i$ is a basis of $\t$, with $B$-dual basis $e^i$, and
$\mf{r}^L,\ \mf{r}^R$ are the left-, right-invariant bivector
fields defined by $\mf{r}$. The symplectic leaves of this Poisson
structure are the connected components of the
intersections of conjugacy classes in $G$ with the orbits of the
action \eqref{eq:saction}.
\end{proposition}

This Poisson structure was first defined by Semenov-Tian-Shansky,
see \cite{se:dr}.
\begin{proof}
The vectors
\[ \hh (e_i\oplus(-e_i)),\ 0\oplus (-e_{-\alpha}),\ \ e_\alpha\oplus 0\]
form basis of $\s$ that is dual (relative to the bilinear form on $\dd=\g\oplus
\ol{\g}$) to the basis
\[ e^i\oplus e^i,\ \  e_\alpha\oplus e_\alpha,\ \ e_{-\alpha}\oplus e_{-\alpha}\]
of the diagonal. Using the formula \eqref{eq:pibasis} for the
bivector field, we obtain
\[ \begin{split}
\wh{\pi}_G&=\hh \sum_i ((e^i)^L-(e^i)^R)\wedge \f{e_i^L+e_i^R}{2}
+\hh \sum_{\alpha\succ 0}   (e_\alpha^L-e_\alpha^R)\wedge
(-e_{-\alpha}^L) +\hh \sum_{\alpha\succ 0}
(e_{-\alpha}^L-e_{-\alpha}^R)\wedge(-e_\alpha)^R
\\
&=\hh \sum_i e_i^L\wedge (e^i)^R -\sum_{\alpha\succ 0}e_{-\alpha}^L\wedge e_\alpha^R +\hh \mf{r}^L+\hh \mf{r}^R.
\end{split}\]
Here we have used that the left- and right-invariant bivector fields generated by
\[ \sum_i e_i \wedge e^i=\sum_i e_i\wedge
e^i+\sum_{\alpha\succ 0}e_{-\alpha}\wedge e_\alpha
+\sum_{\alpha\succ 0}e_{\alpha}\wedge e_{-\alpha}\]
coincide.
\end{proof}

\begin{remark} The Lagrangian subalgebra $\s$ defines a Manin triple
  $(\dd=\g\oplus \overline{\g},\g_\Delta,\s)$, which induces a
  Poisson-Lie group structure on the double $D=G\times G$. The Poisson
  structure $\wh{\pi}_G$ is the push-forward image of this Poisson-Lie
  structure under the natural projection $D\to D/G\cong G$, see e.g.
  \cite[Sec.~3.6]{al:ma}.
\end{remark}

\section{Pure spinors on Lie groups}\label{sec:spingroup}
In the previous section we identified $\TG\cong G\times\dd$ as
Courant algebroids. In particular, we have an identification
$\Cl(\TG)\cong G\times \Cl(\dd)$ of Clifford algebra bundles. In
this section, we will complement this isomorphism of Clifford
bundles by an isomorphism of spinor modules,
\[ \wedge T^*G\cong G\times \Cl(\g),\]
where $\Cl(\g)$ is given the structure of a spinor module over
$\Cl(\dd)$.  The differential $\d+\eta$ on $\Om(G)$ intertwines
with a certain differential $\d_\Cl$ on $\Cl(\g)$.  Hence, given a
pure spinor $x\in\Cl(\g)$ defining a Lagrangian subspace
$\s\subset\mf{d}$, one directly obtains a pure spinor
$\phi^\s\in\Om(G)$ defining $E^\s$. We will also obtain
expressions for $(\d+\eta)\phi^\s$ from the properties of $x$.

\subsection{$\Cl(\g)$ as a spinor module over $\Cl(\g\oplus\ol{\g})$}
Recall from Examples \ref{ex:example1} and \ref{ex:example2} that
for any vector space $V$ with inner product $B$, the Clifford
algebra $\Cl(V)$ may be viewed as a spinor module over
$\Cl(V\oplus \ol{V})$. In the special case that $V=\g$ is a Lie
algebra, with $B$ an invariant inner product, there is more
structure that we now discuss.

Let $\lambda\colon\g\to \wedge^2\g$ be the map defined by the
condition $-\iota(\xi_2)\lambda(\xi_1)=[\xi_1,\xi_2]_\g$ (see
Section \ref{subsec:clif}), and let $\tri\in \wedge^3\g$ be the
structure constants tensor \eqref{eq:Xi}. Then
\[\begin{split}
\{\lambda(\xi_1),\lambda(\xi_2)\}&=\lambda([\xi_1,\xi_2]_\g),\ \
\{\lambda(\xi_1),\xi_2\}=[\xi_1,\xi_2]_\g\\
\{\tri,\xi\}&=-\f{1}{4}\lambda(\xi),\ \ \ \ \{\tri,\tri\}=0
\end{split}\]
for all $\xi_1,\xi_2,\xi\in\g$.
The quantizations of these elements have
similar properties: Let
\begin{equation}\label{eq:qlambda}
\tau\colon \g \to \Cl(\g),\ \  \tau(\xi)=q(\lambda(\xi)).
\end{equation}
Then
\[ \begin{split}
[\tau(\xi_1),\tau(\xi_2)]_\Cl&=\tau([\xi_1,\xi_2]_\g),\ \
[\tau(\xi_1),\xi_2]_\Cl=[\xi_1,\xi_2]_\g,\\
 [q(\tri),\xi]_\Cl&=-\f{1}{4}\tau(\xi),\ \ \ \  [q(\tri),q(\tri)]_\Cl\in \mathbb{K}.\end{split}\]
(One can show (cf. \cite{al:no}) that the constant
$[q(\tri),q(\tri)]_\Cl$ is $\f{1}{24}$ times the trace of the Casimir
operator in the adjoint representation.) This last identity implies
that the derivation
\begin{equation}\label{eq:clifdiff}
\d^\Cl=-4[q(\tri),\cdot]_\Cl:\Cl(\g) \to \Cl(\g)
\end{equation}
squares to $0$. We call $\d^\Cl$ the \emph{Clifford differential}
\cite{al:no,ko:sy2}.

For the Lie algebra $\dd=\g\oplus \ol{\g}$, with bilinear form $B\oplus(-B)$, the
corresponding elements $\tri_\dd$ and $\lambda_\dd$ in
$\wedge\dd=\wedge\g\otimes\wedge\g$ are given by
\[ \tri_\dd=\tri\otimes 1+1\otimes\tri,\ \ \ \
\lambda_\dd(\xi,\xi')=\lambda(\xi)\otimes
1-1\otimes\lambda(\xi'),\;\;  \mbox{ for } (\xi,\xi')\in \dd.
\]
Note  also that $q(\tri_\dd)^2=0$.
Consider the Clifford algebra $\Cl(\g)$ as a spinor
module over $\Cl(\dd)$, with Clifford
action given on generators $\zeta=(\xi,\xi')\in \dd$ by
\[ \varrho^\Cl(\xi,\xi')=l^\Cl(\xi)-r^\Cl(\xi').\]
Then the Clifford differential $\d^\Cl$ is implemented as a
Clifford action:
\[ \d^\Cl=-4 \varrho^\Cl(q(\tri_\dd)).\]
The elements $\tau_\dd(\zeta)=q(\lambda_\dd(\zeta))$ generate a
$\dd$-action on $\Cl(\g)$, with generators
\[ \L^\Cl(\zeta)=l^\Cl(\tau(\xi))-r^\Cl(\tau(\xi'))=
\varrho^\Cl(\tau(\zeta)).\]
Note that
\begin{equation}\label{eq:wehave1}
 \L^\Cl(\zeta)=[\varrho^\Cl(\zeta),\d^\Cl],
\end{equation}
which implies that
\[ [\varrho^\Cl(\zeta_1),[\varrho^\Cl(\zeta_2),\d^\Cl]]
=\varrho^\Cl([\zeta_1,\zeta_2]).\]
%
%
%

Let $\s\subset\dd$ be a Lagrangian subspace, and recall the definition
of $\Upsilon^\s$ given in \eqref{eq:Upsh}. Given a Lagrangian
complement $\mf{p}$ to $\s$, let $\pr_\s\colon \dd\to \s$ be the
projection along $\mf{p}$, and define a linear functional
$\sig^\s\in\s^*$ by
\begin{equation}\label{eq:sig}
\l\sig^\s,\ \xi\r=\hh \on{trace}(\pr_\s\circ \ad_\xi\big|_{\s}),\ \ \xi\in\s.
\end{equation}
If $\s$ is a Lagrangian subalgebra (i.e. $\Upsilon^\s=0$), we may
omit $\pr_\s$ in this formula; in this case $\sig^\s$ equals
$-\hh$ times the modular character of the Lie algebra $\s$.
\begin{proposition}\label{prop:integ}
  Let $\s\subset \dd$ be a Lagrangian subspace, with
  defining pure spinor $x\in \Cl(\g)$. Choose a Lagrangian complement
$\mf{p}\cong\s^*$ to $\s$ to view $\Upsilon^\s$ as an element of the Clifford
algebra $\Cl(\dd)$.
Then
\[ \d^\Cl x=\varrho^\Cl(-\Upsilon^\s+\sig^\s)x.\]
In particular,  $\mf{s}$ is a Lie subalgebra if and only if the defining pure spinor $x$ is
`integrable', in the sense that
\[ \d^\Cl\,x\in \varrho^\Cl(\dd)x.\]
\end{proposition}
\begin{proof}
The choice of a Lagrangian complement identifies $\dd=\s\oplus
\s^*$, with bilinear form given by the pairing. Using a basis
$e_i$ of $\s$ and a dual basis $f^i$ of $\s^*$, we have
\[
\begin{split}
4\tri_\dd =&\ts{\f{1}{6}}\sum_{ijk}B_\dd([e_i,e_j],e_k) f^i\wedge f^j\wedge f^k
+\hh \sum_{ijk} B_\dd([e_j,e_k],f^i)\ e_i\wedge f^j\wedge f^k
\\&+ \hh \sum_{ijk} B_\dd([f^j,f^k],e_i)\ e_j\wedge e_k\wedge f^i
+\ts{\f{1}{6}}\sum_{ijk} B_\dd([f^j,f^k],f^i)\ e_j\wedge e_k\wedge e_i.
\end{split}
\]
The quantization map takes the last two terms into the left ideal
$\Cl(\dd)\s$, and it takes the second term to
\[  \hh \sum_{ik} B_\dd([e_i,e_k],f^i)\ f^k+\hh \sum_{ijk}
B_\dd([e_j,e_k],f^i)\ f^j f^k e_i =-\sig^\s \mod \Cl(\dd)\s.
\]
This gives
\[ -4q(\tri_\dd)=-\Upsilon^\s+\sig^\s \mod \Cl(\dd)\s,\]
from which the result is immediate.
\end{proof}
Let us now assume that the adjoint action
$\Ad\colon G\to \on{O}(\g)$ lifts to a group homomorphism
\begin{equation}\label{eq:tau2} \tau\colon G\to \on{Pin}(\g)\subset
  \Cl(\g)\end{equation}
to the double cover $\on{Pin}(\g)\to \on{O}(\g)$. If $G$ is connected,
this is automatic if $\pi_1(G)$ is torsion free. Note that \eqref{eq:tau2} is
consistent with our previous notation $\tau(\xi)=q(\lambda(\xi))$, since \cite{al:no}
\begin{equation*}
\tau(\xi)=\left.\f{d}{d t}\right|_{t=0}\tau(\exp t\xi).
\end{equation*}

We will write $\sfN(g)=\sfN(\tau(g))=\pm 1$
for the image under the norm homomorphism, and
$|g|=|\tau(g)|$ for the parity
of $\tau(g)$. Since $\tau(g)$ lifts $\Ad_g$, one has
$(-1)^{|g|}=\det(\Ad_g)$. The definition of the Pin
group implies that conjugation by $\tau(g)$ is the \emph{twisted
  adjoint action},
\begin{equation}\label{eq:equivq}
\tau(g)x\tau(g^{-1})=\wt{\Ad}_g(x):=(-1)^{|g||x|}\Ad_g(x)
\end{equation}
(using the extension of $\Ad_g\in \on{O}(\g)$ to an automorphism of the Clifford
algebra). This twisted adjoint action extends to an action of the
group $D$ on $\Cl(\g)$,
\begin{equation} \A^\Cl(a,a')(x)=\tau(a)x\tau((a')^{-1}).\end{equation}
%

\subsection{The isomorphism $\wedge T^*G\cong G\times \Cl(\g)$}
Let us now fix a generator $\mu\in\det(\g)$, and consider the
corresponding star operator $\star \colon \wedge\g^* \to
\wedge\g$, see Remark \ref{rem:reference}(\ref{it:rem-b}). The
star operator satisfies
\begin{equation}\label{eq:equivstar}
\Ad_g\circ\ \star =(-1)^{|g|} \star\circ \Ad_{g^{-1}}^*.
\end{equation}
We use the trivialization by left-invariant forms to identify
$\wedge T^*G \cong G\times\wedge\g^*$. Applying $\star$ pointwise,
we obtain an isomorphism
$q\circ \star\colon \wedge T^*_gG
\stackrel{\sim}{\longrightarrow} \Cl(\g)$
for each $g\in G$. Let us define the linear map
\begin{equation}\label{eq:Req}
 \ca{R}\colon \Cl(\g) \to \Om(G),\ \
\ca{R}(x)|_g= (q\circ \star)^{-1}(x\tau(g)).
\end{equation}
We denote by $\mu^*\in\det(\g^*)$ the dual generator, defined by
$\iota((\mu^*)^\top)\mu=1$, and let $\mu_G$ be the left-invariant
volume form on $G$ defined by $\mu^*$.

\begin{proposition}\label{prop:R}
The map \eqref{eq:Req}
has the following properties:
\begin{enumerate}
\item $\ca{R}$ intertwines the Clifford actions, in the sense that
\[
\ca{R}(\varrho^\Cl(\zeta)x)=\varrho(\sfs(\zeta))\ca{R}(x),\;\;\;
\forall x\in \Cl(\g),\, \zeta \in \dd.
\]
Up to a scalar function, $\ca{R}$ is uniquely characterized by
this property.

\item $\ca{R}$ intertwines differentials:
\[
\ca{R}(\d^\Cl(x))=(\d+\eta)\ca{R}(x),\;\;\; \forall x\in
\Cl(\g).
\]
\item $\ca{R}$ satisfies has the following $D$-equivariance property:
      For any
      $h=(a,a')\in D$, and at any given point $g\in G$,
\[ \A(h^{-1})^*\ca{R}(x)=(-1)^{|a|(|g|+|x|)}\ca{R}(\A^\Cl(h)x).\]
\item $\ca{R}$ relates the bilinear pairings on the Clifford modules
$\Cl(\g)$ and $\Om(G)$ as follows: At any given point $g\in G$, and
for all $x,x'\in \Cl(\g)$,
\begin{equation}\label{eq:chevpair}
 (\ca{R}(x),\ca{R}(x'))_{\wedge T^*G}=(-1)^{|g|(\dim G+1)}\ \sfN(g) \
  (x,x')_{\Cl(\g)}\ \  \mu_G.
\end{equation}
Here the pairing $(\cdot,\cdot)_{\Cl(\g)}$ is viewed as scalar-valued, using the
trivialization of $\det(\g)$ defined by $\mu$.  (Cf. Remark \ref{rem:reference}.)
\end{enumerate}
\end{proposition}
Notice that the signs in part (c), (d) disappear if $G$ is connected.
\begin{proof}
  (a) Given $\xi \in \g$, let $\eps(\xi):\wedge\g \to \wedge \g$ be
  defined by $\eps(\xi)\xi'=\xi\wedge \xi'$. Then
\[ l^\Cl(\xi)\circ
  q=q\circ (\eps(\xi)+\hh \iota(B^\flat(\xi))),\ \ r^\Cl(\xi)\circ
  q=q\circ (\eps(\xi) - \hh \iota(B^\flat(\xi))).\]
Since the star
  operator exchanges exterior multiplication and contraction, we have
\[ \star^{-1} \circ q^{-1}\circ \varrho^\Cl(\xi,\xi')=
\Big(\iota(\xi-\xi')+\eps\left(B^\flat\left({\ts
  \f{\xi+\xi'}{2}}\right)\right)\Big)\circ \,\star^{-1} \circ q^{-1}.
\]
On the other hand,
\[
(\varrho^\Cl(\xi,\xi')x)\tau(g)=(\xi x-(-1)^{|x|} x\xi')\tau(g)
= \varrho^\Cl(\xi,\Ad_{g^{-1}}\xi') (x\tau(g)).
\]
This implies that, at $g\in G$,
\[
\ca{R}(\varrho^\Cl(\xi,\xi')x)=
\Big(\iota(\xi-\Ad_{g^{-1}}\xi')+\eps\left(B^\flat\left({\ts
\f{\xi+\Ad_{g^{-1}}\xi'}{2}}\right)\right)\Big)\ca{R}(x),
  \]
which is precisely the Clifford action of $\sfs(\xi,\xi')$ since
\[ \sfs(\xi,\xi')=(\xi-\Ad_{g^{-1}}\xi')\oplus B^\flat\left({\ts \f{\xi+\Ad_{g^{-1}}\xi'}{2}}\right)\]
under left-trivialization $\TG\cong G\times (\g\oplus \g^*)$. This
shows that $\ca{R}$ intertwines the Clifford actions of
$\Cl(\dd)\cong \Cl(\T_gG)$. By the uniqueness properties of spinor
modules, $\ca{R}$ is uniquely characterized by this property up to
a scalar.

(b) From the global equivariance property in (c), verified below, we obtain the
infinitesimal equivariance:
$\ca{L}(\A(\zeta))\ca{R}(x)=\ca{R}(\ca{L}^\Cl(\zeta) x).$
Since $[\varrho(\sfs(\zeta)),\d+\eta]=\ca{L}(\A(\zeta))$ and
$[\varrho^\Cl(\zeta),\d^\Cl]=\L^\Cl(\zeta)$, this gives
\[ \begin{split}
\varrho(\sfs(\zeta))\big((\d+\eta)\ca{R}(x)-\ca{R}(\d^\Cl x)\big)
&=\L(\A(\zeta))\ca{R}(x)-\ca{R}(\varrho^\Cl(\zeta)\d^\Cl x)\\
&=\L(\A(\zeta))\ca{R}(x)-\ca{R}(\L^\Cl(\zeta) x).
\end{split}\]
That is, the map $(\d+\eta)\circ \ca{R}-\ca{R}\circ \d^\Cl
\colon \Cl(\g)\to \Gamma(\TG)$ intertwines the Clifford actions,
and hence agrees with $\ca{R}$ up to a scalar function.  Since its
parity is opposite to that of $\ca{R}$, that function is zero.

(c) We have to show that for all $a\in G$,
\begin{equation}\label{eq:alteq}
l_a^*\ca{R}(x)=\ca{R}(x\tau(a)) ,\ \ \ r_{a}^*\ca{R}(x)= (-1)^{|a|(|g|+|x|)}\ca{R}(\tau(a)x).
\end{equation}
In terms of the left-trivialization $\wedge T^*G=G\times\wedge\g^*$,
\[ (l_a^*\ca{R}(x))|_g=\ca{R}(x)\big|_{ag},\ \ (r_a^*\ca{R}(x))\big|_g=\Ad_{a^{-1}}^*(\ca{R}(x)\big|_{ga}).\]
(Here $\Ad_{a^{-1}}^*$ stands for the contragredient action on $\wedge\g^*$, not for a
  pull-back on $\Om(G)$.) We compute, using \eqref{eq:equivq} and \eqref{eq:equivstar}:
\[
\begin{split}\Ad_{a^{-1}}^*\big(\ca{R}(x)\big|_{ga}\big)
&= (-1)^{|a|} \star^{-1} \, q^{-1} \Ad_{a}(x\tau(g a))\\
&=(-1)^{|a|} \,(-1)^{|a|(|x|+|g|+|a|)} \star^{-1} \, q^{-1} (\tau(a)x\tau(g))\\
&=    (-1)^{|a|(|x|+|g|)} \ca{R}(\tau(a)x)\big|_g
\end{split}
\]
The equivariance property with respect to left translations is
immediate from the definition.

(d) Use the generator $\mu\in\det(\g)$ and $\mu_G$ to trivialize both
$\det(\g)$ and $\det(\wedge T^*G)$. By Remark
\ref{rem:reference}(\ref{it:rem-b}) and Example \ref{ex:example2}
we have, at $g\in G$,
\[ (\ca{R}(x),\ca{R}(x'))_{\wedge T^*G}=(x\tau(g),x'\tau(g))_{\Cl(\g)}.\]
This is computed as follows:
\[\begin{split}
 \on{str}(\tau(g)^\top x^\top x'\tau(g))
&=(-1)^{|g|(|g|+|x|+|x'|)}\on{str}(\tau(g)\tau(g)^\top x^\top x')
\\&=\sfN(g)\,(-1)^{|g|(1+|x|+|x'|)}\on{str}(x^\top x')
\end{split}\]
Finally, replace $|x|+|x'|$ with $\dim G$, using that $(x,x')_{\Cl(\g)}$ vanishes
unless $ |x|+|x'|=\dim G\mod 2$.
\end{proof}

As an immediate consequence of Propositions~\ref{prop:integ}
and~\ref{prop:R}, we have

\begin{corollary}
 If $x\in
\Cl(\g)$ is a pure spinor defining a Lagrangian subspace $\s
\subset \dd$, then the differential form
$\phi^{\mf{s}}:=\ca{R}(x)\in \Om(G)$ is a pure spinor defining the
Lagrangian subbundle $E^\s$. It satisfies the differential
equation
\begin{equation}\label{eq:diffeq}
(\d+\eta)\phi^{\mf{s}}=\varrho(\sfs(-\Upsilon^\s+\sig^\s))\phi^{\mf{s}},
\end{equation}
where $\sig^\s\in\s^*$ is defined as in \eqref{eq:sig} (using a
complementary Lagrangian subspace $\mf{p}\cong \s^*\subset \dd$).
Let $H\subset D$ be a subgroup preserving $\mf{s}$, and define the
character $u^\s\colon H\to \mathbb{K}^\times$ by
$\A^\Cl(h)x=u^\s(h)x$. Then
\begin{equation}\label{eq:eqqq}
\A(h^{-1})^*\phi^{\mf{s}}=(-1)^{|a|(|g|+|x|)}\,u^\s(h)\phi^{\mf{s}}\end{equation}
for all $h=(a,a')\in H$, and at any given point $g\in G$.
\end{corollary}

We are mainly interested in pure spinors defining the Cartan-Dirac
structure $E_G$ and its Lagrangian complement $F_G$. These are
obtained by taking $x=1$ and $x=q(\mu)$ in the above:
\begin{proposition}\label{prop:phipsi}
Let $\phi_G,\psi_G\in \Om(G)$ be the differential forms
\begin{equation}
 \phi_G=\ca{R}(1),\ \ \psi_G=\ca{R}(q(\mu)).\end{equation}
Then $\phi_G,\psi_G$ are pure spinors defining the Lagrangian
subbundles $E_G,F_G$. They satisfy the differential equations,
\begin{equation}\label{eq:psieqn}
(\d+\eta)\phi_G=0,\ \ (\d+\eta)\psi_G=-\varrho(\sfe(\tri))\psi_G.
\end{equation}
%
The equivariance properties under the adjoint action of $G$ read
\[ \A_{{\sf{ad}}}(a^{-1})^*\phi_G=(-1)^{|a||g|}\phi_G,\ \
   \A_{{\sf{ad}}}(a^{-1})^*\psi_G=(-1)^{|a|(|g|+1)}\psi_G.
\]
\end{proposition}
We will refer to $\phi_G$ as the \emph{Cartan-Dirac spinor}.
\begin{proof}
It is clear that the diagonal $\g_\Delta\subset \dd$ is defined by the pure spinor $x=1$.
Similarly, the anti-diagonal $\mf{p}\subset \dd=\g\oplus\ol{\g}$ is defined by
the pure spinor $q(\mu)\in \Cl(\g)$:
\[\varrho^\Cl(\xi,-\xi)q(\mu)=\xi q(\mu)+
(-1)^{\dim G}q(\mu)\xi=0.\] Hence $\phi_G,\psi_G$ are pure spinors
defining $E_G,F_G$.  The equivariance properties are special cases
of \eqref{eq:eqqq}, since both $\g_\Delta$ and $\mf{p}$ are
preserved under $G_\Delta$.  Here we are using $|1|=0,\
|q(\mu)|=\dim G\mod 2$, while $u^{\mf{p}}(a)=(-1)^{|a|(1+\dim G)}$
by the calculation:
\[\tau(a)q(\mu)\tau(a^{-1})=(-1)^{|a|\dim G}q(\Ad_a(\mu))
=(-1)^{|a|(1+\dim G)}q(\mu).\]
The differential equation for $\phi_G$ follows since $\d^\Cl(1)=0$.
It remains to check the differential equation for $\psi_G$.  Since the
anti-diagonal satisfies $[\mf{p},\mf{p}]_\dd\subset \g_\Delta$, the
element $\sigma^{\mf{p}}\in\mf{p}^*$ is just zero. On the other hand,
the element $\Upsilon^\mf{p}$ is given by $\tri_\Delta$, the image of
$\tri$ under the the map $\wedge\g \stackrel{\sim}{\to}
\wedge\g_\Delta$. Hence $\sfs(\tri_\Delta)=\sfe(\tri)$, confirming
that $\psi_G$ satisfies \eqref{eq:psieqn}.
\end{proof}

\begin{remarks}
\begin{enumerate}
\item
The map $\ca{R}$ depends on the choice of generator $\mu\in
\det(\g)$, via the star operator: Replacing $\mu$ with $t\mu$ changes
$\ca{R}$ to $t^{-1}\ca{R}$. Hence, the definition of
$\psi_G=\ca{R}(q(\mu))$ is independent of the choice of $\mu$.
\item
Since $(1,q(\mu))_{\Cl(\g)}=\mu$, the
 bilinear pairing between $\phi_G,\psi_G$ equals the volume form,
up to a sign:
\[
\big(\phi_G,\psi_G\big)_{\wedge T^*G}=\sfN(g) (-1)^{|g|(\dim
G+1)}\,\mu_G.
\]
\end{enumerate}
\end{remarks}

\begin{proposition}
  Over the open subset $\ca{U}$ of $G$ where $1+\Ad_g$ is invertible,
  the pure spinor $\psi_G$ is given by the formula:
\[ \psi_G={\det}^{1/2}\big({\ts \f{1+\Ad_g}{2}}\big) \exp\Big({\ts \f{1}{4}\,
B\Big(\f{1-\Ad_g}{1+\Ad_g}\ \theta^L,\theta^L\Big)}\Big),
\]
at any given point $g\in \ca{U}$. ( The square root depends on the
choice of lift $\tau\colon G\to \on{Pin}(\g)$.)
\end{proposition}
Note that the exponent in this formula becomes singular where
$1+\Ad_g$ fails to be invertible, but these singularities are
compensated by the zeroes of the factor ${\det}^{1/2}\big({\ts
  \f{1+\Ad_g}{2}}\big)$.  One proof of this formula is given in
\cite{me:lec}; here is an outline of an alternative approach.
\begin{proof}[Sketch of proof]
   One easily checks that over $\ca{U}$,
  $F_G$ coincides with the graph of the 2-form
$\omega_F:=-{\ts \f{1}{4}\, B\Big(\f{1-\Ad_g}{1+\Ad_g}\
    \theta^L,\theta^L\Big)}$.
Hence
  $\psi_G|_{\ca{U}}=f\exp(-\om_F)$ for some nonvanishing function
  $f\in C^\infty(\mathcal{U})$, with $f(e)=1$. Equation
  \eqref{eq:psieqn} reads, after dividing by $f\exp(-\om_F)$,
\[ \d \log(f) +\eta +\exp(\om_F)\varrho(\sfe(\tri))\big(\exp(-\om_F)\big)=0.\]
Taking the form degree  $1$ parts of both sides of this equation, one
obtains the following condition on $f$:
\[ \d \log(f) + \Big(\exp(\om_F) \varrho(\sfe(\tri))\big(\exp(-\om_F)\big)\Big)_{[1]}=0.\]
$f$ is uniquely determined by this Equation with the initial condition
$f(e)=1$. It is straightforward (though slightly cumbersome) to
verify that $f(g)= {\det}^{1/2}\big({\ts
  \f{1+\Ad_g}{2}}\big)$ solves this equation.\end{proof}

%
%
%

%

If $G$ is connected, one has $\det(1+\Ad_g)\not=0$ on a dense open
subset of $G$. However, for a disconnected group $G$ it vanishes on
the components with $\det(\Ad_g)=-1$.

\begin{example}
  Let $G=\on{O}(2)$. Here $\on{O}(\g)=\Z_2$ and $\on{Pin}(\g)=\Z_4$.
  There are two possible lifts $\on{O}(\g)\to \on{Pin}(\g)$. Let
  $\theta\in \Om^1(G)$ be the left-invariant Maurer-Cartan-form (using
  the isomorphism $\g=\R$ defined by a generator $\mu\in \det(\g)=\g$). One finds that on
  $\SO(2)\subset \on{O}(2)$, $\phi_G=\theta$, while $\psi_G=1$. On
  the non-identity component $\on{O}(2)\backslash\SO(2)$, the roles are reversed:
  $\psi_G=\pm\theta$ and $\phi_G=\pm1$. (The signs depend on the
  choice of lift.) Observe that $\phi_G,\psi_G$ given by these
  formulas have the correct equivariance properties.
\end{example}

\subsection{Group multiplication}
In this section, we will examine the composition of the map
$\ca{R}\colon \Cl(\g)\to \Om(G)$ with the pull-back under group
multiplication. It will be convenient to work with the element
$\Lambda\in \Cl(\g)\otimes\Om(G)$, defined by the property
\[ \ca{R}(x)=\on{str}(x\Lambda)\]
where we have extended $\on{str}\colon \Cl(\g)\to \wedge^{[\mathrm{top}]}(\g)=\mathbb{K}$ to the
tensor product with $\Om(G)$. The properties of $\ca{R}$ under the
Clifford action translate into
\[ (l^\Cl(\xi)+\varrho(s^R(\xi)))\Lambda=0,\ \ \
(-r^\Cl(\xi)+\varrho(s^L(\xi)))\Lambda=0.\]
Thus $\Lambda$ is itself a pure spinor for the action of
$\Cl(\dd)\times \Cl(\TG)$, defining a
Lagrangian subbundle of $\dd\times\TG$.
The equivariance properties \eqref{eq:alteq}
of $\ca{R}$ translate into
\[ l_a^* \Lambda=\tau(a^{-1})\Lambda,\ \ \ \ r_{a^{-1}}^*\Lambda=\Lambda \tau(a)\]
The first identity is immediate, while for the second identity is
obtained by the calculation:
\[ \begin{split}
\on{str}(x\,r_{a^{-1}}^*\Lambda)&=r_{a^{-1}}^*\ca{R}(x)=(-1)^{|a|(|x|+|g|)}\ca{R}(\tau(a)x)
\\&=(-1)^{|a|(|x|+|g|)} \on{str}(\tau(a)x\Lambda)
\\&=\on{str}(x\Lambda\tau(a)).\end{split}\]
(Note that $|\Lambda|=|g|$ at $g\in G$.) We finally observe that the
pull-back of $\Lambda$ to the group unit
is simply
\begin{equation}\label{eq:pulle}
i_e^*\Lambda=1\in\Cl(\g).
\end{equation}

Let $\Lambda^1,\Lambda^2\in \Cl(\g)\otimes\Om(G\times G)$ be the
pull-back to the first, second $G$-factor, and recall the 2-form
$\varsigma\in \Om^2(G\times G)$ from \eqref{eq:varsigma}.
\begin{proposition}\label{prop:mult}
The pull-back of $\Lambda$ under group multiplication satisfies
\begin{equation}\label{eq:Lrelated}
e^{\varsigma}\on{Mult}^*\Lambda=\Lambda^1\Lambda^2,\end{equation}
using the product in the algebra $\Cl(\g)\otimes\Om(G\times G)$.
\end{proposition}
\begin{proof}
Using \eqref{eq:rrelated}, we find that both sides of \eqref{eq:Lrelated} are annihilated by
the following operators:
\[ l^\Cl(\xi)+\varrho(s^{R,1}(\xi)),\ \
-r^\Cl(\xi)+\varrho(s^{L,2}(\xi)),\ \ \ \
\varrho(s^{L,1}(\xi)+s^{R,2}(\xi)).\]
Hence the two sides of \eqref{eq:Lrelated} are pure spinors,
defining the same Lagrangian subbundle of $\dd\times \T(G\times
G)$. So the two sides agree up to a scalar function.

  The 2-form $\varsigma$ is
invariant under $l_{a,1}$ (left multiplication by $a$ on the first
factor) and $r_{a^{-1},2}$ (right multiplication by $a^{-1}$ on
the second factor). From the equivariance of $\Lambda$, and since
$\on{Mult}\circ l_{a,1}=l_a\circ \on{Mult}$ and $\on{Mult}\circ
r_{a^{-1},2}=r_{a^{-1}}\circ \on{Mult}$,  we obtain the following
equivariance property of $e^{\varsigma}\on{Mult}^*\Lambda$:
\[ \begin{split}
(l_{a,1})^*(e^{\varsigma}\on{Mult}^*\Lambda)
&=\tau(a^{-1})\ (e^{\varsigma}\on{Mult}^*\Lambda),\\
(r_{a^{-1},2})^*(e^{\varsigma}\on{Mult}^*\Lambda) &=
(e^{\varsigma}\on{Mult}^*\Lambda)\tau(a).\end{split}
\]
The product $\Lambda^1\Lambda^2$ has a similar equivariance property. Hence, to
verify \eqref{eq:Lrelated} it suffices to compare the two sides at
$(e,e)\in G\times G$. But by \eqref{eq:pulle},
both sides of
\eqref{eq:Lrelated} pull back to $1$ at $(e,e)$.
\end{proof}

We will use Proposition \ref{prop:mult} to obtain a formula for
the pull-back of $\psi_G=\ca{R}(q(\mu))$, the pure spinor defining
the Lagrangian subbundle $F_G\subset \TG$. Recall the element
$\gamma\in\wedge^2(\g\oplus\g)$ from \eqref{eq:gamma}.
\begin{theorem}\label{th:psipull}
The pull-back of $\psi_G$ under group multiplication is given by
the formula
\[e^{\varsigma} \, \Mult^*\psi_G
=\varrho(\exp(-{\sfe(\gamma)}))\,(\psi_G^1\otimes \psi_G^2)\]
\end{theorem}
Note that \emph{up to a scalar function}, this identity follows from
Proposition \ref{prop:egamma}.
\begin{proof}
The element $\gamma$ enters the following formula (cf. \cite[Lemma
3.1]{al:no}) , relating the product $\on{Mult}^\Cl$ in $\Cl(\g)$ with
the wedge product $\on{Mult}^\wedge$ in $\wedge(\g)$:
\[ q^{-1}\circ \on{Mult}^\Cl=\on{Mult}^\wedge
\circ \exp(-\iota^\wedge(\gamma))\circ q^{-1}\colon
\Cl(\g\oplus\g)\to \wedge(\g).
\]
Since $\on{str} \circ\, l^\Cl(q(\mu))\circ q\colon \wedge\g\to \mathbb{K}$
is simply the augmentation map, we have

\[\psi_G=\ca{R}(q(\mu))=\on{str}(q(\mu)\Lambda)=
q^{-1}(\Lambda)_{[0]},\] where the subscript indicates the degree
$0$ part with respect to $\wedge\g$. Using \eqref{eq:Lrelated}, we
calculate:
\[\begin{split}
e^{\varsigma}\Mult^*\psi_G&=q^{-1}(\Lambda^1\Lambda^2)_{[0]}\\
&=q^{-1}\circ \big(\on{Mult}^\Cl(\Lambda^1\otimes \Lambda^2)\big)_{[0]}
\\&=\big( \on{Mult}^\wedge\circ \exp(-\iota^\wedge(\gamma))\circ q^{-1}
(\Lambda^1\otimes \Lambda^2)\big)_{[0]}
\\&=\exp(-\sfe(\gamma))\circ \big( \on{Mult}^\wedge\circ  q^{-1}
(\Lambda^1\otimes \Lambda^2)\big)_{[0]}\\
&=\exp(-\sfe(\gamma))\circ (\psi^1_G\otimes \psi^2_G).\end{split}
\]
Here we used that $(\iota^\Cl(\xi)+\varrho(\sfe(\xi)))\Lambda=0$, hence
$(\iota^\wedge(\gamma)-\varrho(\sfe(\gamma)))q^{-1}(\Lambda^1\otimes \Lambda^2)=0$.
\end{proof}

%
\subsection{Exponential map}\label{subsec:exp2}
Let us return to our description (Section \ref{subsec:exp1})
$\T\g^*=\g^*\times\dd_0$ of the Courant algebroid over $\g^*$,
where $\dd_0=\g^*\rtimes \g$.

Let $\wedge\g^*$ be the contravariant spinor module over $\Cl(\dd_0)$
(cf. Section \ref{subsec:contra}), with Clifford action denoted
$\varrho^\wedge$. Let $\d^\wedge$ be the exterior algebra
differential. For all $w=(\beta,\xi)\in \dd_0$ one has
\[ L^\wedge(w):=[\d^\wedge,\varrho^\wedge(w)]=\d^\wedge
\beta-(\ad_\xi)^*.\]
One easily checks that $L^\wedge(w)$ defines an action of the Lie
algebra $\dd_0$. This action exponentiates to an action of the group
$D_0$, given as
\[ \A^\wedge(\beta,g) y=\exp(\d^\wedge\beta)\wedge (\Ad_{g^{-1}})^* y,\]
The function
\[ \tau_0\colon \g^*\to \wedge\g^*,\ \
\tau_0(\beta)=\exp(\d^\wedge\beta)\in\wedge\g^*\]
is the counterpart to the function $\tau\colon G\to \Cl(\g)$. The
$D_0$-action commutes with the differential, and it is straightforward
to check that the Clifford action is
$D_0$-equivariant:
\[ \A^\wedge(\beta,g)\big(\varrho^\wedge(w)y\big)=
\varrho^\wedge(\Ad_{(\beta,g)}w)\big( \A^\wedge(\beta,g) y\big), \]
for $w\in \dd_0,\ (\beta,g)\in D_0,\ y\in \wedge\g^*$.
Choose a generator $\mu\in \det(\g^*)$, and let
$\star\colon\wedge\g\to \wedge\g^*$ be the associated star operator
\footnote{Note that in the previous Section, $\mu$ denoted a
  generator of $\det(\g)$, and hence the star operator went from
  $\wedge\g^*\to \wedge\g$. This change in notation is intended, since
  our aim is to compare the Poisson manifold $\g^*$ with the Dirac
  manifold $G$.}.  Let $X_\pi$ denote the modular vector field of the
Kirillov-Poisson structure $\pi_{\g^*}$, relative to the
translation-invariant volume form $\mu_{\g^*}\in
\Gamma(\det(T^*\g^*))$ defined by the dual generator
$\mu^*\in\det(\g)$.  (Recall that $X_\pi=0$ if $\g$ is unimodular.)
Define a linear map
\[ \ca{R}_0\colon \wedge\g^*\to \Om(\g^*),\]
given at any point $\nu\in\g^*$ by
\[ \ca{R}_0(y)= \star^{-1}  (y\wedge \tau_0(\nu))
\in \wedge\g=\wedge T^*_\nu\g^*
.\]

Parallel to Proposition \ref{prop:R}, we have,
\begin{proposition}\label{prop:R0}
\begin{enumerate}
\item
The map $\ca{R}_0$ intertwines the Clifford actions of $\dd_0$:
 \[\ca{R}_0\circ \varrho^\wedge(w)=\varrho(\sfs_0(w))\circ \ca{R}_0,\ \ w\in
 \dd_0.\]
It is uniquely determined by this property, up to a scalar function.
\item
The map $\ca{R}_0$ intertwines the differentials, up to contraction by
the modular vector field:
\[
 \ca{R}_0\circ \d^\wedge =(\d-\iota(X_\pi))\circ \ca{R}_0.\]
\item
$\ca{R}_0$ has the equivariance property, for all $h=(\beta,a)\in
D_0=\g^*\rtimes G$,
\[
(\A_0(h^{-1}))^*\ca{R}_0(y)=\det(\Ad_a)\ \ca{R}_0(\A^\wedge(h)y).
\]
\item
%
$\ca{R}_0$ preserves the bilinear pairings on the spinor modules
$\wedge\g^*,\ \Om(\g^*)$,  in the sense that
\[(\ca{R}_0(y),\ca{R}_0(y'))_{\wedge T^*\g^*}=(y,\,y')_{\wedge\g^*}\ \mu_{\g^*}\]
for all $y,y'\in \wedge\g^*$.
\end{enumerate}
\end{proposition}
\begin{proof}
  Each of the statements (a),(c),(d) is proved by a direct computation,
  parallel to those in Proposition \ref{prop:R}.
To prove (b), we first note that (c)
  implies the infinitesimal equivariance, for $(\beta,\xi)\in \dd_0$,
\begin{equation}\label{eq:Lieder}
\big(\L(\A_0(\beta,\xi))-\on{tr}(\ad_\xi)\big)\ca{R}_0(y)=
\ca{R}_0(\L_0^\wedge(\beta,\xi)y).\end{equation}
Since $\iota(X_\pi)\l\theta_0,\xi\r=\on{tr}(\ad_\xi)$,
we have
\[\L(\A_0(\beta,\xi))-\on{tr}(\ad_\xi) =[(\d-\iota(X_\pi)),\varrho(\sfs_0(\beta,\xi))].\]
Hence we can re-write \eqref{eq:Lieder} as
\[ [(\d-\iota(X_\pi)),\varrho(\sfs_0(\beta,\xi))]\ \ca{R}_0(y)
=\ca{R}_0\Big([\d^\wedge,\varrho^\wedge(\beta,\xi)]\Big).\]
Together with (a), this implies that the linear map
\begin{equation}\label{eq:themap}
 (\d-\iota(X_\pi))\circ \ca{R}_0-\ca{R}_0\circ \d^\wedge \colon
 \wedge\g^*\to \Om(\g^*)
\end{equation}
intertwines the Clifford actions of $\dd_0$. Since the parity of this
map is opposite to that of $\ca{R}_0$, the uniqueness assertion in (a)
implies that \eqref{eq:themap} is zero.
\end{proof}

As before, we may use this map to construct pure spinors
$\ca{R}_0(y)\in\Om(\g^*)$ from pure spinors $y\in \wedge\g^*$.

The element $y=1$ is the pure spinor defining the Lagrangian
subspace $\g\subset \dd_0$, and its image
$\phi_{\g^*}=\ca{R}_0(1)$ defines the Lagrangian subbundle
$E_{\g^*}$ (spanned by the sections $\sfe_0(\xi)$). The pure
spinor $y=\mu\in\wedge\g^*$ defines a Lagrangian complement
$\g^*\subset \dd_0$, and its image $\psi_{\g^*}=\ca{R}_0(\mu)=1$
defines the Lagrangian subbundle $F_{\g^*}=T\g^*$ (spanned by the
sections $\sff_0(\beta)$). For the bilinear pairing between these
pure spinors, we obtain
\[ (\phi_{\g^*},\psi_{\g^*})_{\wedge T^*\g^*}=\mu_{\g^*}.\]
since $(1,\mu)_{\wedge\g^*}=\mu$.
\begin{lemma}
The pure spinor $\phi_{\g^*}$ is given by the formula
\[ \phi_{\g^*}=(-1)^{n(n-1)/2} e^{-\iota(\pi_{\g^*})}\mu_{\g^*}\]
where $n=\dim G$.
\end{lemma}
\begin{proof}
The Kirillov-Poisson bivector on $\g^*$ is given by
$\pi_{\g^*}|_\nu=-\d^\wedge\nu \in \wedge^2\g^*=\wedge^2 T_\nu\g^*$.
That is, $\tau_0=\exp(-\pi_{\g^*})$. The Lemma follows since
$\star$ intertwines exterior product with contractions, and since
$\star^{-1}(1)=(\mu^*)^\top=(-1)^{n(n-1)/2}\mu^*$.
\end{proof}

Let us now return to our original setting where $\g$ carries an
invariant inner product $B$, used to identify $\g\cong \g^*$. We
take the generators $\mu \in \det(\g)$ (from the last section) and
$\mu\in \det(\g^*)$ (from the present section) to be equal under
this identification.

Let $\mu_\g$ be the translation invariant volume form on
$\g\cong\g^*$, and $\mu_G$ the corresponding left-invariant volume
form on $G$. Let $J\in C^\infty(\g)$ be the Jacobian of the
exponential map, defined by $\exp^*\mu_G=J\,\mu_\g$. Recall
that $\g_\natural\subset \g$ is the dense open subset where $\exp$ is
a local diffeomorphism, i.e where $J\not=0$.
With $\varpi\in\Om^2(\g)$ as in Section \ref{subsec:exp1}, we have:

\begin{proposition}Over the subset $\g_\natural$,
the maps $\ca{R}_0\colon \wedge\g\to \Om(\g)$ and $\ca{R}\colon \Cl(\g)\to
\Om(G)$ are related as follows:
\begin{equation}\label{eq:spinorre}
\exp^*(\ca{R}(x))= J^{1/2} e^{-\varpi}
\varrho(\wt{A}^{-\sfe_0(\varepsilon)})(\ca{R}_0(y)),
\end{equation}
for $x=q(y)$.  Here $\varepsilon\in
C^\infty(\g_\natural,\wedge^2\g)$ is the solution of the classical
dynamical Yang-Baxter equation, cf. Proposition \ref{prop:cdybe},
and $J^{1/2}\in C^\infty(\g)$ is a smooth square root of $J$,
equal to $1$ at the origin.
\end{proposition}

\begin{proof}
The map $\wt{\ca{R}}_0\colon \wedge\g\to \Om(\g_\natural)$ given
as
\[\wt{\ca{R}}_0(y)=  e^{-\varpi} \varrho(\wt{A}^{-\sfe_0(\varepsilon)})\exp^*\ca{R}(q(y))\]
intertwines the $\Cl(\dd_0)$-actions, hence it coincides with
$\wt{\ca{R}}_0=f\,\ca{R}_0$ for a scalar function. To find $f$, we
consider bilinear pairings. Note that
\[\begin{split}
 (\wt{\ca{R}}_0(y),\ti{\ca{R}}_0(y'))_{\wedge T^*\g}
&=
(\exp^*\ca{R}(q(y)),\exp^*\ca{R}(q(y')))_{\wedge T^*\g}\\
&= \exp^* \big(
\ca{R}(q(y)),\,\ca{R}(q(y'))\big)_{\wedge T^*G}.
\end{split}\]
Taking $y'=1,\ y=\mu$ we obtain

\[ f^2\mu_\g=f^2\,({\ca{R}}_0(\mu),\,{\ca{R}}_0(1))_{\wedge T^*\g}=
( \wt{\ca{R}}_0(\mu),\wt{\ca{R}}_0(1))_{\wedge T^*\g}=\exp^*
\mu_G=J\,\mu_\g.\]
This shows that $f^2=J$.
\end{proof}
\begin{remark}
  Of course, $\exp^*(R(x))$ is defined globally on all of $\g$, not
  only on $\g_\natural$. It follows from the Proposition that
  $J^{1/2}\exp(\sfe_0(\varepsilon))$ extends smoothly to all of $\g$.
  Hence, the expression
\[J^{1/2}\exp(\varepsilon)\]
extends smoothly to a global function $\g\to \wedge\g$. For a direct
proof, see \cite{al:cli}.
\end{remark}

Applying the proposition to $y=1$ and $y=\mu$, we find in particular
that
\begin{equation}\label{eq:psiexp}
\begin{split}
 \exp^*\phi_G&=J^{1/2}\,e^{-\varpi} \phi_\g,\\
\exp^*\psi_G&=J^{1/2}\,e^{-\varpi}
\varrho(\wt{A}^{-\sfe_0(\varepsilon)})(1).
\end{split}\end{equation}

\subsection{The Gauss-Dirac spinor}\label{subsec:gdspinor}
We return to the set-up of Section \ref{subsec:gauss}, with
$G=K^\C$ denoting the complexification of a compact Lie group,
with Cartan subgroup $T=T_K^\C$.
 Recall that the Gauss-Dirac structure $\wh{F}_G$ is defined by the Lagrangian
subspace $\s\subset \dd$, with basis the collection of all
$e_\alpha\oplus 0$, $0\oplus e_{-\alpha}$, $e_i\oplus (-e_i)$ where
$\alpha\succ 0$ are positive roots and $i=1,\ldots,l=\on{rank}(G)$.
The element
\begin{equation}\label{eq:gdspinor}
 x=\prod_{\alpha\succ 0}e_{\alpha}e_{-\alpha}\ \prod_i e_i\in \Cl(\g)\end{equation}
is non-zero and is annihilated by the Clifford action of $\s$; hence
it is a pure spinor defining $\s$.
Note that $x$ satisfies
\[ \tau(h_+)x=x,\ \ x\tau(h_-^{-1})=x,\ \ \tau(h_0)x\tau(h_0)=h_0^{2\rho}\,x
\]
for all $h_+\in N_+$, $h_-\in N_-$, $h_0\in T$. Here
$\rho=\hh \sum_{\alpha\succ 0}\alpha$, and
$t\mapsto t^{2\rho}\in\C^\times$ is the character of $T$ defined by the
weight $2\rho$. Hence,
\[ \wh{\psi}_G=\ca{R}(x)\in\Om(G)\]
is a pure spinor defining $\wh{F}_G$. We refer to $\wh{\psi}_G$ as the
\emph{Gauss-Dirac spinor}. Its equivariance
properties are:
\[ l_{h_+}^*\wh{\psi}_G=\wh{\psi}_G,\ \ r_{h_-^{-1}}^*\wh{\psi}_G=
\wh{\psi}_G,\ \ l_{h_0}^*
r_{h_0}^*\wh{\psi}_G=h_0^{2\rho}\wh{\psi}_G.\]
That is, $\wh{\psi}_G$ is invariant up to the character, given by
the group homomorphism $S\to T$ followed by the $2\rho$-character.

Since the big Gauss cell $\O=N_-T N_+\subset G$ is dense in
$G$, the equivariance property, together with the fact that the
pull-back of $\psi_G$ to the group unit is equal to
$\on{str}(x)=1$, completely characterizes the pure spinor
$\psi_G$, and allows us to give an explicit formula. Recall the 2-form
$\om_\O$ on the big Gauss cell, given by \eqref{eq:omformula}:
\begin{proposition}\label{prop:gaussexplicit}
  The restriction of the pure spinor $\wh{\psi}_G$ to the big Gauss
  cell
$\O=j(N_-\times T \times N_+)$ is
  given by the formula,
\[ \wh{\psi}_G\big|_\O=g_0^\rho \exp(-\om_\O).\]
Here $g_0\colon \O\to T$ is the composition of the Gauss decomposition
$j^{-1}\colon \O\to N_-\times T\times N_+$ with projection to the
middle factor.
\end{proposition}
\begin{proof}
  Both sides are pure spinors defining the Gauss-Dirac structure over
  $\O$, with the same equivariance property under $S$, and both sides
  pull back to $1$ at the group unit $e$.
\end{proof}
We now compare the Gauss-Dirac spinor $\wh{\psi}_G$ with the pure spinor
$\psi_G$ from Proposition \ref{prop:phipsi}.
\begin{proposition}
The pure spinors $\psi_G,\wh{\psi}_G$ are related by a twist by the
$r$-matrix $\mf{r}$:
\[ \wh{\psi}_G=\varrho(\exp(-\sfe(\mf{r}))\psi_G. \]
\end{proposition}
\begin{proof}
Let $\mf{r}_\Delta\in \wedge^2\dd$ be the image of $\mf{r}$ under the
diagonal inclusion $\g\hra \dd$. We will show that
\begin{equation}\label{eq:spinoridentity}
x=\varrho^\Cl(\exp(-\mf{r}_\Delta))q(\mu).
\end{equation}
The proposition follows from this identity by applying the map
$\ca{R}$.  Up to a scalar, \eqref{eq:spinoridentity} holds since
both sides are pure spinors defining the same Lagrangian subspace.
To determine the scalar, we apply the super-trace to both sides.
Recall that the spinor action of elements $\xi_\Delta\in
\g_\Delta\subset \g$ is given by Clifford commutator with the
corresponding element $\xi\in\g$. Since the super-trace vanishes
on Clifford commutators, it follows that
\[ \on{str}(\varrho^\Cl(\exp(-\mf{r}))q(\mu))=\on{str}(q(\mu))=1
=\on{str}(x).\]
\end{proof}
Let us next compute the Clifford differential
$\d^\Cl=-4[q(\tri),\cdot]$ of the element \eqref{eq:gdspinor}. Let
$\rho=\hh \sum_{\alpha\succ 0}\alpha\in\t^*$ be the half-sum of
positive (real) roots.
\begin{lemma}The quantization of the structure constant tensor satisfies,
\[ -4q(\tri)=\tpi \rho \mod \n_-\Cl(\g)\n_+.\]
\end{lemma}
Here $B$ is used to identify $\g^*\cong \g$.
\begin{proof}
By definition,
\[ -4q(\tri)=\f{1}{6} \sum B([e^a,e^b],e^c) e_a e_b e_c,\]
using a basis $e_a$ of $\g$, with $B$-dual basis $e^a$. Take this
basis to be the Cartan-Weil basis, and use the Clifford relations to
write factors $e_{-\alpha}$ to the left and factors $e_\alpha$ to the
right. Then
\[ -4q(\tri)\in \Cl(\g)^T\subset \n_-\Cl(\g)\n_+\oplus \Cl(\t).\]
(For a $T$-equivariant element in $\Cl(\g)$, the $T$-weight of the $\n_-$-factors
must be compensated by the $T$ weights of the $\n_+$-factors.) Since
$ -4q(\tri)$ is an odd element of filtration degree $3$, and since
$\tri$ has no component in $\wedge^3\t$, it follows that
\[ -4q(\tri)\in \t\oplus \n_-\Cl(\g)\n_+.\]
To compute the $\t$-component, we calculate the constant component of
\[ [\xi,-4q(\tri)]_\Cl=\d^\Cl \xi=q(\lambda(\xi))\]
for any $\xi\in \t$. We have
\[ \lambda(\xi)=-\sum_{\alpha\succ 0} [\xi,e_{-\alpha}]\wedge e_\alpha
=2 \pi \sqrt{-1} \sum_{\alpha\succ 0} \l\alpha,\xi\r
e_{-\alpha} \wedge e_\alpha,\]
hence (see Sternberg \cite[Equation (9.25)]{st:lie})
\[ q(\lambda(\xi))=2 \pi \sqrt{-1}\ \sum_{\alpha\succ 0} \l\alpha,\xi\r
e_{-\alpha} e_\alpha +2\pi \sqrt{-1}\ \l\rho,\xi\r.\]
\end{proof}
As a consequence, we obtain,
\begin{proposition}\label{prop:cddiffa}
The element $x=\prod_{\alpha\succ 0}e_{\alpha}e_{-\alpha}\ \prod_i e_i$ satisfies,
\[ \Big(\d^\Cl-2\pi \sqrt{-1} \iota^\Cl(\rho)\Big) x=0.\]
\end{proposition}
\begin{proof}
$\d^\Cl$ is given as the Clifford commutator with $-4q(\tri)$. Since
$x$ is annihilated under both left and right multiplication by
elements of $\n_-\Cl(\g)\n_+$, it follows that
\[ \d^\Cl(x)=2\pi \sqrt{-1} [\rho,x]_\Cl.\]
\end{proof}

As a consequence, the Gauss-Dirac spinor satisfies the
differential equation:
\begin{equation}  \label{eq:gaussdiracspinor}
\big(\d+\eta-\tpi \varrho(\sfe(\rho))\big)\wh{\psi}_G=0.
\end{equation}
In fact, there is a more general version of this Equation, stated in
the following Proposition.  For any (real) dominant weight $\lambda$
of $G$ (not to be confused with the map $\lambda$ above), let
$\Delta_\lambda\in C^\infty(G)$ be the function
\[\Delta_\lambda(g)=\f{\l v_\lambda,g\cdot v_\lambda\r}{\l
  v_\lambda,v_\lambda\r},\]
where $v_\lambda$ is a highest weight vector in the irreducible
unitary representation $(V_\lambda,\l\cdot,\cdot\r)$ of highest
weight $\lambda$. The function $\Delta_\lambda$ is invariant under
the left-action of $N_-$, under the right-action of $N_+$, and
under the $T$-action it satisfies
\begin{equation}\label{eq:deltal}
\Delta_\lambda(tg)=\Delta_\lambda(gt)=t^\lambda \Delta_\lambda(g).
\end{equation}
Since $\Delta_\lambda(e)=1$, it follows that $\Delta_\lambda\not=0$ on
the big Gauss cell. We are interested in the product
$\Delta_\lambda\wh{\psi}_G$. Away from the zeroes of $\Delta_\lambda$,
this is a pure spinor defining the Gauss-Dirac structure. Similar to
$\wh{\psi}_G$, it is invariant under the left-action of $N_-$ and the
right-action of $N_+$, and satisfies
\begin{equation}\label{eq:equivprop}
l_t^*(\Delta_\lambda\wh{\psi}_G)=r_t^*(\Delta_\lambda\wh{\psi}_G)=
t^{\lambda+\rho}(\Delta_\lambda\wh{\psi}_G)\end{equation}
for all $t\in T$.
\begin{proposition}\label{prop:gddiff}
  For any dominant weight $\lambda$, the product
  $\Delta_\lambda\wh{\psi}_G$ satisfies the differential equation:
\begin{equation}\label{eq:deq}
 \big(\d+\eta-\tpi
 \varrho(\sfe(\lambda+\rho))\big)\Delta_\lambda\wh{\psi}_G=0,
\end{equation}
where $B$ is used to identify $\g^*\cong\g$.
\end{proposition}
\begin{proof}
Let $\s\subset \dd$ be the Lagrangian subalgebra \eqref{eq:Gaussh}
defining the Gauss-Dirac structure. We have, for all
$\zeta=(\xi,\xi')\in\s$,
\[
\begin{split}
\lefteqn{\varrho(\sfs(\zeta)) \big( \d+\eta-\tpi
\varrho(\sfe(\lambda+\rho))\big)\,\Delta_\lambda\wh{\psi}_G}\\
&=
\Big[\varrho(\sfs(\zeta)) ,\ \d+\eta-\tpi
\varrho(\sfe(\lambda+\rho))\Big]\Delta_\lambda\wh{\psi}_G \\
&=\big(\L(\xi^L-(\xi')^R))-\tpi B(\xi-\xi',\lambda+\rho)\big)
\Delta_\lambda\wh{\psi}_G=0,\end{split}
\]
where the last equality follows from the equivariance properties
\eqref{eq:equivprop} of $\Delta_\lambda\wh{\psi}_G$.  (Note that for
the elements of the form $\zeta=(\xi,0)$ with $\xi\in\n_+$ or
$\zeta=(0,\xi)$ with $\xi\in\n_-$, the inner product with
$\lambda+\rho\in \t$ vanishes.)  Hence, the left hand side of
\eqref{eq:deq} is annihilated by all $\sfs(\zeta)$, for $\zeta\in\s$.
Hence it is a function times $\wh{\psi}_G$, and thus vanishes since it
has parity opposite to that of $\wh{\psi}_G$.
\end{proof}

\begin{remark} \label{rem:gausscartan}
  The holomorphic Dirac structure $\wh{F}_G$ on $G=K^\C$ restricts to
  a complex Dirac structure $\wh{F}_G|_K=\wh{F}_K$ on the real Lie
  group $K$, with defining pure spinor the pull-back (restriction)
  $\wh{\psi}_K$ of $\wh{\psi}_G$. On the other hand,
  $E_G|_K=(E_K)^\C$.  In the notation of Section
  \ref{subsec:diraccoh}, applied to the Gauss-Cartan-splitting $(\T
  K)^\C=E_K^\C\oplus \wh{F}_K$, we have $\sigma=\tpi
  e(\rho)\in\Gamma((\T K)^\C)$, thus $\dirac_\pm=\d+\eta\pm \tpi
  \varrho(\sfe(\rho))$. As usual, $\dirac_+\phi_K=0,\
  \dirac_-\wh{\psi}_K=0$ (the second equation is the pull-back of
  \eqref{eq:gaussdiracspinor} to $K$). Let $\mu$ be the bi-invariant
  (real) volume form on $K$ defined by $\phi_K,\wh{\psi}_K$. Since
  $\dirac_\pm^2=\pm\tpi \, \L(\A_{\on{ad}}(\rho))$, the Dirac
  cohomology groups $H_\pm(E_K^\C, \wh{F}_K, \mu)$ are the cohomology
  groups of $\dirac_\pm$ on the space of
  $\A_{\on{ad}}(\rho)$-invariant complex-valued differential forms on
  $K$. These may be computed by the standard localization argument
  (\cite{be:ze}, see also \cite{hu:ext}): The set of zeroes of the
  vector field $\A_{\on{ad}}(\rho)$ on $K$ is just the maximal torus
  $T_K$, and the pull-back to $T_K$ intertwines $\dirac_\pm$ with
  $d\pm\tpi B(\theta_T, \rho))$, with $\theta_T$ the Maurer-Cartan
  form on $T_K$. Hence, by localization the pull-back to $T_K$ induces
  an isomorphism,
\[ H_\pm(E_K^\C, \wh{F}_K, \mu) \cong H(\Omega(T_K)^\C, \d\pm\tpi B(\theta_T, \rho))
\]
Since $\rho$ is a weight, it defines a $T_K$-character $t^\rho$, and
the operators $\d\pm\tpi B(\theta_T, \rho)$ are obtained from $\d$ by
conjugation by $t^{\pm\rho}$. Hence $H_\pm(E_K^\C, \wh{F}_K, \mu)
\cong H(T_K)^\C$.
\end{remark}

\vskip.3in
\section{q-Hamiltonian $G$-manifolds}\label{sec:qham}
In this section, we use the techniques developed in this paper to
extend the theory of group-valued moment maps, as developed in
\cite{al:mom,al:du} for the case of compact Lie groups, to more
general settings.

\subsection{Dirac morphisms and group-valued moment maps}
We briefly recall the definitions.
\begin{definition}\label{def:qham}
A \emph{quasi-Hamiltonian $\g$-manifold} (or simply
\emph{q-Hamiltonian $\g$-manifold}) is a manifold $M$ with a Lie
algebra action $\A_M\colon \g\to
\mf{X}(M)$, a 2-form
$\omega$, and a $\g$-equivariant \emph{moment map} $\Phi\colon
M\to G$ such that
\begin{equation}\label{eq:qHam}
\begin{split}
  \d\omega=\Phi^*\eta\ \ \ \ \ &\mbox{}\\
  \iota(\A_M(\xi))\omega=\Phi^*
  B(\xi,\textstyle{\f{\theta^L+\theta^R}{2}}) \ \ \ \ \
&\mbox{(moment map condition)}
  \\
  \ker(\omega_m)=\{\A_M(\xi)_m|\ \Ad_{\Phi(m)}\xi=-\xi\}\ \ \ \ \
&\mbox{(minimal degeneracy condition)}.
\end{split}
\end{equation}
If the action of $\g$ extends to an action of the Lie group $G$,
and if $\omega$ and $\Phi$ are equivariant for the action of $G$,
we speak of a \emph{q-Hamiltonian $G$-manifold}.
\end{definition}
The first two conditions in \eqref{eq:qHam} imply that $\omega$ is
$\g$-invariant (see \cite{al:mom}).  As shown by Bursztyn-Crainic
\cite{bur:di}, the definition of a q-Hamiltonian space may be
restated in Dirac geometric terms (see also  Xu \cite{xu:mom} for
another interpretation).
\begin{theorem}
There is a 1-1 correspondence between q-Hamiltonian
$\g$-manifolds,  
and manifolds $M$ together with a strong Dirac morphism
\begin{equation}\label{eq:qHamm}
 (\Phi,\omega)\colon (M,TM,0)\to (G,E_G,\eta).\end{equation}
More precisely, $(M,\A_M,\omega,\Phi)$ satisfies the first two
conditions if and only if $(\Phi,\omega)$ is a Dirac morphism, and
in this case the third condition is equivalent to this Dirac
morphism being \emph{strong}.
\end{theorem}
\begin{proof}
Let $(M,\A_M,\omega,\Phi)$ be a q-Hamiltonian $\g$-space. Given
$m\in M$, let $E'_{\Phi(m)}$ be
  the forward image of $T_mM$ under $((\d\Phi)_m,\omega_m)$:
\[ E'_{\Phi(m)}=\{(\d\Phi(v),\alpha)|\ v\in T_mM,\ (\d\Phi)_m^*\alpha=\iota(v)\omega_m\}.\]
Taking $v$ of the form $\A_M(\xi)_m$ for $\xi\in\g$, and using the
moment map condition, we see $E'_{\Phi(m)}\supset(E_G)_{\Phi(m)}$.
In fact, one has equality since both are Lagrangian subspaces.
This shows that $(\Phi,\omega)$ is a Dirac morphism. In
particular,
\[(\d\Phi)_m(\ker(\om_m))=\ker((E_G)_{\Phi(m)})=\{\A_{\sf{ad}}(\xi)_{\Phi(m)}| \Ad_{\Phi(m)}\xi=-\xi\}.\]
Hence, the minimal degeneracy condition holds if and only
$(\d\Phi)_m$ restricts to an isomorphism on $\ker(\om_m)$, i.e. if
and only if $(\Phi,\omega)$ is a strong Dirac morphism.
Conversely, given a strong Dirac morphism \eqref{eq:qHamm}, the
associated map $\mf{a}$ defines a $\g$-action
$\A_M(\xi)=\mf{a}(\Phi^*\sfe(\xi))$ on $M$, for which the map
$\Phi$ is $\g$-equivariant. The above argument then shows that
$(M,\A_M,\omega,\Phi)$ is a q-Hamiltonian $\g$-space.
\end{proof}

\begin{remark}
As a consequence of this result (or rather its proof), we see that
if $(M,\A_M,\omega,\Phi)$ satisfies the first two conditions in
\eqref{eq:qHam}, then the third condition (minimal degeneracy) is
equivalent to the transversality property \cite{bur:int,xu:mom}
\[ \ker(\omega)\cap \ker(\d\Phi)=\{0\}.\]
\end{remark}
\begin{remark}
There is a similar result for q-Hamiltonian $G$-manifolds. Here, it is
necessary to assume the existence of a $G$-action on $M$ for which
the Dirac morphism $(\Phi,\omega)$ is equivariant, and such that the infinitesimal action coincides with that defined by $\mf{a}$.
\end{remark}

\begin{example}
By Example \ref{ex:leaves}, the inclusion of the conjugacy classes
$\Co$ in $G$, with 2-forms defined by the Cartan-Dirac structure,
defines a strong Dirac morphism $(\iota_\Co,\om_\Co)$. Thus,
conjugacy classes are q-Hamiltonian $G$-manifolds.
\end{example}

Using our results on the Cartan-Dirac structure, it is now
straightforward to deduce the basic properties of q-Hamiltonian
spaces $(M,\A_M,\omega,\Phi)$. In contrast with the original
treatment in \cite{al:mom}, the discussion works equally well for
non-compact Lie groups, and also in the holomorphic category.

\begin{theorem}[Fusion]\label{th:fusion}
Let $(M,\A_M,\Phi,\om)$ be a q-Hamiltonian $G\times G$-manifold.
Let $\A_{\on{fus}}$ be the diagonal $G$-action,
$\Phi_{\on{fus}}=\on{Mult}\circ \Phi$, and
$\om_{\on{fus}}=\om+\Phi^*\varsigma$, with
$\varsigma\in\Om^2(G^2)$ the 2-form defined in
\eqref{eq:varsigma}. Then
$(M,\A_{\on{fus}},\Phi_{\on{fus}},\om_{\on{fus}})$ is a
q-Hamiltonian $G$-manifold. (An analogous statement holds for
q-Hamiltonian $\g\times\g$-manifolds.)
\end{theorem}
\begin{proof}
Since
\[(\Phi_{\on{fus}},\om_{\on{fus}})=(\on{Mult},\varsigma)
\circ (\Phi,\omega)\]
is a composition of two strong Dirac morphism, it is itself a
strong Dirac morphism from $(M,TM,0)$ to $(G,E_G,\eta)$. The
induced map $M\times\g=\Phi_{\on{fus}}^*E_G\to TM$ is a
composition of the map $\on{Mult}^*E_G\to E_{G\times G}$ defined
by the strong Dirac morphism $(\on{Mult},\varsigma)$, with the map
$\Phi^* E_G\times E_G\to TM$ given by the strong Dirac morphism
$(\Phi,\omega)$. If we use the sections $\sfe(\xi)$ to identify
$E_G \cong G\times \g$, the latter map is the $\g\times\g$-action
on $M$, while the former is the diagonal inclusion $\g\to
\g\times\g$. This confirms that the resulting action is just the
diagonal action.
\end{proof}

If $M=M_1\times M_2$ is a direct product of two q-Hamiltonian
manifolds, the quadruple
$(M,\A_{\on{fus}},\Phi_{\on{fus}},\om_{\on{fus}})$ is called the
\emph{fusion product} of $M_1,M_2$. In particular we obtain products
of conjugacy classes as new examples of q-Hamiltonian $G$-spaces.

Suppose $(M,\A_M,\omega_0,\Phi_0)$ is a Hamiltonian $\g$-manifold: That
is, $\omega_0$ is symplectic, and $\Phi_0\colon M\to \g^*$ is the
moment map for a Hamiltonian $\g$-action on $M$. As is well-known,
this is equivalent to $\Phi_0$ being a Poisson map from the
symplectic manifold $(M,\omega_0)$ to the Poisson manifold
$(\g^*,\pi_{\g^*})$. But this is also equivalent to
\[(\Phi_0,\omega_0)\colon (M,TM,0)\to (\g^*,\on{Gr}_\pi,0)\]
being a strong Dirac morphism. A \emph{Hamiltonian $G$-manifold}
comes with a $G$-action on $M$ integrating the $\g$-action, and
such that the Dirac morphism $(\Phi_0,\omega_0)$ is equivariant.
Given an invariant inner product $B$ on $\g$, used to identify
$\g^*\cong\g$, we may compose the Dirac morphism
$(\Phi_0,\omega_0)$ with the Dirac morphism $(\exp,\varpi)$ from
Theorem \ref{th:exp}, and obtain:
\begin{theorem}[Exponentials] \label{th:exp2}
Suppose $(M,\A_M,\omega_0,\Phi_0)$ is a Hamiltonian $G$-manifold,
and let $\om=\om_0+\Phi_0^*\varpi$, $\Phi=\exp\circ \Phi_0$. Then
$(M,\A_M,\omega,\Phi)$ satisfies the first two conditions in
\eqref{eq:qHam}. On $M_\natural=\Phi_0^{-1}(\g_\natural)$, the
third condition (minimal degeneracy) holds as well, thus
$(M_\natural,\A_M,\omega,\Phi)$ is a q-Hamiltonian $G$-manifold.
(Similar statements hold for q-Hamiltonian $\g$-manifolds.)
\end{theorem}

\subsection{Volume forms}
Any symplectic manifold $(M,\om)$ carries a distinguished volume form,
given as the top degree component $\exp(\om)^{[\dim
  M]}=\f{1}{n!}\om^n$. For a q-Hamiltonian $G$-manifold
$(M,\A_M,\om,\Phi)$, the 2-form $\om$ is usually degenerate, hence
$\exp(\om)^{[\mathrm{top}]}$ will have zeroes. Nevertheless, any
q-Hamiltonian $G$-manifold carries a distinguished volume form,
provided the adjoint action $\on{Ad}\colon G\to \on{O}(\g)$ lifts to
$\on{Pin}(\g)$:
\begin{theorem}[Volume forms] Suppose the adjoint action $\on{Ad}\colon
G\to \on{O}(\g)$ lifts to $\on{Pin}(\g)$, and let
$\psi_G\in\Om(G)$ be the pure spinor defined by this lift. For any
q-Hamiltonian $G$-manifold $(M,\A_M,\om,\Phi)$, the differential
form
\begin{equation}\label{eq:volume}
\mu_M=(\exp(\omega)\wedge \Phi^*\psi_G)^{[\dim M]}
\end{equation}
is a volume form. It has the equivariance property
$\A_M(g)^*\mu_M=\det(\Ad_g)\,\mu_M$. More generally, if
$(M,\A_M,\om,\Phi)$ satisfies the first two conditions in
\eqref{eq:qHam}, the form $\mu_M$ is non-zero exactly at those
points where $\omega$ satisfies the minimal degeneracy condition.
\end{theorem}
Of course, the factor $\det(\Ad_g)=\pm 1$ is trivial if $G$ is
connected.

\begin{proof}
Since $\psi_G$ is a pure spinor defining the complementary
Lagrangian subbundle $F_G$, and since $(\Phi,\omega)$ is a strong
 Dirac morphism, the pull-back $\Phi^*\psi_G$ is non-zero
everywhere. Furthermore, $\exp(\omega)\Phi^*\psi_G$ is a pure
spinor defining the backward image $F$ of $F_G$ under the Dirac
morphism $(\Phi,\omega)$. Since $F$ is transverse to $TM$ (see
Proposition \ref{prop:transI}), the top degree part of
$\exp(\omega)\Phi^*\psi_G$ is nonvanishing. More generally, if
$(M,\A_M,\om,\Phi)$ only satisfies the first two conditions in
\eqref{eq:qHam}, then the above argument applies at all points of
$M$ where $(\Phi,\omega)$ is a strong Dirac morphism. But these
are exactly the points where $\Phi^*\psi_G$ is non-zero.

The equivariance property of $\mu_M$ is a direct consequence of the
equivariance properties of $\phi_G$ and $\psi_G$ described in
Proposition \ref{prop:phipsi}.
\end{proof}

The volume form $\mu_M$ is called the \emph{Liouville volume form}
of the q-Hamiltonian $G$-manifold $(M,\A_M,\om,\Phi)$. Let
$|\mu_M|$ be the associated measure.  If the moment $\Phi$ is
proper, the push-forward $\Phi_*|\mu_M|$ is a well-defined measure
on $G$, called the \emph{Duistermaat-Heckman measure}.

\begin{remark}
  For the case of compact Lie groups, the q-Hamiltonian Liouville
  forms and Duistermaat-Heckman measures were introduced in
  \cite{al:du}. The fact that $\mu_M$ is a volume form was verified by
  `direct computation'. However, the argument in \cite{al:du} does not
  extend to non-compact Lie groups.
\end{remark}

\begin{remark}
The expression $\exp(\om)\Phi^*\psi_G$ entering the definition of
the volume form $\mu_M$ satisfies the differential equation
\begin{equation}\label{eq:psieqn3}
(\d+ \iota(\A_M(\tri)))\Big(\exp(\om)\Phi^*\psi_G\Big)=0.\end{equation}
This follows from the differential equation \eqref{eq:psieqn} for
$\psi_G$ together with Remark \ref{rem:reference}\eqref{it:aaa}.
\end{remark}
\begin{proposition}
  Suppose $(M,\A_M,\omega,\Phi)$ is a q-Hamiltonian $G$-manifold, and
  that $\on{\Ad}$ lifts to the Pin group. Then $M$ is even-dimensional
  if $\det(\Ad_\Phi)=+1$, and odd-dimensional if $\det(\Ad_\Phi)=-1$.
In particular, it is even-dimensional when $G$ is
connected, and in this case $M$ carries a canonical orientation.
\end{proposition}
\begin{proof}
  The construction of $\psi_G$ in terms of the map $\ca{R}$ (see
  Proposition \ref{prop:phipsi}) shows that the form $\psi_G$ has even
  degree at points $g\in G$ with $\det(\Ad_g)=1$, and odd degree at
  points with $\det(\Ad_g)=-1$.
Hence, the parity of the volume form $\mu_M$ is determined by the
parity of $\det(\Ad_\Phi)$. If $G$ is connected, the lift of $\on{Ad}$
(which exists by assumption) is unique, and $\det(\Ad_g)\equiv 1$.
\end{proof}
Without the existence of a lift to $\on{Pin}(\g)$, the form
$\psi_G$ is only defined locally, up to sign. That is, we still
obtain a $G$-invariant \emph{measure} on $M$, given locally as $(e^{\omega}
\Phi^*\psi_G)^{[\mathrm{top}]}$. It is interesting to specialize these
results to conjugacy classes:

\begin{theorem}
  Suppose $G$ is a connected Lie group, whose Lie algebra carries an
  invariant inner product $B$. Then:
\begin{enumerate}
\item Every conjugacy class $\Co\subset
G$ carries a distinguished invariant measure
(depending only on $B$).
\item The conjugacy class $\Co$ of $g\in G$ is even-dimensional if and only if
      $\det(\Ad_g)=+1$.
\item
If the adjoint action $G\to
\on{O}(\g)$ lifts to $\on{Pin}(\g)$, then every conjugacy class carries a
distinguished orientation.
\end{enumerate}
\end{theorem}

\begin{example}
  Consider the conjugacy classes of $G=\on{O}(2)$: If $g\in \SO(2)$,
  the conjugacy class of $g$ is zero-dimensional, consisting of either
  one or two points. On the other hand, the circle
  $\on{O}(2)\backslash \SO(2)\cong S^1$ forms a single conjugacy class.
  Similarly, for $G=\on{O}(3)$, the elements $g\in G$ with
  $\det(g)=-1$ have $\det(\Ad_g)=1$. Each of these form a single
  2-dimensional conjugacy class. The group $\on{SO}(3)$ is the simplest
  example where the adjoint action $G\to \SO(\g)$ (which in this case
  is just the identity map) does not lift to the spin group.
  Indeed the conjugacy class of rotations by $180^o$ is isomorphic to
  $\R P(2)$, hence non-orientable.
\end{example}

\begin{example}\label{ex:homspace}
  Suppose $G$ carries an involution $\sig$, such that the
  corresponding involution of $\g$ preserves $B$. Form the semi-direct
  product $G\rtimes \Z_2$, where the action of $\Z_2$ is generated by
  the involution $\sig$. The $G\rtimes\Z_2$-conjugacy class of the
  element $(e,\sigma)$ is isomorphic to the homogeneous space
  $M=G/G^\sig$, which therefore is an example of a q-Hamiltonian
  $G\rtimes \Z_2$-space. The 2-form on $M$ is just zero. Let us
  compute the Liouville measure on $M$, for the case that the
  restriction of $B$ to $\g^\sig=\ker(\sig-1)$ is still
  non-degenerate. Let $e_1,\ldots,e_n$ be a basis of $\g$, with
  $B(e_i,e_j)=\pm \delta_{ij}$, such that $e_1\ldots,e_k$ are a basis
  of $\g^{\sig}$. Then
\[ \wt{\sig}=2^{(n-k)/2} e_{k+1}\cdots e_n\in \on{Pin}(\g)\]
is a lift of $\sig$. Note that $\wt{\sig}^2=\pm 1$, with sign
depending on $n-k$. Taking $\mu=e_1\wedge \cdots \wedge e_n$ as
the Riemannian volume form on $\g$, we have
\[ \wt{\sig}q(\mu)=\pm 2^{(n-k)/2} e_1\ldots e_k\]
so $\star q^{-1}(\wt{\sig}q(\mu))=\pm 2^{(n-k)/2} e_{k+1}\wedge \cdots
\wedge e_n$.
We conclude that the Liouville measure on $M=G/G^\sig$ coincides with
the $G$-invariant measure defined by the metric on
$(\g^\sig)^\perp\subset \g$.
\end{example}

\begin{proposition}[Volume form for `fusions']\label{prop:fusvol}
The volume form of a q-Hamiltonian $G\times G$-manifold
$(M,\A_M,\om,\Phi)$ (as in Theorem \ref{th:fusion}) coincides with
the volume form of its fusion
$(M,\A_{\on{fus}},\om_{\on{fus}},\Phi_{\on{fus}})$:
\[
(\exp(\om)\, \Phi^*\psi_{G\times G})^{[\dim M]}
=(\exp(\om_{\on{fus}})\, \Phi_{\on{fus}}^* \psi_G)^{[\dim M]}.
\]
\end{proposition}
\begin{proof}
Using Theorem \ref{th:psipull}, we have
\[ \begin{split}
\exp(\om_{\on{fus}})\, \Phi_{\on{fus}}^* \psi_G&=
\exp(\om+\Phi^*\varsigma)\, \Phi^*\on{Mult}^*\psi_G\\
&=\exp(\om)\, \Phi^*\big(\varrho(\exp(-\sfe(\gamma)))\psi^1_G\otimes \psi^2_G\big)\\
&=\exp(-\iota(\A_M(\gamma))) (\exp(\om)\Phi^*\psi_{G\times G}),
\end{split}
\]
where we used Remark \ref{rem:reference}\eqref{it:aaa} for the last equality. Since the
operator $\exp(-\iota(\A_M(\gamma)))$ does not affect the top
degree part, the proof is complete.
\end{proof}

\begin{example}
  An important example of a q-Hamiltonian $G$-space is the double
  $D(G)=G\times G$, with moment map the commutator
  $\Phi(a,b)=aba^{-1}b^{-1}$. As explained in \cite{al:du} the double
  is obtained by fusion, as follows: Start by viewing the Lie
  group $G$ as a homogeneous space $G=G\times G/G_\Delta$, where $G_\Delta$ is
  the
  diagonal subgroup. Since $G_\Delta$ is the fixed point set for the
  involution $\sig$ of $G\times G$ switching the two factors, we see
  as in Example~\ref{ex:homspace} that $G$ is a q-Hamiltonian $(G\times
  G)\rtimes \Z_2$-space, with moment map $a\mapsto (a,a^{-1},\sig)$.
  The Liouville measure is simply the Haar measure on $G$. Fusing two
  copies, the direct product $G\times G$ becomes a q-Hamiltonian
  $G\times G$-space. Finally, passing to the diagonal action one
  arrives at the double $D(G)$. By Proposition \ref{prop:fusvol}, the
  resulting Liouville measure on $D(G)$ is just the Haar measure.
\end{example}

\begin{proposition}[Volume form for `exponentials']\label{prop:volexp}
Let $(M,\A_M,\Phi_0,\om_0)$ be a Hamiltonian $G$-space, and
$(M,\A_M,\Phi,\om)$ its `exponential', as in Theorem
\ref{th:exp2}. Then
\[ (\exp(\om)\Phi^*\psi_G) ^{[\dim M]}=\Phi_0^* J^{1/2}\
\exp(\om_0)^{[\dim M]}.\]
\end{proposition}
\begin{proof}
Using the relation \eqref{eq:psiexp} between $\exp^*\psi_G$ and
$\psi_\g=1$, we find
\[ \begin{split}
\exp(\om)\, \Phi^* \psi_G&=
\exp(\om_0+\Phi_0^*\varpi) \Phi_0^* \exp^* \psi_G\\
&=\exp(\om_0) \Phi_0^* J^{1/2} \varrho(\wt{A}^{-\sfe_0(\varepsilon)})(1)\\
&= \Phi_0^*J^{1/2}\exp(-\iota(\A_M(\varepsilon))) \exp(\om_0).
\end{split}\]
Since $\exp(-\iota(\A_M(\varepsilon)))$ does not affect the top
degree part, the proof is complete.
\end{proof}


\subsection{The volume form in terms of the Gauss-Dirac
spinor}\label{subsec:gaussvol} Suppose now that $K$ is a
\emph{compact} Lie group, with complexification $G=K^\C$, and let
$B\colon\g\times\g\to \C$ be the complexification of a positive
definite inner product on $\k$. In this case, as discussed in
Section~\ref{subsec:gauss}, $E_G$ has a second Lagrangian
complement $\wh{F}_G$, defined by the Gauss-Dirac spinor
$\wh{\psi}_G$.
Its
pull-back to $K \subset G$, denoted by $\wh{\psi}_K$, is thus a
complex-valued pure spinor defining a (complex) Lagrangian
complement $\wh{F}_K\subset (\T K)^\C$.

Given a q-Hamiltonian $K$-space $(M,\A_M,\Phi,\om)$, the complex
differential form $\exp(\omega)\Phi^* \wh{\psi}_K$ is related to
$\exp(\omega)\Phi^*\psi_K$ by the $r$-matrix,
\[ \exp(\omega)\Phi^* \wh{\psi}_K=
\exp(-\iota(\A_M(\mf{r})))\Big(\exp(\om)\Phi^*\psi_K\Big).\]
Since $\exp(-\iota(\A_M(\mf{r})))$
does not affect
the top degree part, it follows that we can write our volume form
also in terms of $\wh{\psi}_K$:
\[ \mu_M=\Big(\exp(\omega)\Phi^* \wh{\psi}_K\Big)^{[\dim M]}.\]
\begin{remark}
%
%
  Let $\ti{F}_M$ be the backward image of $\wh{F}_K$ under the strong
  Dirac morphism $(\Phi,\om)\colon (M, TM, 0) \to (K, E_K, \eta)$.
  Since $\ti{F}_M$ is transverse to $TM^\C$, it is given by a graph of
  a (complex-valued) bivector $\pi$, and $H_-(TM^\C,\ti{F}_M,
  \mu_M)\cong H_\pi(M)= H(\Omega(M)^{X_{\pi}}, d-\iota(X_{\pi}))$. A
  simple calculation shows that $X_{\pi}=\tpi\A_M(\rho)$ (where $B$ is
  used to identify $\k^*\cong\k$).

  The pure spinors $\phi_M=1$ and $\phi_K$
  satisfy $d\phi_M=0$ and $(d+\eta)\phi_K=0$. Hence, by Proposition
  \ref{prop:psipsi'} the map $e^\omega \Phi^*$ descends to Dirac
  cohomology, $H_-(E_K^\C, \wh{F}_K, \mu_K) \to H_\pi(M)$. In particular, $\dirac_-\wh{\psi}_K=0$ implies
  that $\exp(\omega)\Phi^* \wh{\psi}_K$ is closed under
  $\d-\tpi\iota(\A_M(\rho))$.
 For $M$ is compact, the class $[e^\omega \Phi^*\wh{\psi}_K]$ in
  $H_\pi(M)$ is nonvanishing because
  its integral is $\int_M{\mu_M}>0$.
\end{remark}

Let $\Delta_\lambda\colon G\to \C$ be the holomorphic functions
introduced in Section \ref{subsec:gdspinor}.

\begin{proposition}\label{prop:gaussvol}
For any dominant weight $\lambda$, the complex differential form
$\exp(\omega)\Phi^* (\Delta_\lambda \wh{\psi}_K)$ satisfies the
differential equation
\begin{equation}\label{eq:diff-eqn} (\d
  -\tpi \iota(\A_M(\lambda+\rho))\big)\Big(\exp(\omega)\Phi^* (\Delta_\lambda
  \wh{\psi}_K)\Big)=0.
\end{equation}
Her $B_K$ is used to identify $\k^*\cong\k$.
\end{proposition}
\begin{proof}
This follows from the differential equation for the Gauss-Dirac spinor,
Proposition \ref{prop:gddiff}, together with Remark \ref{rem:reference}\eqref{it:aaa}.
\end{proof}
As remarked
in \cite{al:gr}, the orthogonal projection of $\dim V_\lambda
\Delta_\lambda|_K$ to the $K$-invariant functions on $K$ coincides
with the irreducible character $\chi_\lambda$ of highest weight
$\lambda$. Thus,
\[ \begin{split}
  \int_M \exp(\omega)\Phi^* (\Delta_\lambda \wh{\psi}_K)&= \int_M
  |\mu_M| \Phi^*\Delta_\lambda \\&=\int_K
  \Phi_*|\mu_M|\,\Delta_\lambda
  =(\dim V_\lambda)^{-1} \int_K
  \chi_\lambda\, \Phi_*|\mu_M|.\end{split}\]
On the other hand, by \eqref{eq:diff-eqn} the integral may be computed
by localization \cite{be:ze} to the zeroes of the vector field
$\A_M(\lambda+\rho)$. As shown in \cite{al:gr}, the 2-form $\om$ pulls
back to symplectic forms $\om_Z=\iota_Z^*\om$ on the components $Z$ of
the zero set, and the restriction $\Phi_Z=\iota_Z^*\Phi$ takes values
in $T$. Since
$\iota_T^*(\Delta_\lambda\wh{\psi}_K)(t)=t^{\lambda+\rho}$ for $t\in
T$, one obtains the following formula for the Fourier coefficients of
the q-Hamiltonian Duistermaat-Heckman measure:
\[ \int_K \chi_\lambda \Phi_*|\mu_M|=\dim V_\lambda \sum_{Z\subset
\A_M(\lambda+\rho)^{-1}(0)}\int_Z
\f{\exp(\omega_Z)(\Phi_Z)^{\lambda+\rho}}{\Eul(\nu_Z,\tpi
(\lambda+\rho))}.\]
Here $\Eul(\nu_Z,\cdot)$ is the $T$-equivariant Euler form of the
normal bundle. This formula was proved in \cite{al:gr}, using a
more elaborate argument. Taking $\lambda=0$, one obtains a formula
for the volume $\int_M |\mu_M|$ of $M$.


\subsection{q-Hamiltonian q-Poisson $\g$-manifolds}
Just as any symplectic 2-form
determines a Poisson bivector $\pi$, any q-Hamiltonian $G$-manifold
carries a distinguished bivector field $\pi$. However, since $\om$ is
not non-degenerate $\pi$ is not simply obtained as an inverse, and
also $\pi$ is not generally a Poisson structure.

Suppose $(M,\A_M,\om,\Phi)$ is a q-Hamiltonian $\g$-manifold, or
equivalently that $(\Phi,\omega)$ is a strong Dirac morphism
$(M,TM,0)\to (G,E_G,\eta)$. Let $\wt{F}\subset \TM$ be the
backward image of $F_G$ under this Dirac morphism. It is a
complement to $TM$, hence it is of the form $\wt{F}=\on{Gr}_\pi$
for some $\g$-invariant bivector field $\pi\in\mf{X}^2(M)$. By
Proposition \ref{prop:diracmaps}(c), the Schouten bracket of this
bivector field with itself satisfies
\begin{equation}\label{eq:qHP1}
\hh [\pi,\pi]_{\on{Sch}}=\A_M(\tri).
\end{equation}
Let $\sfp'\colon \TG\to E_G$ be the projection along $F_G$. Let
$\{v_a\}$ and $\{v^a\}$ be bases of $\g$ with
$B(v_a,v^b)=\delta_a^b$. Then
$\sfp'(x')=\sum_a \l x,\,\sff(v_a)\r \sfe(v^a)$ for all $x'\in \Gamma(\TG)$. For
$\alpha'\in \Om^1(G)\subset \Gamma(\TG)$, we have $\l \alpha',\,\sff(v^a)\r=\hh\l
\alpha',v_a^L+v_a^R\r\ \sfe(v^a)$. Hence, \eqref{eq:relation} shows
that
\begin{equation}\label{eq:qHP2}
 \pi^\sharp \Phi^*\alpha'=-\sum_a \Phi^* \l
\alpha',{\ts \f{v_a^L+v_a^R}{2}}\r\ \A_M(v^a),\ \
\alpha'\in\Om^1(G),
\end{equation}
and, by \eqref{eq:Erange}, we have:
\begin{equation}\label{eq:qHP3}
 \on{ran}(\A_M)+\on{ran}(\pi^\sharp)=TM.
\end{equation}
This last condition can be viewed as a counterpart to the
invertibility of a Poisson bivector defined by a symplectic form.
Dropping this condition, one arrives at the following definition:
\begin{definition}\cite{al:ma,al:qu}
  A \emph{q-Hamiltonian q-Poisson $\g$-manifold} is a manifold $M$,
  together with a Lie algebra action $\A_M\colon \g\to \mf{X}(M)$, a
  $\g$-invariant bivector field $\pi$, and a $\g$-equivariant
  \emph{moment map} $\Phi\colon M\to G$, such that conditions
  \eqref{eq:qHP1} and \eqref{eq:qHP2} are satisfied. If the
  $\g$-action on $M$ integrates to a $G$-action, such that $\pi,\Phi$
  are $G$-equivariant, we speak of a \emph{q-Hamiltonian q-Poisson
  $G$-manifold}.
\end{definition}

\begin{example}
  The basic example of a Hamiltonian Poisson $G$-manifold is provided
  by the coadjoint action on $M=\g^*$, with $\pi=\pi_{\g^*}$ the Kirillov
  bivector and moment map the identity map. Similarly, the quadruple
  $(G,\A_{\sf{ad}},\pi_G,\on{id})$, with $\pi_G$ the bivector field
  \eqref{eq:pig}, is a q-Hamiltonian q-Poisson $G$-manifold.
\end{example}
The techniques in this paper allow us to give a much simpler proof to
the following theorem from \cite{bur:di}:
\begin{theorem}\label{the:equiv}
  There is a 1-1 correspondence between q-Hamiltonian q-Poisson
  $\g$-manifolds $(M,\A_M,\pi,\Phi)$, and Dirac manifolds $(M,E_M,\eta_M)$
  equipped with a strong Dirac morphism
\begin{equation}
(\Phi,0)\colon (M,E_M,\eta_M)\to (G,E_G,\eta).
\end{equation}
Under this correspondence,
$\on{ran}(E_M)=\on{ran}(\A_M)+\on{ran}(\pi^\sharp)$.
\end{theorem}
\begin{proof}
  Suppose $(\Phi,0)\colon (M,E_M,\eta_M)\to (G,E_G,\eta)$ is a strong
  Dirac morphism.  Consider the bundle map $\mf{a}\colon \Phi^*E_G\to TM$
  defined by $\Phi$ (see Section \ref{subsec:morph}).  By
  Proposition \ref{prop:diracmaps}(c), the vector fields
  $\A_M(\xi)=\mf{a}(\sfe(\xi))\in \mf{X}(M)$ define a Lie algebra action of
  $\g$ on $M$ for which $\Phi$ is equivariant. Note also that since
  $\on{ran}(\mf{a})\subset \on{ran}(E_M)$, this action preserves the
  leaves $Q\subset M$ of $E_M$. In fact, the bundle $E_M$ is
  $\g$-invariant: If $E_M=\on{Gr}_\omega$ this follows from the $\g$-invariance of
  $\omega$ (see comment after Def.~\ref{def:qham}), and
  in the general case it follows since $E_M|_Q$ is invariant, for
  any leaf
  $Q$.  Let $F_M$ be the backward image of $F_G$ under $(\Phi,0)$, and
  $\pi\in\mf{X}^2(M)$ be the bivector field defined by the splitting
  $\TM=E_M\oplus F_M$. Then $\pi$ is $\g$-invariant (since $E_M,F_M$ are). Equation
  \eqref{eq:qHP1} follows from Proposition \ref{prop:diracmaps}(d), while
  Equation \eqref{eq:qHP2} is a consequence of Theorem
  \ref{prop:reconstruct}, Equation \eqref{eq:relation}.

  Conversely, given a quasi-Poisson $\g$-manifold $(M,\A_M,\pi,\Phi)$,
  let $\mf{a}\colon \Phi^* E_G\to TM$ be the bundle map given on
  sections by $\Phi^* \sfe(\xi)\mapsto \A_M(\xi)$.  The $\g$-equivariance
  of $\Phi$ implies that $\Phi\circ
  \mf{a}=\pr_{\Phi^*TG}|_{\Phi^*E_G}$.  Theorem~\ref{prop:reconstruct}
  provides a Lagrangian splitting $\TM=E_M\oplus F_M$ such that $F_M$ is the
  backward image of $F_G$ and $E_G$ is the forward image of $E_M$. It
  remains to check the integrability condition of $E_M$ relative to the
  3-form $\eta_M=\Phi^*\eta$. Let $\Upsilon^E\in\Gamma(\wedge^3 F_M)$ be
  the Courant tensor of $E_M$. We have to show that $\Upsilon^E=0$, or
  equivalently that $\Gamma(E_M)$ is closed under the $\eta_M$-twisted Courant bracket.
  Recall that $E_M$ is spanned by the sections of two types:
\[
\wh{\A}_M(\xi):=\wh{\mf{a}}(\Phi^* \sfe(\xi))=\A_M(\xi)\oplus \Phi^*
B({\ts \f{\theta^L+\theta^R}{2}},\xi)\] for $\xi\in \g$, and
 sections $h(\alpha)$, for $\alpha\in \Om^1(M)$, where the map $h$ is
 defined as in \eqref{eq:hdef}, with $\V$ replaced with $\TM$, and with
 $\om=0$.  Since $\wh{\mf{a}}$ is a comorphism of Lie algebroids
 (cf. Proposition \ref{prop:comor}), we have
\begin{equation}\label{eq:eqn1}
 \Cour{\wh{\A}_M(\xi_1), \wh{\A}_M(\xi_2)}_{\eta_M}=
\wh{\A}_M([\xi_1,\xi_2]).\end{equation}
Furthermore, since $\pi$ is $\g$-invariant, it follows from
\eqref{eq:hdef} that the map $h$ is $\g$-equivariant, and
therefore
\[ [\varrho(h(\alpha)),[\varrho(\wh{\A}_M(\xi)),\d+\Phi^*\eta]]
=[\varrho(h(\alpha)),\L(\A_M(\xi))]=
-\varrho(h(\L(\A_M(\xi)\alpha)).\]
Thus
\begin{equation}\label{eq:eqn2}
 \Cour{\wh{\A}_M(\xi), h(\alpha)}_{\eta_M}=h(\L(\A_M(\xi)) \alpha)
\end{equation}
by definition of the Courant bracket. Equations \eqref{eq:eqn1}
and \eqref{eq:eqn2} show that $ \Cour{\wh{\A}_M(\xi),
\cdot}_{\eta_M}$ preserves $\Gamma(E_M)$. Thus
$\Upsilon^E(x_1,x_2,x_3)$ vanishes if one of the three sections
$x_i\in\Gamma(E_M)$ lies in the range of $\wh{\A}_M$.  It remains to
show that $\Upsilon^E(h(\alpha_1),h(\alpha_2),h(\alpha_3))=0$ for
all 1-forms $\alpha_i$, or equivalently that $h^* \Upsilon^E=0$,
where $h^*\colon F_M\to TM$ is the dual map to $h\colon T^*M\to
E_M=F_M^*$.  Since $h=\sfp|_{T^*M}$, where $\sfp:\TM\to E_M$ is the
projection along $F_M$ (see \eqref{eq:hdef}), we have
$h^*=\pr_{TM}|_{F_M}$.  Thus, we must show that
$\pr_{TM}\Upsilon^E=0$. By Proposition \ref{prop:diracmaps}(b),
and the defining property of q-Hamiltonian q-Poisson spaces, we
have
\[\pr_{TM}(\Upsilon^F) =\mf{a}(\Phi^*\Upsilon^{F_G})=\A_M(\tri)=
\hh[\pi,\pi]_{\on{Sch}}.
\]
On the other hand, Theorem \ref{th:bivector}(a) gives
$\pr_{TM}(\Upsilon^E)+\pr_{TM}(\Upsilon^F) -
\hh[\pi,\pi]_{\on{Sch}} =0$. Taking the two results together, we
obtain $\pr_{TM}(\Upsilon^E)=0$ as desired.
\end{proof}

As an immediate consequence, the data $(M,\A_M,\pi,\Phi)$ defining
a q-Hamiltonian q-Poisson $G$-manifold are equivalent to the data
of a $G$-equivariant Dirac manifold $(M,E_M,\eta_M)$, equipped with
a $G$-equivariant Dirac morphism $(\Phi,0)$, for which the
$G$-action on $M$ integrates the $\g$-action defined by the Dirac
morphism.

\begin{proposition}[Fusion] Suppose $(M,\A_M,\pi,\Phi)$ is a
  q-Hamiltonian q-Poisson $\g\times\g$-manifold. Let $\A_{\on{fus}}$
  be the diagonal $\g$-action, $\Phi_{\on{fus}}=\on{Mult}\circ\Phi$,
  and $\pi_{\on{fus}}=\pi+\A_M(\gamma)$. Then
  $(M,\A_{\on{fus}},\pi_{\on{fus}},\Phi_{\on{fus}})$ is a
  q-Hamiltonian q-Poisson $\g$-manifold.
\end{proposition}
\begin{proof}
By Theorem \ref{the:equiv}, the given q-Poisson
$\g\times\g$-manifold corresponds to a Dirac manifold
$(M,E_M,\eta_M)$ such that $(\Phi,0)$ is a Dirac morphism into
$(G,E_G,\eta)\times (G,E_G,\eta)$. Thus,
$\eta_M=\Phi^*(\eta_G^1+\eta_G^2)$.  The bivector field $\pi$ is
defined by the Lagrangian splitting $\TM=E_M\oplus F_M$, where $F_M$ is
the backward image of $F_G^1\oplus F_G^2$ under $(\Phi,0)$.
Composing with $(\on{Mult},\varsigma)$ (cf.
Thm.~\ref{prop:multcardir}), we obtain a strong Dirac morphism,
\[
(\Phi_{\on{fus}},\Phi^*\varsigma)\colon (M,E_M,\eta_M)\to
(G,E_G,\eta),
\]
which in turn defines a q-Hamiltonian q-Poisson $\g$-manifold.
%
%
Let $\wt{F}_M$ be the backward image of $F_G$ under this Dirac
morphism. By Proposition \ref{prop:egamma}, $\wt{F}$ is related to
$F$ by the section $\wh{\A}_M(\gamma)\in \Gamma(\wedge^2 E_M)$, where
$\wh{\A}_M\colon \g\times\g \to E_M$ is the map defined by the Dirac
morphism $(\Phi,0)$.  Hence, by Proposition \ref{prop:gaugee}, the
bivector for the new splitting $\TM=E_M\oplus \wt{F}_M$ is
$\pi_{\on{fus}}=\pi+\A_M(\gamma)$.
\end{proof}

\begin{proposition}[Exponentials] Suppose $(M,\A_M,\pi_0,\Phi_0)$ is a
  Hamiltonian Poisson $\g$-manifold. That is, $\A_M$ is a $\g$-action on
  $M$, $\pi_0$ is a $\g$-invariant Poisson structure, and
  $\Phi_0\colon M\to \g$ is a $\g$-equivariant moment map generating
  the given action on $M$. Assume that $\Phi_0(M)\subset \g_\natural$,
  and let
\[\Phi=\exp\circ \Phi_0,\ \ \ \pi=\pi_0+\A_M(\Phi_0^*\varepsilon)\]
where $\varepsilon\in  C^\infty(\g_\natural,\wedge^2\g)$ is the
solution of the CDYBE defined in
  Section \ref{subsec:exp1}.  Then $(M,\A_M,\pi,\Phi)$ is a q-Hamiltonian q-Poisson
  $\g$-manifold.
\end{proposition}
\begin{proof}
It is well-known that $(M,\A_M,\pi_0,\Phi_0)$ is a Hamiltonian
$\g$-manifold if and only if $\Phi_0\colon M\to \g^*$ is a Poisson
map, i.e., if and only if
\[(\Phi_0,0)\colon (M,E_M,0)\to (\g^*,E_{\pi_{\g^*}},0)\]
is a strong Dirac morphism, with $E_M=\on{Gr}_{\pi_0}$ and
$E_{\g^*}=\on{Gr}_{\pi_{\g^*}}$. Using $B$ to identify
$\g^*\cong\g$, and composing with the strong Dirac morphism
$(\exp,\varpi)$, one obtains the strong Dirac morphism
\[
(\Phi,\Phi_0^*\varpi)\colon (M,E_M,0)\to (G,E_G,\eta),
\]
which in turn gives rise to a q-Hamiltonian q-Poisson
$\g$-manifold $(M,\A_M,\pi,\Phi)$. The backward image
$\wt{F}_M\subset \TM$ of $F_G$ under the Dirac morphism
$(\Phi,\Phi_0^*\varpi)$ is a Lagrangian complement to
$E_M=\on{Gr}_\pi$.  Let $\wh{\a}\colon \Phi_0^*E_\g\to E_M$ be defined
by the Dirac morphism $(\Phi_0,0)$, and put $ \wh{\A}_M(\xi)=
\wh{\a}\circ \Phi_0^*\sfe_0(\xi)$. As explained in Section
\ref{subsec:exp1}, $\wt{F}_M$ is related the Lagrangian complement
$F_M=TM$ by the section $\wh{\A}_M(\Phi_0^* \varepsilon)$. Hence,
$\pi=\pi_0+\A_M(\Phi_0^* \varepsilon)$.
\end{proof}


\subsection{$\k^*$-valued moment maps}
Let $K$ be any Lie group.  An ordinary Hamiltonian Poisson
$K$-manifold is a triple $(M,\pi,\Phi)$ where $M$ is a
$K$-manifold, $\pi\in\mf{X}^2(M)$ is an invariant Poisson
structure, and $\Phi\colon M\to \k^*$ is a $K$-equivariant map
satisfying the moment map condition,
\[ \pi^\sharp(\d\l\Phi,\xi\r)=\A_M(\xi).\]
The moment map condition is equivalent to $\Phi$ being a Poisson
map. The following result implies that $\k^*$-valued moment maps
can be viewed as special cases of $G=\k^*\rtimes K$-valued moment
maps. Let $\g=\k^*\rtimes\k$ carry the invariant inner product
given by the pairing.

\begin{proposition}
  The inclusion map $j\colon \k^*\hra \k^*\rtimes K=G$ is a strong
  Dirac morphism $(j,0)$, as well as a backward Dirac morphism,
   relative to the Kirillov-Poisson structure on $\k^*$ and the
  Cartan-Dirac structure on $G$. The backward image of $F_G$ under
  this Dirac morphism is $F_{\k^*}=T\k^*$.
The pure spinor $\psi_G$ on $G=\k^*\rtimes K$ satisfies
\[j^*\psi_G=1.\]
\end{proposition}
\begin{proof}
  The Cartan-Dirac structure $E_G$ is spanned by the sections
  ${\sf{e}}(w)$ for $w=(\beta,\xi)\in\g$, while $E_{\k^*}$ is spanned
  by the sections ${\sf{e}_0}(\xi)$ for $\xi\in\k$. The first part of
  the Proposition will follow once we show that
  ${\sf{s}}_0(\beta,\xi)\sim_{(j,0)} {\sf{s}}(\beta,\xi)$, i.e.
\begin{equation}\label{eq:easyenough}
 {\sf{e}}_0(\xi)\sim_{(j,0)} {\sf{e}}(\beta,\xi),\ \
 {\sf{f}}_0(\beta)\sim_{(j,0)} {\sf{f}}(\beta,\xi)
.\end{equation}
The vector field part of the first relation follows since the
inclusion $j\colon \k^*\hra \k^*\rtimes K$ is equivariant for the
conjugation action of $G=\k^*\rtimes K$. (Here, the
$\k^*$-component of $G$ acts trivially on $\k^*$, while the
$K$-component acts by the co-adjoint action.) For the 1-form part,
we note that the pull-back of the Maurer-Cartan forms
$\theta^L,\theta^R\in\Om^1(G)\otimes\g$ to the subgroup
$\k^*\subset G$ is the Maurer-Cartan form for additive group
$\k^*$, i.e.
\[ j^*\theta^L=j^*\theta^R=\theta_0\]
where the `tautological 1-form' $\theta_0\in
\Om^1(\k^*)\otimes\k^*$ is defined as in
Section~\ref{subsec:exp1}. Thus
\[ j^* B\big({\ts \f{\theta^L_G+\theta^R_G}{2}},(\beta,\xi)\big)
=B(\theta_0,(\beta,\xi))=\l\theta_0,\xi\r.\]
This verifies the first relation in \eqref{eq:easyenough}; the
second one is checked similarly.

Since the adjoint action $\Ad\colon G\to \on{O}(\g)$ is trivial
over $\k^*$, the lift $\tau\colon G\to \on{Pin}(\g)\subset
\Cl(\g)$ satisfies $\tau|_{\k^*}=1$. It follows that the pure
spinor $\psi_G=\ca{R}(q(\mu))$ satisfies $j^*\psi_G=1$.
\end{proof}

\begin{corollary}
  Let $(M,\pi)$ be a Poisson manifold. Then $\Phi\colon M\to \k^*$ is
  a Poisson map if and only if the composition $j\circ \Phi\colon M\to
  G$ is a strong Dirac morphism
\[(j\circ \Phi,0) \colon (M,\on{Gr}_\pi,0)\to (G,E_G,\eta).\]
\end{corollary}

Put differently, Hamiltonian Poisson $K$-manifolds are
q-Hamiltonian q-Poisson $\k^*\rtimes K$-manifolds for which the
moment map happens to take values in $\k^*$.

As a special case, a Hamiltonian $K$-manifold $(M,\om,\Phi)$ (with
$\om$ a symplectic 2-form, and $\Phi$ satisfying the moment map
condition $\iota(\A_M(\xi))\om=\d\l\Phi,\xi\r$) is equivalent to a
q-Hamiltonian $G=\k^*\rtimes K$-space for which the moment map
takes values in $\k^*$. Since $j^*\psi_G=1$, its q-Hamiltonian
volume form coincides with the usual Liouville form
$(\exp\om)^{[\on{top}]}$.


\section{$K^*$-valued moment maps}\label{sec:kstar}
For a Poisson Lie group $K$, J.-H. Lu \cite{lu:mo} introduced
another type of group-valued moment map, taking values in the dual
Poisson Lie group $K^*$. For a compact Lie group $K$, with its
standard Poisson structure, this moment map theory turns out to be
equivalent to the usual $\k^*$-valued one. In this Section, we
will re-examine this equivalence using the techniques developed in
this paper.

\subsection{Review of $K^*$-valued moment maps}
The theory of Poisson-Lie groups were introduced by Drinfeld in
\cite{dr:qu}, see e.g. \cite{ch:gu}
for an overview and bibliography. Suppose $K$ is a
connected Poisson Lie group, with Poisson structure defined by a
Manin triple $(\g,\k,\k')$. (That is, $\g$ is a Lie algebra with
an invariant split inner product, and $\k,\k'$ are complementary
Lagrangian subalgebras.) Use the paring to identify $\k'=\k^*$,
and let $K^*$ be the associated dual Poisson Lie group. We assume
that $\g$ integrates to a Lie group $G$ (the double) such that
$K,K^*$ are subgroups and the product map $K\times K^*\to G$ is a
diffeomorphism.  The left action of $K$ on $G$ descends to a
\emph{dressing action} $\A_{K^*}$ on $K^*$ (viewed as a
homogeneous space $G/K$).  The Poisson structure on $K^*$, or
equivalently its graph $E_{K^*}=\on{Gr}_{\pi_{K^*}}\subset \T
K^*$, may be expressed in terms of the infinitesimal dressing
action, as the span of sections
 \[ e_{K^*}(\xi)=\A_{K^*}(\xi)\oplus \l\theta^R_{K^*},\xi\r\]
 for $\xi\in\k$. Here $\theta^R_{K^*}\in\Om^1(K^*)\otimes\k^*$ is the
 right-invariant Maurer-Cartan form for $K^*$. Note that as a Lie
 algebroid, $E_{K^*}$ is just the action algebroid.

 For the remainder of this Section \ref{sec:kstar}, we will assume
 that $K$ is a compact real Lie group.  The \emph{standard Poisson
   structure} on $K$ is described as follows. Let $G=K^\C$ be the
 complexification, with Lie algebra $\g$, and let
\[ \g=\k\oplus\a\oplus \n,\ \ \ G=KAN\]
be the Iwasawa decompositions. Here $\mf{a}=\sqrt{-1}\t_K,\
A=\exp\mf{a}$ and $N=N_+$ (using the notation from Section
\ref{subsec:gauss}).  We denote by $B_K$ an invariant inner
product on $\k$, and let $\l\cdot,\cdot\r$ be the imaginary part
of $2 B_K^\C$. Then $(\g,\k,\a\oplus \n)$ (where $\g$ is viewed as
a real Lie algebra) is a Manin triple. Thus $K$ becomes a Poisson
Lie group, with dual Poisson Lie group $K^*=AN$.

A \emph{$K^*$-valued Hamiltonian $\k$-manifold}, as defined by Lu
\cite{lu:mo}, is a symplectic manifold $(M,\om)$ together with a
Poisson map $\Phi\colon M\to K^*$.  Equivalently,
$(\Phi,\om)\colon (M,TM,0)\to (K^*,E_{K^*},0)$ is a strong Dirac
morphism. The Poisson map $\Phi$ induces a $\k$-action on $M$, and
if this action integrates to an action of $K$ we speak of a
\emph{$K^*$-valued Hamiltonian
  $K$-manifold}.  An interesting feature is that $\om$ is not
$K$-invariant, in general: Instead, the action map $K\times M\to M$ is
a Poisson map.  Accordingly, the volume form $(\exp\om)^{[\on{top}]}$
is not $K$-invariant. However, let $\Phi^A\colon M\to A$ be the
composition of $\Phi$ with projection $K^*=AN\to A$, and
$(\Phi^A)^{2\rho}\colon M\to \R_{>0}$ its image under the homomorphism
$T\to \C^\times,\ t\mapsto t^{2\rho}$ defined by the sum of positive
roots. By \cite[Theorem 5.1]{al:li}, the product
\begin{equation} \label{eq:invtvol}(\Phi^A)^{2\rho}\ (\exp\om)^{[\on{top}]}
\end{equation}
is a $K$-invariant volume form. The proof in \cite{al:li} uses a
tricky argument; one of the goals of this Section is to give a
more conceptual explanation.

\subsection{$P$-valued moment maps.}
To explain the origin of the volume form \eqref{eq:invtvol}, we
will use the notion of a $P$-valued moment map introduced in
\cite[Section 10]{al:mom}. Let $g\mapsto g^c$ denote the complex
conjugation map on $G$, and let
\[ I(g)\equiv g^\dagger =(g^{-1})^c.\]
On the Lie algebra level, let $\xi\mapsto \xi^c$ denote
conjugation, and $\xi^\dagger=-\xi^c$. We have $K=\{g\in G|\
g^\dagger=g^{-1}\}$. Let
\[ P=\{g^\dagger g|\ g\in G\} \]
denote the subset of `positive definite' elements in $G$. Then $P$
is a submanifold fixed under $I$, and the product map defines the
Cartan decomposition $G=KP$. Let $E_G$ be the (holomorphic) Dirac
structure on $G$ defined by the inner product
\[ B:={\ts \f{1}{\sqrt{-1}}}B_K^\C.\]
Since $(\theta^L)^\dagger=I^*\theta^R,\ \
(\theta^R)^\dagger=I^*\theta^L$, the Cartan 3-form on $G$
satisfies satisfies $\eta^c=I^*\eta$, thus $\eta_P:=\iota_P^*\eta$
is real-valued. Similarly, the pull-backs of the 1-forms $B({\ts
  \f{\theta^L+\theta^R}{2}},\xi)$ for $\xi\in\k$ are real-valued. It
follows that the sections
\[ {\sf{e}}_P(\xi):={\sf{e}}(\xi)|_P\]
are real-valued. Letting $E_P\subset \T P$ be the subbundle
spanned by these sections, it follows that $(P,E_P,\eta_P)$ is a
real Dirac manifold, with $(E_P)^\C=E_G|_P$. As a Lie algebroid,
$E_P$ is just the action algebroid for the $K$-action on $P$.
Similarly, the sections ${\sf{f}}_P(\xi):={\sf{f}}(\xi)|_P$ are
real-valued, defining a complement $F_P$ to $E_P$. The bundle
$F_P$ is defined by the (real-valued) pure spinor,
$\psi_P:=\iota_P^*\psi_G\in\Om(P)$.
\begin{remark}
  Since $\det(\Ad_g+1)>0$ for $g\in P$ (all
  eigenvalues of $\Ad_g$ are strictly positive), one finds that
  $\ker(E_P)=\{0\}$. Hence $E_P$ is the graph of a bivector $\pi_P$
  with $\hh [\pi_P,\pi_P]=\pi_P^\sharp(\eta_P)$.
\end{remark}
A \emph{$P$-valued Hamiltonian $\k$-manifold} \cite[Section
10]{al:mom} is a manifold $M$ together with a strong Dirac
morphism $(\Phi_1,\om_1)\colon (M,TM,0)\to (P,E_P,\eta_P)$.  For
any such space we obtain, as for the q-Hamiltonian setting, an
invariant volume form
\begin{equation}\label{eq:Pvolume}
 (\exp(\om_1)\wedge \Phi_1^*\psi_P)^{[\on{top}]}.
\end{equation}
Here $\psi_P$ may be replaced by $\wh{\psi}_P$, the pull-back of
the Gauss-Dirac spinor.\footnote{In Section \ref{subsec:gaussvol},
$B$ was
  taken as the complexification of $B_K$, while here we have an extra
  factor $\sqrt{-1}$. This amounts to a simple rescaling of the
  bilinear form $B_K^\C$, not affecting any of the results.}  By
Proposition \ref{prop:gaussvol},   the expression $\exp\om_1
\wedge \Phi_1^*(\Delta_\lambda\ \wh{\psi}_P)$ is closed under the
differential $\d-2\pi \iota(\A_M(\lambda+\rho)))$, for
any dominant weight $\lambda$.

\subsection{Equivalence between $K^*$-valued and $P$-valued moment maps}
To relate the $K^*$-valued theory with the $P$-valued theory, we
use the $K$-equivariant diffeomorphism
\[ \kappa\colon K^*\to P,\ g\mapsto g^\dagger g.\]
Note that this map takes values in the big Gauss cell,
$\O=N_-KN\subset G$. Let $\varpi_\O$ denote the (complex) 2-form
on the big Gauss cell, and $\varpi_{K^*}=\kappa^*\om_\O$. It is
easy to check that $\varpi_{K^*}$ is real-valued. One can check
that
\[ e_{K^*}(\xi)\sim_{(\kappa,\varpi_{K^*})} e_P(\xi)\]
for all $\xi\in\k$: The vector field part of this relation is
equivalent to the $\k$-equivariance, while the 1-form part is
verified in \cite[Section 10]{al:mom}. It follows that
$(\kappa,\varpi_{K^*})$ is a Dirac isomorphism from
$(K^*,E_{K^*},0)$ onto $(P,E_P,\eta_P)$.

  Thus, if $(M,\om,\Phi)$ is a $K^*$-valued Hamiltonian $\k$-manifold,
  then $(M,\om_1,\Phi_1)$ with   $\om_1=\om+\Phi^*\varpi_{K^*}$ and $\Phi_1=\kappa\circ\Phi$ is a $P$-valued
  Hamiltonian $\k$-manifold. In particular, we obtain an invariant
  volume form on $M$,
\[ \big(\exp(\om+\Phi^*\varpi_{K^*})\wedge
\Phi^*\kappa^*\wh{\psi}_P\big)^{[\on{top}]}.
\]
Using the explicit formula (Proposition \ref{prop:gaussexplicit})
for the Gauss-Dirac spinor, we obtain
\[ \kappa^*\wh{\psi}_P=a^{2\rho}\exp(-\varpi_{K^*}),\]
where $a\colon K^*\to A$ is projection to the $A$-factor. Hence,
\[\exp(\om+\Phi^*\varpi_{K^*})\wedge
\Phi^*\kappa^*\psi_P= (\Phi^A)^{2\rho} \exp(\om),\]
identifying the volume form for the associated $P$-valued space
with the volume form \eqref{eq:invtvol}.
\begin{proposition}
  For any $K^*$-valued Hamiltonian $\k$-space $(M,\om,\Phi)$, the
  volume form $(\Phi^A)^{2\rho}(\exp\om)^{[\on{top}]}$ is $\k$-invariant.
  Moreover, for all dominant weights $\lambda$ the
  differential form
\[ (\Phi^A)^{2(\lambda+\rho)} \exp(\om)\]
is closed under the differential $\d-2\pi
\A_M(B_K^\sharp(\lambda+\rho))$.
\end{proposition}
\begin{proof}
  Invariance follows from the identification with the volume form for
  the associated $P$-valued space. The second claim follows from
  Proposition \ref{prop:gaussvol}, since the function
  $\Delta_\lambda$ from Section \ref{subsec:gdspinor} satisfies
  $\kappa^*\Delta_\lambda=a^{2\lambda}$.
\end{proof}
The differential equation permits a computation of the integrals
$\int_M (\Phi^A)^{2(\lambda+\rho)} (\exp(\om))^{[\on{top}]}$ by
localization \cite{be:ze} to the zeroes of the vector field
$\A_M(B_K^\sharp(\lambda+\rho))$, similar to the formula in
\ref{subsec:gaussvol}.

\subsection{Equivalence between $P$-valued and $\k^*$-valued moment maps}
Finally, let us express the correspondence \cite[Section
10]{al:mom} between $P$-valued moment maps and $\k^*$-valued
moment maps in terms of Dirac morphisms. The exponential map for
$G=K^\C$ restricts to a diffeomorphism
\[ \exp_{\mf{p}}\colon \mf{p}:=\sqrt{-1}\k\to P:=\exp(\sqrt{-1}\k).\]
Let $\varpi\in\Om^2(\g)$ be the primitive of $\exp^*\eta$ defined
in \eqref{eq:varpi}, and $\varpi_{\mf{p}}$ its pull-back to
$\mf{p}$. Since $\eta_P$ is real-valued, so is $\varpi_{\mf{p}}$,
and
  $\d\varpi_{\mf{p}}=(\exp|_{\mf{p}})^*\eta_P$. Similarly,
  $J_{\mf{p}}:=J|_{\mf{p}}>0$. The formulas for
  $\varpi_{\mf{p}}$ and $J_{\mf{p}}$ are similar to those for the Lie
  algebra $\k$, but with $\sinh$ functions replaced by $\sin$
  functions. Use $B^\sharp=\sqrt{-1}B_K^\sharp$ to identify $\k^*\cong
  \mf{p}$. By Proposition \ref{prop:exp},
\[
{\sf{e}}_0(\xi)\sim_{(\exp_{\mf{p}},\varpi_{\mf{p}})}{\sf{e}}_P(\xi),\
\ \ \xi\in\k.\]
Hence $(\exp_{\mf{p}},\varpi_{\mf{p}})$ is a Dirac (iso)morphism
from $(\k^*,E_{\k^*},0)$ to $(P,E_P,\eta_P)$. This sets up a 1-1
correspondence between $P$-valued and $\k^*$-valued Hamiltonian
$\k$-spaces. Thinking of the latter as given by strong Dirac
morphisms $(\Phi_0,\om_0)$ to $(\k^*,E_{\k^*},0)$, the
correspondence reads
\[ (\Phi_1,\om_1)=(\exp_{\mf{p}},\varpi_{\mf{p}})\circ (\Phi_0,\om_0).\]
The volume forms are related by
$(\exp(\om_1)\wedge\Phi_1^*\psi_P)^{[\on{top}]}=J_{\mf{p}}^{1/2}\exp(\om_0)^{[\on{top}]}$.

\bibliographystyle{amsplain}

\def\cprime{$'$} \def\polhk#1{\setbox0=\hbox{#1}{\ooalign{\hidewidth
  \lower1.5ex\hbox{`}\hidewidth\crcr\unhbox0}}} \def\cprime{$'$}
  \def\cprime{$'$} \def\polhk#1{\setbox0=\hbox{#1}{\ooalign{\hidewidth
  \lower1.5ex\hbox{`}\hidewidth\crcr\unhbox0}}} \def\cprime{$'$}
  \def\cprime{$'$}
\providecommand{\bysame}{\leavevmode\hbox to3em{\hrulefill}\thinspace}
\providecommand{\MR}{\relax\ifhmode\unskip\space\fi MR }
\providecommand{\MRhref}[2]{%
  \href{http://www.ams.org/mathscinet-getitem?mr=#1}{#2}
}
\providecommand{\href}[2]{#2}

\end{document}